\documentclass[13pt]{amsart}
    \usepackage{amsmath}
    \usepackage{}
    \usepackage{amsfonts}
    \usepackage{amsfonts,latexsym,rawfonts,amsmath,amssymb,amsthm}
    \usepackage[plainpages=false]{hyperref}

    \usepackage[a4paper]{geometry}
    \usepackage{esint}
    \usepackage{graphics,graphicx}
    \usepackage{tikz}

    \numberwithin{equation}{section}

    \newcommand{\beq}{\begin{equation}}
    \newcommand{\eeq}{\end{equation}}
    \newcommand{\beqs}{\begin{eqnarray*}}
    \newcommand{\eeqs}{\end{eqnarray*}}
    \newcommand{\beqn}{\begin{eqnarray}}
    \newcommand{\eeqn}{\end{eqnarray}}
    \newcommand{\beqa}{\begin{array}}
    \newcommand{\eeqa}{\end{array}}

    \def\lra{\longrightarrow}

    \def\bc{\begin{center}}
    \def\ec{\end{center}}

    \let\textPolishL\l
    \def\l{\lambda}

    \def\begeq{\begin{equation}}
    \def\endeq{\end{equation}}
    \def\and{\quad{\rm and}\quad}

    \let\lra=\longrightarrow

    \def\mapright\#1{\,\smash{\mathop{\lra}\limits^{\#1}}\,}

    \newtheorem{prop}{Proposition}[section]
    \newtheorem{theo}[prop]{Theorem}
    \newtheorem{lem}[prop]{Lemma}
    \newtheorem{claim}[prop]{Claim}
    \newtheorem{cor}[prop]{Corollary}
    \newtheorem{rem}[prop]{Remark}
    
    \newtheorem{defi}[prop]{Definition}
    \newtheorem{conj}[prop]{Conjecture}

\begin{document}

 	\title[Singular  K\"ahler-Ricci solitons ]{Singular  K\"ahler-Ricci solitons  and complete gradient shrinking   K\"ahler-Ricci solitons on   non-compact  manifolds}
 	  \author{Feng Wang  and Xiaohua  ${\rm ZHU}^{\dag}$}
 	
\address{School of Mathematical Sciences, Zhejiang
University, Hangzhou,  310028, China.}
\email{wfmath@zju.edu.cn}

 	\address{School of Mathematical Sciences, Capital Normal
 		University, Beijing 100048, China.}
 	\email{xhzhu@math.pku.edu.cn}
 	
 	\thanks {$\dag$Partially supported by National Key R\&D Programs of China 2023YFA1009900  and  NSFC 12271009.}
 	\begin{abstract}   We  construct a class of complete gradient  shrinking K\"ahler-Ricci solitons
	on non-compact  toric manifolds by exhausting a family of singular K\"ahler-Ricci solitons associated to  multiple divisors.  The construction  method will be also extended on 
	 non-compact  $\mathbf G$-manifolds. 

	\end{abstract}
 	\date{}
 	\maketitle
 	 \tableofcontents

	\section{Introduction}
	
	As a generalization of Einstein metric, Ricci soliton was extensively studied in variant  aspects. 
	One of motivations is  from  the study of  Hamilton's  Ricci flow \cite{Ha82}. 
	As a class of  important singular models,  gradient Ricci solitons play a  crucial  role in the singularity analysis  of Ricci flow (\cite{Ha95,  Pe02}).  For examples,  near any singular point of  tpye I,    Ricci flow
	will develop a    gradient  shrinking Ricci soliton (cf. \cite{Ni,  CLN, EMT}).   This  encourages  people to classify such kind of   gradient Ricci solitons.

	However,  we still understand   the Ricci soliton not  so much.  Actually,  almost examples  come from K\"ahler metrics with large symmetry  (cf. \cite{WZ04, SZ, FIK, FW, Del, LZZ,  LZhu,  LTZ}, etc.).   
	In K\"ahler geometry,  a famous Hamilton-Tian conjecture confirms that a K\"ahlrer-Ricci (abbreviated by KR) flow  on a Fano manifold  with an initial metric in the  canonical K\"ahler class must converge to a  KR soliton  (possibly with singularities) \cite{Ti97, Ti15, TZh, WZ20, TZZZZ}.  The conjecture  tells us that   the KR soliton is   general  canonical model  more than  the  K\"ahler-Einstein metric for Fano manifolds.
	
	 As we know,   different to expanding  Ricci soliton and steady   Ricci soliton,    non-compact    gradient  shrinking  Ricci solitons  have  more  strong rigidities  (cf. \cite{MW-positive, KW-rigidity}).   Up to know, we only know two non-trivial  classes of such  Ricci solitons (but not Einsteins) from K\"ahler geometry,  one is called FIK metrics constructed by ODE method in \cite{FIK}, while another is from the blow-up solution of  tpye I singularity in Ricci flow on   toric surface $\mathbb CP^2\# 2\overline {\mathbb CP^2}$ constructed  by Bamler--Cifarelli--Conlon--Deruelle in \cite{BCCD}. The latter is a  gradient  shrinking  KR soliton  (abbreviated by GSKRS)
	on  toric surface,  the blow-up space  $BL_p(\mathbb CP^1\times \mathbb C)$ of $\mathbb CP^1\times \mathbb C$ at one point.  Actually,   with help of   these two kind of  examples,  a complete classification of GSKRS on K\"ahler surfaces can be     finished  by combining  a series of works \cite{WZ04, CCD-classification, CCD-type1,  BCCD, LW}.
	
	The purpose of this paper, we  provide an approach  to produce   GSKRS  on a class of non-compact  K\"ahler  manifolds $M$ with large symmetric group
	via PDE method, in particular on  non-compact   toric manifolds.  The idea  is similar with  the Calabi-Yau metric case \cite{TY1, TY2}, we consider $M=\bar M\setminus D$ as an open manifold  of compact K\"ahler manifold $\bar M$, where $D$ is a  divisor of $\bar M$. However, because of lack of apriori estimates,  we are not able to  solve a global   GSKRS equation directly,  while replacing by  an    exhausting  method,  we  return to   solve a family of   singular  KR soliton equations for large cohomology K\"ahler class $\Omega_k$ in $\bar M$,  where integer  $k\gg  1$,  see Section 2 below. 
	
	In the case of  toric manifold $M=\bar M\setminus D$,  we  assume that $\bar M$ is a compact toric manifold while $D$ is a  simple nonnegative  multiple     toric divisor.  Then the toric structure of $M$  is determined by  a   non-compact polytope   $P$.
	Inspired by a recent work of Cifarelli \cite{Ci}  on  the uniqueness of  GSKRS  on
	non-compact  toric manifolds,  we  assume 
	$P$ is a 
	a positive   Fano cone ( see Definition \ref{cone-polytope}), which means  
	that there exists a   vector $A_0\in \mathbb R^n$ and a positive constant $c_0$ such that for any $y\in P$ with  $|y|\gg 1$ it holds
	\begin{align}\label{behavior-P}
		<A_0, y>\le -c_0. 
	\end{align}
	The condition (\ref{behavior-P}) guarantees  that $H_0$-invariant (or log $H$-invariant \cite{TZZZ}) introduced in \cite{TZ02}  is well-defined for  the case of 
	non-compact  toric manifolds.   As a consequence,  one can uniquely determine a soliton vector field so that  associated  modified Futaki-invariant vanishes  as in the case of  toric Fano manifolds \cite{WZ04}.
	
	Set a subset  in $\mathbb R^{n-1}$ by
	\begin{align}\label{project-p}P'=\{x'\in \mathbb R^{n-1}|~(x', x^n)\in P\}.
	\end{align}
	It is easy to see that $P$ satisfies   (\ref{behavior-P})  when $P'$  is a bounded subset in $\mathbb R^{n-1}$. For simplicity, we call $P$ parabolic  in this  special case while $P$ is called 
	positive if  $(\ref{behavior-P})$ is satisfied.  	Clearly,   the  associated  $2d$ Fano cone to Bamler--Cifarelli--Conlon--Deruelle's example $BL_p(\mathbb CP^1\times \mathbb C)$ satisfies  the parabolic  condition.  In   higher dimensional case,   
	the    associated   Fano cone to  the blow-up space of $ \mathbb CP^{n-1}\times \mathbb C$  at one  point  is  an example of  parabolic type.   One can see  more examples  that   all associated   Fano cone   to  the blow-up space of $ N\times \mathbb C$     at one  point is    parabolic by assuming the  blow-up space  is toric, where $N$ is a compact  toric manifold. 
	
	The following is  our main result in this paper.
	
	\begin{theo}\label{main-resullt}
		Let $M^n$  be a non-compact  toric manifold associated to a  parabolic Fano  polytope  cone $P$. 
		Then there exists a GSKRS on  the torus orbit of $M$ with infinite volume, which can be extended to a complete GSKRS on   $M$.
	\end{theo}

	By Definition \ref{cone-polytope},   a   Fano cone polytope  $P$  implies the existence of  Delzant polytopes $P_k$ for any positive integers $k\ge k_0$, where $k_0>0$ is an integer and  $P_k$ is associated to a same compact  toric manifold $\bar M$ with various K\"ahler class $\Omega_k$. We will first prove the existence of singular KR solitons $\omega_{\varphi_k}$   on $(\bar M, \Omega_k)$ associated to the multiple  divisor  $kD$  in sense of  Definition \ref{kr-soliton} ,  see  Theorem \ref{existence-1}  in Section 5.  We note that the corresponding equation  of singular KR soliton  will  reduce a certain  real Monge-Amp\`ere equation in our case.  
	Thus we can use the variation method  to solve the existence    as in  \cite{BB, LZZ,  LZhu,  LTZ} for Fano varieties with large torus actions.   A  regularity result for such  singular KR solitons will be also proved by  using  a  fundamental work of   Berman, et al.  in \cite{BBEGZ},  see Theorem \ref{partial-smooth}.
	
	Next,  we want to show that $\omega_{\varphi_k}$ converges to a GSKRS  $\omega_\infty$ on  the torus orbit of $M$ with associated moment  polytope  $P$.  The proof  will depend on various  a priori   estimates for the convex solutions  $\varphi_k$ of   real Monge-Amp\`ere equations,  see Section 6. 
	
	At last,  we shall prove that the   GSKRS  $\omega_\infty$ can be  extended to a complete GSKRS on   $M$. Our method is to
	get a local  uniform $C^2$-estimate for  $\omega_{\varphi_k}$
	   with respect to a chosen  backgroud K\"ahler metric for any open manifold $\bar M\setminus  U_{D,\epsilon}$, where $  U_{D,\epsilon}$ is a small $\epsilon$-neighborhood of $D$.    
	 The argument for the regularity  proofs in  Section 2 and Section 6 will be extended.

	We would like to emphasize that   we shall assume that  $P$ is parabolic in last two steps. We hope this assumption can be removed. Namely, we conjecture that Theorem \ref{main-resullt} is true for any positive Fano cone. Actually, the same result can be generalized to $\mathbf G$-manifolds, see Conjecture \ref{conj-g}.   The conjecture  will be also verified  for  $\mathbf G$-manifolds  with   parabolic Fano cones as a generalization of Theorem \ref{main-resullt}, see Theorem \ref{g-para}.
	
	The paper is organized as follows.  
In Section 2, we introduce singular   KR solitons   on a compact  K\"ahler manifold associated to a divisor and prove a regularity theorem for those singular metrics, see Theorem \ref{partial-smooth}. In Section 3, we introduce  a  class  of non-compact toric manifolds associated to Fano polytope cone and prove the existence of soliton vector fields while in  Section 4 we construct 
a class of admissible  K\"ahler metrics by modifying classical Gulliemin functions. In Section 5, we prove the existence of singular   KR solitons via the variation method, see Theorem \ref{existence-1}.
In Section 6, we give various estimates in order to prove the convergence of singular   KR solitons  constructed in Theorem \ref{existence-1}, see Theorem \ref{limit-solution}.  Theorem~\ref{main-resullt}will be proved in Section 7 by establishing  the global  regularity for the limit solution in Theorem \ref{limit-solution}. In Section 8, we generalize Theorem~\ref{main-resullt} to $\mathbf G$-manifolds, see Theorem \ref{g-para}. The last Section 9 is an appendix, where we construct  several  examples of  GSKRS on non-compact toric manifolds  and   $\mathbf G$-manifolds.

	\vskip5mm
	
	$\mathbf {Note}$:  The work on the first part of Theorem~\ref{main-resullt} was listed as in  preparation in a  survey paper \cite{TZ-applications}. The main result was also talked by the second author  in  School and Workshop on Geometric Analysis, ICTP, 24 August-4 September, 2026.
	We also note that
	Collins-Huang-Tong-Xu recently posted an overlap  result as in  Theorem \ref{main-resullt} via AI \cite{CHTX}. Our approach and  method in the present  paper are both  different to theirs.

\section{Singular  K\"ahler-Ricci soliton}

Let $(\bar M, \omega_g)$  be a  compact K\"ahler manifold with
a K\"ahler form
$$\omega_g=\sqrt{-1} g_{i,\bar j} dz^i\wedge d\bar z^j.$$
Suppose that  $D$ is  a divisor on $\bar M$ such  that
\begin{align}\label{D-condition}  [\omega_g]-2\pi c_1([D])\in 2\pi c_1(\bar M) .
\end{align}
In this paper, we always  assume that $D$ is an effective  semi-positive  divisor  with a holomorphic defining section on $\bar M$.

\begin{defi}\label{admissable-0} Let  $D$ be a divisor  which satisfies  (\ref{D-condition}).  We call $\omega\in [\omega_g]$  with  $L^\infty$ K\"ahler potential  
 an admissible  K\"ahler metric     
  associated to  $D$ if there is a bounded  function $h$ on $\bar M$ such that
\begin{align}\label{admi-equ}{\rm Ric}(\omega)-\omega+ [D]=\sqrt{-1}\partial\bar\partial h.
\end{align}
$h$ is called a twisted Ricci potential associated to $D$ in sense of current.

\end{defi}

\begin{defi}\label{kr-soliton} Let $D$ as in Definition \ref{admissable-0}.  Let
$$\omega=\omega_{\phi}=\omega_g+\sqrt{-1}\partial\bar\partial\phi.$$
 We call $\omega\in [\omega_g]$  a singular  K\"ahler-Ricci (KR) soliton  with a  K\"ahler potential
$\phi\in L^\infty(\bar M)$
 associated to $D$  if there exists a  holomorphic vector field  (HVF) $X$ on $\bar M$   with 
  $X(\phi)\in L^\infty(\bar M)$ such that  $\omega$ satisfies
\begin{align}\label{twisted-equ-0}{\rm Ric}(\omega)-\omega+[D]= L_X\omega.
\end{align}
\end{defi}

In literature,  $\omega$ in (\ref{twisted-equ-0}) is also called a KR soliton  associated to  a Log Fano pair $(\bar M, -D)$, see \cite{HL,  BLXZ}.  This class of  metrics include  conical KR solitons with big conical angle in the case of Fano manifold $\bar M$  when $\omega\in \alpha c_1(M)$ and $c_1(D)\in (\alpha-1) c_1(M)$,  where $\alpha\ge 1$. To see   examples  of conical KR solitons  with cone angles  less  than $2\pi$,   we refer the reader to  \cite{DGSW, DeS20, WZZ}, etc..

\subsection{Complex Monge-Amp\`ere equations  associated to (\ref{twisted-equ-0})}

Suppose  that $(\bar M, \omega_g)$  admits an  admissible  K\"ahler metric $\omega\in [\omega_g]$.  Then    (\ref{twisted-equ-0}) for a singular  KR soliton  $\omega=\omega_{\phi}=\omega+\sqrt{-1}\partial\bar\partial\phi$ is equivalent to the following complex Monge-Amp\`ere (CMA) equation,
\begin{align}\label{MAE}\omega_{\phi}^n= e^{h-\theta_X-X(\phi)-\phi}
 \omega^n,
\end{align}
where $\theta_X$ is a potential of $X$ associated to   the background metric $\omega$ such that
\begin{align}\label{potential-X} L_X\omega= \sqrt{-1}\partial\bar\partial \theta_X.
\end{align}
In geometry, (\ref{MAE}) is
equivalent to
\begin{align}\label{ricci-potential}h_\phi=\theta_X(\phi)=\theta_X+X(\phi),\end{align}
where
$h_\phi$ is the twisted Ricci potential of $\omega_\phi$
 Actually, one can check that for any $\omega_\phi$ it holds
$$h_\phi=\log [\frac{\omega_{\phi}^{n}}{\omega_{g}^{n}}e^{\phi- h_{g}}].$$

On the other hand, by (\ref{admi-equ}), we can choose  a hermitian metric   $H$ on $D$ such that for a smooth K\"ahler metric $\omega_0\in  [\omega_g]$ it holds
\begin{align}\label{h-potential}{\rm Ric}(\omega_0)-\omega_0+c_1(L_D, H)= \sqrt{-1}\partial\bar\partial h_0,
\end{align}
where $h_0$ is a smooth function on $M$ as a twisted Ricci potential of $\omega_0$.  It follows that there is a  $L^\infty$-function $\hat h$ such that
$$\omega_g^n= e^{\hat h}|S|_H^2\omega_0^n,$$
where $S$ is a defining function of $D$. Thus   (\ref{MAE}) can rewritten  as
\begin{align}\label{MAE-deg}\omega_{\phi}^n= e^{h_0-\theta_X-X(\phi)-\phi}|S|_H^2\omega_0^n,
\end{align}
where $\omega_{\phi}=\omega_0+\sqrt{-1}\partial\bar\partial\phi$ and K\"ahler potential $\phi $    is  different to one in (\ref{MAE}) with  a $L^\infty$-function,   and  $\theta_X$ is a potential of $X$ associated to  $\omega_0$ as in (\ref{potential-X}).

\subsection{Regularity of solutions on  (\ref{MAE-deg})}
By studying (\ref{MAE-deg}),  we have the following regularity result for the singular KR soliton.

\begin{theo}\label{partial-smooth} Any singular KR soliton   $\omega_\varphi$ of (\ref{MAE}) in the sense of Definition \ref{kr-soliton} is smooth on $\bar M\setminus D$.
\end{theo}

To prove Theorem \ref{partial-smooth},  we need the following variant of Yau's  $C^2$-estimate  using an idea in \cite{BBEGZ}.

\begin{lem}\label{C^2}
Let $dV$ be a smooth volume volume on $\bar M$ and $X$ be a holomorphic vector field. Let $\psi^{+}$ and $\psi^{-}$ be   two smooth functions of $\bar M$ which satisfy:
\begin{align}
i)~ &\sqrt{-1}\partial \bar\partial\psi^{\pm} \geq -C \omega_0\notag;\\
ii) ~&|\psi^{-}|+|X(\psi^{-})| \leq C;\notag\\
iii)~& e^{\psi^{+}} \leq C,\notag
\end{align}
 where   $C>0$ is a uniform constant.  Assume that $\phi$ is a smooth $\omega_0$-psh function such that $$\left(\omega_0+\sqrt{-1}\partial \bar\partial \phi\right)^n=e^{-X(\phi)+\psi^{+}-\psi^{-}}dV.$$
Then there exists $A>0$ only depending on $\omega_0, X, dV$ and $C$ such that
\begin{align}\label{yau-estimate}
0 \leq \omega_0+\sqrt{-1}\partial \bar\partial \phi \leq A \omega_0.
\end{align}
\end{lem}

\begin{proof}
We can normalize $\phi$ such that $\sup\phi=0$. By $L^\infty$ estimate, we have $|\phi|\leq A$. Denote $\omega_\phi=\omega_0+\sqrt{-1}\partial \bar\partial \phi$.
At first, we have (cf. \cite{Yau} and \cite[pp.~98--99]{Siu}),

\begin{align}\label{Yau's C^2}
\Delta_{\omega_{\phi}} \log \operatorname{tr}_{\omega_0} \omega_{\phi} \geq -\frac{\operatorname{tr}_{\omega_0} \operatorname{Ric}\omega_{\phi}}{\operatorname{tr}_{\omega_0}(\omega_{\phi})}-B \operatorname{tr}_{\omega_{\phi}}(\omega_{0}),
\end{align}
 where $B$ depends on the lower bound of the bisectional curvature of $\omega_0$. From the equation, we have
\begin{align}
-\operatorname{Ric}\omega_{\phi}&=-{\rm Ric}\,dV+\sqrt{-1}\partial \bar\partial\psi^+-\sqrt{-1}\partial \bar\partial\psi^- -\sqrt{-1}\partial \bar\partial X(\phi)\\ \notag
&\geq -A\omega_0-\sqrt{-1}\partial \bar\partial X(\phi)-\sqrt{-1}\partial \bar\partial\psi^-
\end{align}
for some constant $A>0$. It follows that
\begin{align}\label{low}
\Delta_{\omega_{\phi}} \log \operatorname{tr}_{\omega_0} \omega_{\phi} \geq -A{tr}_{\omega_\phi} \omega_{0}-\frac{\operatorname{tr}_{\omega_0}\sqrt{-1}\partial \bar\partial \psi^-}{{tr}_{\omega_0} \omega_{\phi}}-\frac{\operatorname{tr}_{\omega_0}\sqrt{-1}\partial \bar\partial X(\phi)}{{tr}_{\omega_0} \omega_{\phi}}.
\end{align}
Since $C \omega_{0}+\sqrt{-1}\partial \bar\partial\psi^-\geq 0$,  i.e., 
 $$
0 \leq C \omega_{0}+\sqrt{-1}\partial \bar\partial\psi^- \leq \operatorname{tr}_{\omega_{\phi}}\left(C \omega_{0}+\sqrt{-1}\partial \bar\partial\psi^-\right)\omega_{\phi}, 
$$
we have
$$
0 \leq C n+\Delta_{\omega_{0}} \psi^{-} \leq\left(C \operatorname{tr}_{\omega_{\phi}}\left(\omega_{0}\right)+\Delta_{\omega_{\phi}} \psi^{-}\right) \operatorname{tr}_{\omega_{0}}\left(\omega_{\phi}\right) .
$$
Hence,  by (\ref{low}),  we get
\begin{align}\label{subsolution}\Delta_{\omega_{\phi}} (\log \operatorname{tr}_{\omega_0} \omega_{\phi}+\psi^--2A\phi)\geq A{\rm tr}_{\omega_\phi} \omega_{0}-2An- \frac{\Delta_{\omega_0}X(\phi)}{{\rm tr}_{\omega_0} \omega_{\phi}}.
\end{align}

At the point $x_0$ where $G(\cdot)=\log \operatorname{tr}_{\omega_0} \omega_{\phi}+\psi^--2A\phi$ attains the maximum, we have $$\phi_{i\bar i j}(x_0)=(n+\Delta_{\omega_0}\phi)(2A\phi_j-\psi^-_j)(x_0).$$ 
It follows that
$$ |X_j \phi_{i\bar i j}|\le (n+\Delta_{\omega_0}\phi)|(2A|X(\phi)|+|X(\psi^-)|).$$
Since $\Delta_{\omega_0} X(\phi)=X_{ij}\phi_{\bar i j}+ X_j \phi_{i\bar i j}$ in normal coordinates, we get
\begin{align}\label{laplace-x}
|\Delta_{\omega_0} X(\phi)(x_0)|\le 
 (n+\Delta_{\omega_0}\phi)(2A|X(\phi)|+|X(\psi^-)|+C_1)+C_2.
\end{align}
 Note that  $|X(\phi)|\leq C$ by \cite[Corollary~5.3]{Zhu}. Thus we obtain 
 $$|\Delta_{\omega_0} X(\phi)(x_0)|\leq A'( \operatorname{tr}_{\omega_0} \omega_{\phi}(x_0)+1)$$
  for another constant $A'>0$.  
 Now by  (\ref{subsolution}) with the help of    the fact
  $$\operatorname{tr}_{\omega_0} \omega_{\phi}\leq  ({\rm tr}_{\omega_\phi} \omega_{0})^{n-1}\frac{\omega_\phi^n}{ \omega_{0}^n},$$  
   we can apply  the maximum principle to  the function $G(\cdot)$ to prove  the upper bound  in (\ref{yau-estimate}) (cf. \cite{Yau, TZ1}).
\end{proof}

\begin{proof}[Proof of Theorem \ref{partial-smooth}]
Consider the regularization $\psi^+_\delta:=\log (|S|_H^2+\delta)$ of $\log |s_D|_h^2$ and choose a sequence of smooth $\omega$-psh functions  $\tilde \phi_j$  which is
 $L^p$-convergent  to  $\phi$. 
   By  \cite{Zhu},  we can solve the following equation of $\phi_j$,
\begin{align}\label{Yau-solution} \left(\omega_0+\sqrt{-1}\partial \bar\partial \phi_j\right)^n=c_je^{h_0-\theta_X-X(\phi_j)+\psi_{\delta_j}^{+}-\tilde \phi_j}\omega_0^n,\end{align}
 where $c_j$ is chosen such that 
$\int_M c_je^{\psi_{\delta_j}^{+}-\tilde \phi_j}e^{h_0-\theta_X}\omega_0^n={\rm vol}(M)$ and $\delta_j=\frac{1}{j}$.
Note that $X(\phi_j)$ is uniformly bounded \cite{Zhu}.
 Thus, by  Kolodziej's result \cite{Kol98},  $ \phi_j$ is uniformly bounded.  Moreover,  we may also assume that   $X(\tilde \phi_j)$ is uniformly bounded,  otherwise we replace  $\tilde \phi_j$ by $\phi_j$ and take  the new solution of  (\ref{Yau-solution}) instead of $\phi_j$.

 On the other hand, it is easy to see that
\begin{align}
&\sqrt{-1}\partial\bar\partial\log(\delta+|S|_H^2)\notag\\
&=\sqrt{-1}\partial\frac{\bar\partial |S|_H^2}{\delta+|S|_H^2} \nonumber\\
&=\sqrt{-1}\partial\frac{S\cdot \overline{\nabla S}+|S|_H^2\bar\partial \log H}{\delta+|S|_H^2} \nonumber\\
&=\sqrt{-1}\frac{\nabla S\cdot \overline{\nabla S}+|S|_H^2\partial\bar\partial \log H}{\delta+|S|_0^2}-\sqrt{-1}\frac{(\nabla S\cdot \bar S)(S\cdot \overline{\nabla S})}{(\delta+|S|_H^2)^2} \nonumber\\
&=-\frac{|S|_H^2}{\delta+|S|_H^2} \omega_H+\delta\sqrt{-1}\frac{ \nabla S\cdot \overline{\nabla S}}{(\delta+|S|_H^2)^2}\notag\\
&\ge -C_0 \omega_0.\notag
\end{align}
 Then applying  Lemma \ref{C^2} with $\psi_{\delta_j}^+$ as $\psi^+$, $\tilde \phi_j$ as $\psi^-$ and $e^{h_0-\theta_X}\omega_0^n$ as the smooth volume form, we know that $\Delta_{\omega_0}\phi_j$ is uniformly bounded  away from $D$. Thus,  $\phi_j$ converges to some continuous function $\tilde \phi$ which satisfies 
$$\left(\omega_0+\sqrt{-1}\partial \bar\partial \tilde \phi \right)^n=e^{h_0-\theta_X-X(\tilde \phi)-\phi}|S|_H^{2}\omega_0^n,$$ 
and $\Delta_{\omega_0}\tilde \phi$ is bounded away from $D$.   Moreover,   the solution of the above equation is unique  in the $L^\infty$-space.  Consequently,  compared with (\ref{MAE-deg}),  we get $\tilde \phi=\phi$. Hence,  $\Delta_{\omega_0}\phi$ is  bounded.  By Evans-Krylov theory, we can also obtain the higher order estimates   away from $D$. 
\end{proof}

\section{Setting on   non-compact toric manifolds}

In the following sections, we concentrate the K\"ahler manifold 
$\bar M$ on toric manifolds and the divisor $D$ chosen as a multiple of smooth divisor of $\bar M$. 
Let $M=\bar M\setminus D$.
As the Fano case,  
studied in \cite{WZ04}, 
 it is natural to study torus invariant singular KR solitons    for a torus soliton VF $X$ on 
 a non-compact toric   manifold 
 $ \bar M$.  Then under an affine coordinates  $z=(z_1,..., z_n)$  on  a tours orbit of  $  M$,  where $z_i=x_i+\sqrt{-1} \theta_i$ and $x=(x_1,..., x_n)\in \mathbb R^n$, 
 $$X=\sum_i a^i\frac{\partial}{\partial z_i}$$
 for a vector $A_0=(a^1,...,a^n)\in \mathbb R^n$.
Thus  by writing $\omega=\sqrt{-1}\partial \bar\partial \varphi$ on $(\mathbb C^*)^n$,  where $\varphi$ is a smooth convex function on $\mathbb R^n$, 
 (\ref{twisted-equ-0}) can be  reduced to the following real MA  equation on $\mathbb R^n$,
\begin{align}\label{RMA-2-non}\det(\nabla^2 \varphi) =e^{-\varphi-<A_0, \nabla\varphi>}.
\end{align}
Note that 
$$\int_{\mathbb R^n}\det(\nabla^2 \varphi) e^{<A_0, \nabla\varphi>}
=\int_Pe^{<A_0,y>} dy,$$
where $P=\nabla\varphi(\mathbb R^n)$ is a non-compact  convex set in $\mathbb R^n$. 
Hence,
\begin{align}\label{integral}\int_Pe^{<A_0,y>} dy=\int_M e^{-\varphi}\le C_0.
\end{align}
This means that  $\int_Pe^{<A_0,y>} dy$
is well defined.
Actually one can show that $P$ is a non-compact polytope, which is an  intersection of several half spaces  in $\mathbb R^n$ given by
\begin{align}\label{polytope-complete}P=\cap_{\alpha=1,..., d} \{l_\alpha=1-<a_\alpha,x> >0\},
\end{align}
where all $a_\alpha\in  \mathbb Z^n$ (cf. \cite{Ci}).

Recall that an   $n$-dimension polytope $P^n$ is an  intersection of several half spaces  in $\mathbb R^n$ given by
\begin{align}\label{polytope}P=\cap_{\alpha=1,..., d} \{l_\alpha=\lambda_\alpha-<a_\alpha,x> >0\},
\end{align}
where all $a_\alpha\in  \mathbb Z^n$.

\begin{defi}\label{cone-polytope} An   $n$-dimension  non-compact polytope $P^n$  is called satisfying a  cone condition  if there is  an $a_0\in \mathbb Z^n$  and an integer $k_0\in \mathbb Z^+$ such that  for any integer $k\ge k_0$
\begin{align}\label{pk-l}P_k=P\cap \{l_0=k-<a_0,x> >0\}
\end{align}
is a Delzant polytope. In addition, all  constant $\lambda_\alpha$ are equal to $1$ in (\ref{polytope}),  i.e., in the case of  (\ref{polytope-complete}),  $P^n$ is called a Fano cone.

\end{defi}

By Definition \ref{cone-polytope}, each $P_k$ corresponds to a same  compact toric manifold $\bar M$.  For convenience, we call this toric manifold $\bar M$  associated to the pair $(P, a_0)$.  In  the case of  Fano cone $P$,  if one can choose $k_0=1$, $P_1$  just corresponds to a toric  Fano  manifold. In our setting, we focus on the case of Fano cone, but not necessary to assume that $\bar M$ is a Fano  manifold  associated  to $P_k$. On the other hand, with respect to the condition (\ref{integral}),  we  will  further  assume that  there exists a   vector $A_0\in \mathbb R^n$ and a positive constant $c_0$ such that for any $y\in P$ with  $|y|\gg 1$ it holds
\begin{align}\label{behavior-P-0}
<A_0, y>\le -c_0. 
\end{align}
(\ref{behavior-P-0}) implies that the angle between
the vectors $-A_0$ and $y$ is  strictly less than $\frac{\pi}{2}$ for any $|y|\gg 1$ in  $P$. For simplicity,  we call    $P$  a positive  polytope cone.

 As in \cite{TZ02}, we introduce a $H_0(\cdot)$-functional  on $\mathbb R^n$ by
 \begin{align}\label{h0-inv-infty}
H_0(B)=\int_{P} e^{<B, y>} dy=0, ~i=1,...,n.
\end{align}
By the assumption (\ref{behavior-P-0}), it is easy to see that there is a minimizer $A_0$ (we still use the same notation) such that
\begin{align}\label{minimizer}H_0(A_0)=\inf_{B\in \mathbb R^n} H_0(B).
\end{align}
Thus, we get
\begin{align}\label{futaki-inv-0}
f_{A_0}(e_i)=\int_{P} y_i e^{<A_0, y>} dy=0, ~i=1,...,n.
\end{align}
Actually,   the functional   $f_{A_0}(\cdot)$  is just the modified Futaki-invariant introduced in \cite{TZ02} when $P=P_1$ in the case of Fano polytope.  In this case,  it was proved that there always exists such an  $A_0=  (a^1,...,a^n)$ which satisfies (\ref{futaki-inv-0}). Hence, without confusion, we also call 
 $$X=\sum_i a^i\frac{\partial}{\partial z_i}$$
a toric soliton VF on a  non-compact toric  manifold $M$,  where $A_0$ is determined by (\ref{minimizer}).

Let $D$ be the divisor  associated to  the linear function $l_0$ in (\ref{pk-l}) and   $\Omega_k$ the  K\"ahler class  associated to $P_k$ ($k\ge k_0$) on   the toric manifold $\bar M=M_k$.  Then
$$ \Omega_k= 2\pi c_1(M)+2(k-1)\pi c_1([D]).$$ 
In particular,  if $\bar M$ is a  toric  Fano manifold,   $k=1$,   $\Omega_k$ is same to $2\pi c_1(\bar M)$.  By a result of Wang-Zhu \cite{WZ04}, there always exists a KR soliton on a such a manifold $\bar M$.  Thus we hope to generalize the above result to 
get a singular KR soliton  $\omega_k$  associated to  multiple divisor $(k-1)D$ on each  K\"ahler class $\Omega_k$ in  sense  of Definition \ref{admissable-0}. Namely,  there exists a   a torus soliton VF $X$ on $\bar M$ such that $\omega_k$ satisfies
\begin{align}\label{KR-soliton-toric}{\rm Ric}(\omega)-\omega+ (k-1)[D]= L_X\omega, ~{\rm on}~ \bar M.
\end{align}
Clearly,  if  the  family of singular KR solitons $\omega_k$ above  is convergent  away from $ D$,  we will get a limit KR soliton on  $M$ with non-compact moment polytope.
Hence,  in the next step,   we  shall do a prior estimates for  $\omega_k$ to  prove the convergence and completeness.  Both of two steps  will depend on the toric structure of $\omega_k$. Actually,  each   $\omega_k$  corresponds to a convex  solution of certain real MA equation 
as in (\ref{RMA-2-non}), to see (\ref{RMA-2-p}) below.

\subsection{Reduced MA equation} 

For a torus $K=(S^1)^n$  invariant metric $\omega_0\in \Omega_k=\in 2\pi c_1(M)+(k-1)[D]$,  there is a smooth convex function $\psi_0$ in $\mathbb R^n$ with
$\nabla\varphi(\mathbb R^n)=P_k$  such that
\begin{align}\label{convex-potential}\omega_0=\sqrt{-1}\partial \bar\partial \psi_0, ~{\rm on} ~(\mathbb C^*)^n.
\end{align}
 By (\ref{h-potential}), we have
\begin{align}\label{h0-potential}e^{h_0}|S|_H^{2(k-1)}\omega_0^n=e^{-\psi_0}\Pi_{i=1}^n d z_i \wedge d \bar z_i, ~{\rm on} ~(\mathbb C^*)^n,
\end{align}
where $h_0$ is the twisted Ricci potential of $\omega_0$ and $H$ is a Hermitian metric on $D$.
Note that  for a  a torus HVF $X$  there is a vector   $A=(a^i)$ in $\mathbb R^n$ such that
$X=a^i\frac{\partial}{\partial z^i}$ on $M$.
Thus,  for a  $K$-invariant  weak KR soliton  $\omega\in \Omega_k$ by
$$\omega_{\phi}=\sqrt{-1}\partial \bar\partial (\varphi)=\sqrt{-1}\partial \bar\partial (\psi_0+\phi),$$
where  $\phi$ is a  $K$-invariant K\"ahler potential of  $\omega_{\phi}$, we get from  (\ref{twisted-equ-0}), 
\begin{align}
&e^{<A, \nabla\varphi>}\omega^n=(|S|_H^{2(k-1)}) e^{ h_0-\phi}\Pi_{i=1}^n d z_i \wedge d \bar z_i\notag\\
& =e^{-\varphi}\Pi_{i=1}^n d z_i \wedge d \bar z_i, ~{\rm on} ~(\mathbb C^*)^n.  \notag
\end{align}
This means that    (\ref{twisted-equ-0})  can be  reduced to the following real MA  equation  as (\ref{RMA-2-non}) on $\mathbb R^n$,
\begin{align}\label{RMA-2-p}\det(\nabla^2 \varphi) =e^{-\varphi-<A, \nabla\varphi>}.
\end{align}

Conversely,   if there is  a function $\varphi=\psi_0+\phi$ satisfying $(\ref{RMA-2-p})$ and $|\phi|\leq C$ for some constant $C>0$,  then by (\ref{h0-potential}),   we have 
\begin{align}\label{inverse}
(\omega_0+\sqrt{-1}\partial \bar \partial \phi)^n&=e^{-\varphi-<A, \nabla\varphi>}\Pi_{i=1}^n 
d z_i \wedge d \bar z_i\notag\\
&=e^{-\phi-X(\phi)-\theta_0}e^{h_0}|S|_H^{2(k-1)}\omega_0^n,  ~{\rm on} ~(\mathbb C^*)^n,
\end{align}
where $\theta_0=<A, \nabla\psi_0>$.
Since both sides are measures charging no mass on analytic subsets,  the above relation can be extended to  the total space on $M$, i.e.,  $\phi$ is a weak solution of 
\begin{align}\label{singkr-k}
(\omega_0+\sqrt{-1}\partial \bar \partial \phi)^n=e^{-\phi-X(\phi)-\theta_0}e^{h_0}|S|_H^{2(k-1)}\omega_0^n,  ~{\rm on} ~M,
\end{align}
 which is the same as (\ref{MAE-deg}).

By (\ref{RMA-2-p}), it is easy to see that the vector $A=A_k$ should
satisfy
\begin{align}\label{futaki-inv}
f_A(e_i)=\int_{P_k} y_i e^{<A, y>} dy=0, ~i=1,...,n.
\end{align}
Actually,   the quantity  $F_A(\cdot)$  is just the modified Futaki-invariant introduced in \cite{TZ02}.  Thus $A$ is  uniquely   determined by a minimal  critical point of the following $H_0$-invariant for vectors in $\mathbb R^n$
defined by
\begin{align}\label{h0-inv}
H_0^k(B)=\int_{P_k} e^{<B, y>} dy.
\end{align}

As in \cite{BB, LZ, LTZ, TZ24},  we will use a variation method to solve the equation \ref{RMA-2-p}. Although each polytope $P_k$  does  not  satisfy the Fano condition at present, we  can also  prove that  the Ding-energy  $F(\cdot)$ associated to  (\ref{RMA-2-p})  is proper, thanks  the existence  of admissible K\"ahler metrics constructed in Section 3.1.

The following is a basic observation which will be used in Section 6 below.

\begin{lem}\label{Ak-constant}
\begin{align}\label{Ak-limit}\lim_k A_k=A_0.
\end{align}
\end{lem}

\begin{proof}
(\ref{Ak-limit})By the assumption  (\ref{behavior-P}), 
 $$\lim_k \inf_{B\in \mathbb R^n} H_0^k(B)=\inf_{B\in \mathbb R^n} H_0(B)=H_0(A_0). $$
Then (\ref{Ak-limit}) follows  immediately from   the fact  that  the critical point $A_k$  of  $H_0^k(B)$ is unique.
\end{proof}

\subsection{The space $\mathcal E^1(P_k)$}

As in   \cite{BBEGZ, HL},   the  $F(\cdot)$-functional     (or  Ding-energy) can be defined  on
a space $\mathcal E^1(\bar M,\Omega_k)$ by
\begin{align}\mathcal E^1(\bar M, \Omega_k)\,=\,\{\phi| ~&\phi ~{\rm has  ~full ~MA~mass~and}  \notag\\
& \sup_{\bar M}\phi =0,~ I(\phi)=\int_{\bar M}-\phi\omega_\phi^n<\infty\}.\notag
\end{align}
On a toric  manifold,    it is natural to restrict    $\mathcal E^1(M, \Omega_k)$  on  a subspace  by
  \begin{align}\label{e-space}&\mathcal E_K^1(\bar M, \Omega_k)\notag\\
=&\{\phi\in \mathcal E^1(\bar M ,\Omega_k)|~ \varphi=\psi_0+\phi  ~{\rm is~ a ~convex~ function~ on }~\mathbb R^n \}.
  \end{align}
Here  $\psi_0$ is a convex function  associated to a  K\"ahler metric  in  $\Omega_k$ as in (\ref{convex-potential}).
It is clear that  $\mathcal E^1_{K}(\bar M,\Omega_k)$ is locally precompact in terms of convex functions on $\mathbb R^n$.  We can also prove the  completeness $\mathcal E^1_{K}(\bar M,\Omega_k)$ by using the Legendre dual as in \cite{BB, LTZ}.

For any $\phi\in \mathcal E_K^1(\bar M, \Omega_k)$, we define   the Legendre function $u_\phi$ of $\varphi$ on $P_k$ by
\begin{align}\label{legendre-transform-u}u_\varphi(y)\,=\,\sup_{x\in\mathbb R^n}\{\langle x, y\rangle-\varphi(x)\}.
 \end{align}
 In the case of  K\"ahler potential $\phi\in \Omega_k)$,  $u_\varphi(\cdot)$  is just 
 $$u_\varphi(y)\,=\langle x,  \nabla\varphi (x)\rangle-\varphi(x),$$
 where $y=\nabla\varphi(x)$.
We can normalize $\psi_0$ by 
\begin{align}\label{normalization}
\inf_{\mathbb R^n} \psi_0\,=\, \psi_0(o)\,=\,0
\end{align}
as  well as  for any   $\varphi=\psi_0+\phi$ with $\phi\in \mathcal E_K^1(\bar M, \Omega_k)$,
 where $o$ is the origin in  $\mathbb R^n$.
As a consequence,   $u_\varphi(y)$  satisfies that
\begin{align}\label{normalization-u}
\inf_{P_k}u_\varphi\,=\,u_\varphi(o)\,=\,0.
\end{align}

We set a class of  convex functions on $P_k$ by
\begin{align}\mathcal E^1(P_k)\,=\,\{u| ~&\text{$u$ is convex on $P_k$ which satisfies
(\ref{normalization-u})  and  }
\notag\\
&\int_{P_k} u\,dy<+\infty\}.\notag
\end{align}
 Then for each $u\in \mathcal E^1(P_k)$,  we can define a  Legendre  dual function  $\varphi_u$ of $u$ on $R^n$ by
\begin{align}\label{legendre-u}\varphi_u(x)\,=\,\sup_{y\in\mathbb P_k}\{\langle x, y\rangle-u(y)\}.
 \end{align}
 It is easy to see that $\varphi_u$ is a convex function and its  Legendre  function  on $P_k$ is just $u$.


Let  
$$\mathcal E^0(P_k) =\{u\in \mathcal E^1(P_k)|~ u~{\rm is ~uniformly ~bounded}\}.$$   
Then one can  check that $\phi\in \mathcal E_K^1(\bar M, \Omega_k)$ is uniformly if and only if
$u_\varphi\in \mathcal E^0(P_k)$, where $ \varphi=\psi_0+\phi  $ (cf. \cite{Do02, SZ}).
The following result is due to \cite{LTZ}.

\begin{prop}\label{E1-legendre}
 A  K\"ahler potential $\phi\in\mathcal E^1_{K}(\bar M, \Omega_k)$ with
 normalized $\varphi$ satisfying  \eqref{normalization} if and only if the Legendre function $u_\varphi$ of $\varphi$ lies in $\mathcal E^1(P_k)$.
As a consequence, $u_\varphi$ is locally bounded in  Int$(P_k)$ if $\phi\in\mathcal E^1_{ K}(\bar M,  \Omega_k)$.
\end{prop}

\section{Construction of admissible K\"ahler metric}

In the following, we modify the Gulliemin metrics  in (\cite{Gu}) to construct a class of  admissible K\"ahler metrics  on $\bar M$  associated to  each $P_k$ for a  Fano cone $P$.

Let $l_\alpha$ and $l_0$ be linear functions in
 (\ref{polytope-complete}) and (\ref{pk-l}), respectively. 
Set
\begin{align}\label{u0-convex}
u_0=\sum l_\alpha\ln l_\alpha 
\end{align}
and 
\begin{align}\label{G-convex-F}
u_k=\sum l_\alpha\ln l_\alpha +\frac{1}{k} l_0\ln l_0.
\end{align}
Then
\begin{align}\label{2-deri}
\nabla^2u_k=\sum \frac{<\nabla l_\alpha, \nabla l_\alpha>}{l_\alpha} + \frac{1}{k}   \frac{<\nabla l_0, \nabla l_0>}{l_0}.
\end{align}

For a vertex $p\in \partial P_k\cap  \{l_0=<\alpha_0,x>- k=0\}$,   there are n facets  given by  $\{l_\alpha'=0\}$ which   intersect  at $p$  such that
the  determinant of   matrix  $T$  corresponding to $n$ normal vectors $\alpha'\in Z^n$   is  equal to one.  Thus by the transformation of  $T^{-1}$, we may assume  the n facets  $\{l_\alpha'=0\}$ are given by
\begin{align}\label{facets}
\{\lambda_i-y_i=0\},~i=2,...,n~{\rm and}~\{k-y_1=0\}.
\end{align}
Hence, any $y\in P_k$  near the facet $ \{l_0=k-y_1=0\}$, but away from  other facets,  we have
$$\nabla u_k(y)=(-\frac{1}{k} (\ln l_0+1), 0)+x',$$
where $x'\in \mathbb R^{n}$ is a uniformly  bounded vector.
It follows
\begin{align}\label{linear-func}
<\nabla u_k(y), y>=y_1x_1+O(1)
\end{align}
and
$$y_1=k+o(1)~{\rm and}~ x_1= -\frac{1}{k} (\ln l_0+1) +O(1).$$
The last relation implies
\begin{align}\label{x-y}l_0=\exp\{-kx_1+O(1)\},
\end{align}
and consequently,
$$y_1=k- \exp\{-kx_1+O(1)\},$$
where $O(1)$ denotes a smooth uniformly bounded function on $P_k$.
Therefore, again by (\ref{linear-func}), we get
\begin{align}\label{yk}<\nabla u_k(y), y>=kx_1-x_1\exp\{-kx_1+O(1)\}=kx_1+O(1).
\end{align}

By the above argument,  we also see that
    any $y\in P_k$  near each facet  $\{\lambda_i-y_i=0\},~i=2,...,n$,  but away from  other facets, it holds
\begin{align}\label{yi}    <\nabla u_k(y), y>=\lambda_ix_i-x_i\exp\{-x_i+O(1)\}=\lambda_ix_i+O(1).
\end{align}
Thus combining it with (\ref{yk}),  we conclude that near the  vertex $p\in \partial P_k\cap  \{l_k=<\alpha_0,x>- k=0\}$, it holds
\begin{align}\label{c0-potential}    <\nabla u_k(y), y>= kx_1+ \sum_{i=2, ..., n} \lambda_ix_i+O(1).
\end{align}

\begin{lem}\label{ricci-behavior}Suppose that $P$ is a Fano cone. Let
\begin{align}\label{legendre-phi}\varphi_k=   <\nabla u_k(y), y>-u_k
\end{align}
be the Legendre function of $u_k$ in (\ref{G-convex-F}).
Then
\begin{align}\label{ke-condition} (\ln {\rm det}(\nabla^2\varphi_k)+\varphi_k)=h_k
\end{align}
is a smooth uniformly bounded function.
\end{lem}

\begin{proof}It suffices to verify (\ref{ke-condition}) near each vertex.  Moreover,    the n facets intersected at the vertex satisfy (\ref{facets}).
Thus by (\ref{c0-potential}),  in the case of  Fano polytope $P_0$, we see
   \begin{align}\label{c0-potential-phi}    \varphi_k= kx_1+ \sum_{i=2, ..., n} x_i+O(1).
\end{align}
On the other hand, by (\ref{2-deri}),
$${\rm det}(\nabla^2\varphi_k)=({\rm det}(\nabla^2 u_k)^{-1}=ka(y)(k-y_1)\cdot (1-y_2)\cdot ...\cdot(1-y_n),$$
where $a(y)$ is a  uniformly bounded function on $P_k$.
Then by (\ref{x-y}) together  other $(n-1)$ similar relations
$$ (1-y_i)=\exp\{-x_i+O(1)\}, ~i=2,...,n,$$
it follows
\begin{align}\label{determinant}{\rm det}(\nabla^2\varphi_k)=\exp\{ -kx_1-\sum_{i=2,...,n}x_i +O(1)\}.
\end{align}
Hence, combining  the above relation with (\ref{c0-potential-phi}), we get (\ref{ke-condition}).

\end{proof}

\begin{rem}\label{ricci-class} By (\ref{determinant}), near the divisor $D$  associated to the facet $\{l_0=0\}$,   but away from  other facets, it holds
\begin{align}\label{near-lk}
{\rm Ric}(\omega_{\varphi_k})=&- \sqrt{-1} \partial\bar\partial\log(\det(\nabla^2\varphi_k ))+ (k-1)\sqrt{-1} \partial\bar\partial \ln  |z_1|^2\notag\\
&=k \sqrt{-1} \partial\bar\partial (x_1)+ (k-1)\sqrt{-1} \partial\bar\partial \ln  |z_1|^2+O(1)\notag\\
&=- \sqrt{-1} \partial\bar\partial \ln  |z_1|^2+O(1).
\end{align}
where $(z_1, z')$ are local coordinates near  a  neighborhood $U$ of regular part of $D$ such that  $D\cap U= \{z_1=0\}\cap U$.
Similarly, we have a same formula  of  Ricci form near each  divisor $D_\alpha$ associated to  facet $\{l_\alpha=0\}$ as (\ref{near-lk}). Thus we get \begin{align}\label{ricci-cohom}{\rm Ric}(\omega_{\varphi_k})\in  2\pi \sum_\alpha  c_1([D_\alpha])+ 2\pi  c_1([D])=2\pi  c_1(\bar M).
\end{align}
By  Lemma \ref{ricci-behavior} or directly by (\ref{c0-potential-phi}), we also prove
$$\omega_{\varphi_k} \in 2\pi \sum_\alpha  c_1([D_\alpha])+  2\pi k c_1([D]).$$
Consequently, 
$$\omega_{\varphi_k} \in 2\pi   c_1(\bar M)+  2\pi (k-1) c_1([D]).$$
\end{rem}

By the formula (\ref{2-deri}), we have the following explicit behavior of metric $\omega_{\phi_k}$ near  the   divisor   $D$.

 \begin{lem}\label{metric-behavior} Near the divisor $D$  associated to the facet $\{l_k=0\}$,   but away from  other facets, there are  local coordinates such that
\begin{align}\label{metric-behaviorD}
 \omega_{\varphi_k} &= (c_0+o(1)) |z_1|^{2(k-1)}dz_1\wedge d\bar z_1\notag\\
  &+ \sum_{i=1,...n} (a_i+o(1))  \frac{ |z_1|^{2k-1}}{|z_i|}  dz_1\wedge d\bar z_i + \sum_{i=1,...n} (\bar a_i +o(1)) \frac{ |z_1|^{2k-1}}{|z_i|}  dz_i\wedge d\bar z_1 +\omega',
\end{align}
where 
 $c_0>0$ and $a_i$ are  constants and
 $\omega'$ is a $(1,1)$-form such that
 \begin{align}\label{omega'}\omega'(\frac{\partial}{\partial z_1},\frac{\partial}{\bar\partial z_i})=0,~  \forall ~i=1,...,n;
~|\omega'(\frac{\partial}{\partial z_i}, \frac{\partial}{\bar\partial z_j})|\le A,~\forall ~ (i, j)\neq (1,1).
\end{align}

\end{lem}

\begin{proof}By  (\ref{2-deri}) and (\ref{x-y})  , we have
$$(u_k)_{1,1}=\frac{1}{kl_k}+o(1)=\exp\{kx_1+O(1)\}$$
and
$$|(u_k)_{i,j}|<A, ~\forall ~(i, j)\neq (1,1).$$
Note
$$(\nabla^2 \varphi_k)=(\nabla^2 u_k)^{-1}.$$
Then   it is easy to see that
$$(\varphi_k)_{1,1}=\exp\{-kx_1+O(1)\}; ~ (\phi_k)_{1,i}=\exp\{-kx_1+O(1)\}, ~\forall ~i=2,...,n; $$
and
$$ |(\varphi_k)_{i,j}|<A, ,~\forall ~ i, j=2,...,n.$$
By the coordinates transformation,
$$x_i=-\ln |z_i|^2, ~i=1,...,n,$$
we prove  that there is $(1,1)$-form  $\omega'$,  which satisfies  (\ref{omega'}),  such that
\begin{align} \omega_{\varphi_k}&= (c_0+o(1))\sqrt{-1}  |z_1|^{2k} d\ln z_1\wedge  d\overline {\ln z_1}\notag\\
&+ \sum_{i=1,...n} (a_i+o(1))  \sqrt{-1} |z_1|^{2k}   d\ln z_1\wedge d\overline{ \ln z_i} + \sum_{i=1,...n} (\bar a_i +o(1))  |z_1|^{2k}   d\ln z_i\wedge d\overline{ \ln z_1}+ \omega'\notag\\
& =(c_0+o(1))\sqrt{-1}  |z_1|^{2(k-1)}dz_1\wedge d\bar z_1\notag\\
&+ \sum_{i=1,...n} (a_i+o(1)) \sqrt{-1}  \frac{ |z_1|^{2k-1}}{|z_i|}  dz_1\wedge d\bar z_i + \sum_{i=1,...n} (\bar a_i +o(1))\sqrt{-1}  \frac{ |z_1|^{2k-1}}{|z_i|}  dz_i\wedge d\bar z_1+ \omega'.\notag
\end{align}
\end{proof}

\begin{cor}\label{metric-comparison} Let $\psi_k$ be the Gullimin metric associated to  the Gullimin function on $P_k$ by
$$
u_k^0=\sum l_\alpha\ln l_\alpha + l_k\ln l_k,$$
where $l_\alpha, l_k$ are same linear functions as in (\ref{G-convex-F}). Then there is a  uniformly bounded smooth  function $\varphi_k$ on $\mathbb R^n$ such that
\begin{align}\label{two-metric}\omega_{\varphi_k}-\omega_{\psi_k}=\sqrt{-1}\partial\bar\partial \varphi_k.
\end{align}
As a consequence,
$\omega_{\psi_k}\in [\omega_{\psi_k}].$
Moreover, 
\begin{align}\label{Gullimin-func}\omega_{\varphi_k}^n= e^{h_k}|z|_h^{2k-2}\omega_{\psi_k}^n,
\end{align}
where $h_k$is a uniformly smooth  bounded function and $h$ is a hermitian metric on the divisor $D$  with a defining section.
\end{cor}

\begin{proof}As in the computation of (\ref{yk}),
for any $y\in P_k$  near the  facet  $\{\l_k=k-y_1=0\}$,  but away from  other facets, it holds
\begin{align}\label{yk-0}    <\nabla u_k(y), y>=kx_1-x_1\exp\{-x_1+O(1)\}=kx_1+O(1).
\end{align}
Thus  in the case of  Fano polytope $P_0$,  as in (\ref{c0-potential-phi}), we get
   \begin{align}\label{c0-potential-phi-0}    \varphi_k= kx_1+ \sum_{i=2, ..., n} x_i+O(1).
\end{align}
It follows
$$\omega_{\varphi_k}=\sqrt{-1}\partial\bar\partial (kx_1+ \sum_{i=2, ..., n} x_i+h),$$
where $h$ is a  uniformly bounded smooth  function $\varphi_k$ on $\mathbb R^n$.
Hence, combining it with  (\ref{c0-potential-phi}), we obtain (\ref{two-metric}) immediately.

Similarly, as in the computation of (\ref{determinant}),  one can get
 \begin{align}\label{determinant-0}{\rm det}(\nabla^2\varphi_k^0)=\exp\{ -x_1-\sum_{i=2,...,n}x_i +O(1)\}.
\end{align}
Then under the local coordinates $(z_1, z')$ in Remark \ref{ricci-class},  there  is  a uniformly smooth  bounded function $\bar h_k$ such that
 $${\rm det}(\nabla^2\varphi_k)= e^{\bar h_k}{\rm det}(\nabla^2\varphi_k^0)|z_1|^{2k-2}.$$
The above  is equivalent to
$$\omega_{\varphi_k}^n= e^{\bar h_k}|z_1|^{2k-2}\omega_{\varphi_k^0}^n.$$
Thus by  choosing  a hermitian metric  $h$ on  the divisor $D$,  we get immediately,
(\ref{Gullimin-func}).
\end{proof}

\begin{rem}\label{admissable} By (\ref{Gullimin-func}), we have
$${\rm Ric}(\omega_{\varphi_k})-\omega_{\varphi_k}+ ( k-1)[D]= \sqrt{-1}\partial\bar\partial h_k,$$
which means  $\omega_{\varphi_k}$ is an admisiable metric on $\bar M$ in the K\"ahler class $\Omega_k=[\omega_k^0]$.

\end{rem}

\section{Solutions of (\ref{RMA-2-p}) and  (\ref{twisted-equ-0})}

Let  $v_{P_k}(x)$ be the support function on $\mathbb R^n$ defined  by
\begin{align}\label{support-function}
v_{P_k}(x)\,=\,\sup_{y\in P_k} \langle x,y\rangle.
\end{align}
Then by  Proposition  \ref{E1-legendre}, we see that for any $u\in{\mathcal  E}^1_{K}(P_k)$, its Legendre function
$$\varphi_u(x)\,=\,\sup_{y\in P_k}\{\langle x, y\rangle-u(y)\}\,\le\, v_{P_k}(x)$$
  corresponds to  a $K$-invariant weak K\"ahler potential $\phi_u=\varphi_u-\psi_0$ which belongs to $\mathcal E^1_{ K}(M,   \Omega_k)$.
   As we know,  $e^{-\phi_u}\in L^p(\omega_0)$ for any $p\ge 0$. Thus $\int_{\mathbb R^n}e^{-\varphi_u}dx$ is well-defined.

We introduce the following  
functional   on $\mathcal E^1(P)$, where $P=P_k$ by
  $$\mathcal D(u)\,=\,\mathcal L(u)+\mathcal F(u),$$
where
$$\mathcal L(u)\,= \int_{P}u e^{\sum a^iy_i}dy$$
and
$$\mathcal F(u)\,=-\log\left(\int_{\mathbb R^n}  e^{-\varphi_u}dx\right).$$
For a  K\"ahler potential  $\phi$  in  the  K\"ahler class $\Omega_k$, one can check 
 that
$$\mathcal L(u_\varphi)\,=\,-
 \sum_{k=0}^n\int_M \int_0^1 (1-t)^{n-k} t^k\phi e^{\theta_X+tX(\phi)}\omega_\phi^k\wedge\omega_0^{n-k} $$
and $\mathcal D(u_\phi)$ is just the modified  Ding-energy  $F(\phi)$   in  the  K\"ahler class $\Omega_k$ introduced in \cite{TZ02}.
For convenience, we call $\mathcal D(\cdot)$ the reduced Ding-energy  on $\mathcal E^1(P)$.

By (\ref{futaki-inv}),  for an affine function $l(y)$ on $\mathbb R^n$ it holds
$$\mathcal L(u+l(y))=\mathcal L(u).$$
Thus we can normalize $u\in \mathcal E^1(P)$ so that
$$u(o)=\inf_P u(\cdot)=0.$$
 We  will denote this class of  convex functions $u$ by   $\mathcal E_0^1(P)$.

The following result implies  that   $\mathcal D(u)$ is proper in  $\mathcal E^1(P)$.

\begin{prop}\label{ding-proper}   There are constants $\delta$ and $C_\delta$ such that
\begin{align}\label{proper-d}
\mathcal D(u)\,\geq\,  \delta \int_{P}u\,dy+ C_\delta,~\forall u\in\mathcal E^1_{0}(P).
\end{align}

\end{prop}

By Proposition \ref{ding-proper}, we are able to  study a  minimizer of  $\mathcal D(u_\varphi)$  to get a solution of  (\ref{RMA-2-p}) as in \cite{BB, LTZ}.  As a consequence, we prove

\begin{theo}\label{existence-1}  Let $P$ be a positive Fano polytope cone and $\bar M$ a toric manifold associated to the pair  $(P, a_0)$ defined in Definition  \ref{cone-polytope}.    Then on  any  K\"ahler class $(\bar M, \Omega_k)$ ($k\ge k_0$) 
 there exists a singular KR soliton $\omega_{\phi_k}\in \Omega_k$.
\end{theo}

\subsection{ Proof  of  Proposition \ref {ding-proper} }
 
The proof follows the line of argument in \cite{LZ}. We  note that $u_\phi$ satisfies the normalized condition
$u\geq u(0)=0$. Then
\begin{align}\label{linear-prop}
\mathcal L(u)\,\geq\, \lambda\int_{P}u\,dy,~\forall ~u\in {\mathcal  E}^1_0(P).
\end{align}

For the non-linear term $\mathcal F(\cdot)$, we can also get an analogy of \cite[Lemma~6.4]{LTZ} as follows.
Let $\varphi_0,\varphi_1\in\mathcal E^1_{ K}(M,\Omega_k)$ and $u_0,u_1$ be two Legendre functions of $\psi_0+\phi_0$ and $\psi_0+\phi_1$, respectively.
Let
 \begin{align}\label{linear-path}u_t=(1-t)u_0 +tu_1,~(t\in [0,1])
 \end{align}
  be a linear path connecting $u_0$ to $u_1$. Then by Proposition  \ref{E1-legendre},
the corresponding Legendre functions  $\psi_t$ of $u_t$  give a path in $\mathcal E^1_{K}(\bar M,\Omega_k)$.  The following
lemma shows that $\mathcal F(\psi_t)$ is convex in $t$ (cf. \cite[Lemma~6.4]{LTZ}).

\begin{lem}\label{F-convex} Let
$$\hat{\mathcal F}(t)\,=\,-\log\int_{\mathbb R^n}e^{-\psi_t}dx,~t\in [0,1].$$
Then $\hat{\mathcal F}(t)$ is convex in $t$ and so is $\mathcal F(\psi_t)$.
\end{lem}

\begin{proof}[Proof of Proposition  \ref{ding-proper}]Let $u_k$ be the Gullimin function as defined in (\ref{G-convex-F}) and $\varphi_k$ the Legendre function of $u_k$. Then 
Ricci potential of   $\omega_k$  (cf. Remark \ref{admissable}) is  uniformly  bounded on $M$. We set
\begin{eqnarray*}
A(y)\,=\,{\frac{V}{\int_{\mathbb R^n}e^{-\varphi_k}\mathbf dx}}e^{h_k(\nabla u_k(y))}, 
\end{eqnarray*}
which is $L^\infty$ function, where $y(x)=\nabla\varphi_k(x)$. Then the functional
\begin{eqnarray*}
\mathcal D_A(u)\,=\,\mathcal L^0_A(u)+\mathcal F(u),
\end{eqnarray*}
is well-defined  on   $\mathcal E^1(P)$, where
$$\mathcal L^0_A(u)\,=\,{\frac1V}\int_{P}u A(y)dy.$$
It is easy to see that $u_k$ is a critical point of $\mathcal D_A(\cdot)$.
On the other hand, by Lemma \ref{F-convex}, $\mathcal F(\cdot)$ is convex along any path in $\mathcal E^1_{K}(M, \Omega_k)$ determined by their Legendre functions  as in (\ref{linear-path}).
Note that  $\mathcal L_A^0(\cdot)$ is convex in $\mathcal E^1(P)$.  Hence
\begin{eqnarray*}
\mathcal D_A(u)\,\geq\,\mathcal D_A(u_0),~\forall u\in\mathcal E^1(P).
\end{eqnarray*}
Now  
 we can apply arguments in the proof of \cite[Proposition 4.9]{LZ} to prove that  there is a constant $C>0$ such that for any $u\in\mathcal E^1(P)$,
\begin{align}
\mathcal D(u)\,\geq\,{\frac{ C\lambda}{1+C}}\int_{P}u\,dy+\mathcal D_A(u_0)-n\log(1+C).\notag
\end{align}
Therefore, we get (\ref{proper-d}).

\end{proof}

\subsection{Proof of Theorem \ref{existence-1}} 

\begin{lem}
$\mathcal D(u)$ is lower semi-continuous on $\mathcal E^1_{0}(P)$.
\end{lem}
\begin{proof}

We will use the argument of Proposition 6.5 in \cite{LTZ}. By Fatou's lemma, $\mathcal L(u)$ is l.s.c.. Now assume that $u_i$ converges to $u\in \mathcal E^1_{0}(P)$ locally uniformly. Denote the Legendre transformation of $u_i, u$ by $\varphi_{u_i},\varphi_u$. For any $\epsilon_0 \in(0,1)$ there is a constant $M_{\epsilon_0}>0$ such that (cf. \cite[Lemma~2.7]{CGSZ}),
$$
\varphi_{u}(x) \geq\left(1-\epsilon_0\right) v_{ P}(x)-M_{\epsilon_0}, \forall x \in \mathbb R^n.
$$
The Legendre function $\varphi_{u_i}$ also converges locally uniformly to $\phi_{u}$. Then
$$
\partial \varphi_{u_i} \rightarrow \partial \varphi_{u}
$$
almost everywhere. Since
$$
\varphi_{u_i}(o)=\varphi_{u}(o)=0,
$$
we have
$$
\varphi_{u_i}(x) \geq\left(1-\epsilon_0\right) v_{ P}(x)-M_{\epsilon_0}, \forall x \in \mathbb R^n,
$$
if $i$ is large enough.  By the  dominant convergence, we know that 
$$\lim_{i \rightarrow \infty} \int_{\mathbb R^n}e^{-\varphi_{u_i}(x)}dx=\int_{\mathbb R^n} e^{-\varphi_{u}(x)}dx.$$ The lemma is proved.
\end{proof}

\begin{proof}[Proof of Theorem \ref{existence-1}]
By Proposition \ref{ding-proper}, we can choose a minimizing sequence $u_i$ of $\mathcal D(u)$ on $\mathcal E^1_{0}(P_k)$. Then $\int_{P_k} u_i dy$ is uniformly bounded. So $u_i$ converges to a convex function $u\in \mathcal E^1_{0}(P_k)$ locally by subsequence. By the lemma above, we know that $u$ is the minimizer of $\mathcal D(u)$ on $\mathcal E^1_0(P_k)$. Now $\mathcal D(u)$ is translation invariant, so $u$ is the minimizer on $\mathcal E^1(P_k)$. The Legendre transformation  $\varphi_u$ of $u$ is the minimizer of the  modified Ding functional $F(\cdot)$. 
Thus  $\varphi_u$ satisfies the Euler-Lagrange equation (cf. \cite{LTZ}),
\begin{align}\label{solution-rn} e^{<A_k, \nabla\varphi_u>}{\rm det}(\nabla^2 \varphi_u) =c_0e^{-\varphi_u},
\end{align}
where $c_0= \frac{\int_{\mathbb R^n}e^{<A_k, \nabla\varphi_u>}{\rm det}(\nabla^2 \varphi_u)  }{\int_{\mathbb R^n} e^{-\varphi_u}}.$
Consequently,   $\varphi_u$ is a solution of  
 (\ref{RMA-2-p}) by difference of a constant. 
  
  On the other hand,  since  $u_\infty$ is finite everywhere in  ${\rm Int}(P_k)$  by the locally uniform convergence and its Legendre function $\varphi_u\le v_{2P}$,  then  for any $\epsilon_0\in(0,1)$ there is a constant $M_{\epsilon_0}>0$ such that
(cf.  \cite[Lemma 2.6]{BB},  \cite[Lemma~2.7]{CGSZ}),
\begin{align}\label{13201}
\varphi_u(x)\,\geq\,(1-\epsilon_0)v_{2P}(x)-M_{\epsilon_0},\forall x\in\mathbb R^n.
\end{align}
It follows that $|x|^p{\rm det}(\nabla ^2 \varphi_u)$ is integrable for any $p>1$.  Choosing some $p>n$, we  get
$$\int_{P_k} |\nabla u|^p dy=\int_{\mathbb R^n} |x|^p {\rm det}(\nabla ^2 \varphi) dx<+\infty.$$  
 Thus $u$ is a continuous function   in $P_k$ by  the Sobolev embedding.   In particular,  $u$ is uniformly bounded, and so as   $\phi_u=\varphi_u-\psi_0$.  Hence,     $\phi_u$ can be extended as a  bounded K\"ahler potential function  on $M$.  
 Moreover,   the  real Monge-Amp\`ere measure  ${\rm det}(\nabla^2 \varphi_u)$ is same as the complex Monge-Amp\`ere measure $\omega_\varphi^n=(\sqrt{-1}\partial\bar\partial\varphi_u)^n$ (cf. \cite[Lemma 4.5]{LTZ}).  
 
 As in  (\ref{inverse}), we  have
 $$e^{-\varphi_u-<A_k, \nabla\varphi_u>}\Pi_{i=1}^n d z_i \wedge d \bar z_i
=e^{-\varphi_u- <A_k, \nabla\varphi_u>-\theta_0}e^{h_0}|S|_H^{2(k-1)}\omega_0^n,  ~{\rm on} ~(\mathbb C^*)^n.
$$
Since both of $\varphi_u$ and $<A_k, \nabla\varphi_u>$ are uniformly bounded,   by (\ref{solution-rn}),   $\phi_u$ can be  extended as a  weak solution of  CMA equation on $M$ as in (\ref{singkr-k}),  
$$
(\omega_0+\sqrt{-1}\partial \bar \partial \phi_u)^n=e^{-\phi_u-X_k(\phi_u)-\theta_0}e^{h_0}|S|_H^{2(k-1)}\omega_0^n,
$$
where $X_k$ is an extended torus HVF on $\bar M$ by the vector $A_k$.
  Namely,   $\phi_u$ satisfies 
 (\ref{MAE-deg}), and so   we prove that $\sqrt{-1}\partial \bar \partial \varphi_u$  can  be  extended as  a singular KR soliton in $\Omega_k$  on $\bar M$.
\end{proof}

\section{A priori estimates and local convergence of $\varphi_k$}

 Let   $\varphi=\varphi_k$ be the solution in  Theorem   \ref{existence-1} constructed in Section 5 and 
 $$m_k=\min_{\mathbb R^n} \varphi_k= \varphi_k(x_m).$$
 By the torus transformation, we may assume that $x_m=0$ and so
 $$\nabla\varphi_k(o)=0~ {\rm and }~\inf_{\mathbb R^n}\varphi_k=0,$$
    namely,   $\varphi_k$ attains the minimum $0$ at $o$.   Actually this is true according to the minimizer $u_k$ of $\mathcal D(\cdot)$ belonged to $\mathcal E^1_{0}(P)$, 
 Then $\varphi=\varphi_k$ is a 
  family of solutions   of Monge-Amp\`ere equations on $\mathbb R^n$,
\begin{align}\label{RMA-2}\det(\nabla^2 \varphi) =e^{-\varphi-m_k-<A_k, \nabla\varphi>}.
\end{align}
 Here $A_k$ are constants, which  converge to  $A_0$ as $k\to \infty$ by Lemma \ref{Ak-constant}.
 Moreover,   $\varphi$ satisfies   the following  further conditions:
\begin{align}\label{two-condition}
&1) ~<A_k, \nabla\varphi>(x', t)<0, ~{\rm as}~t\to\infty;\notag\\
&2)~|\nabla_n\varphi  (x', t)|\le C(t_0), ~{\rm as}~t\le t_0;\notag\\
&3)~|\nabla_{x'}\varphi|(x', t)|\le C_1;\notag\\
&4)~P_k=\nabla \varphi_k(\mathbb R^n)\to P, ~{\rm as}~k\to\infty.
\end{align}

In this section, we always assume that
$$P'=\{y'\in \mathbb R^{n-1}|~(y', y_n)\in \mathbb R^{n}\}$$
is a bounded set, i.e.  $P$ is a parabolic Fano cone.
Let $m_0=\lim_{k\to\infty} m_k$. 
We   purpose to prove

\begin{theo}\label{limit-solution} $\varphi_k$ converges  to a global solution  $\varphi_\infty$  of Monge-Amp\`ere equation on $\mathbb R^n$,
\begin{align}\label{RMA-limit}\det(\nabla^2 \varphi_\infty) =e^{-\varphi-m_0-<A_0, \nabla\varphi>},
\end{align}
which satisfies
\begin{align}\label{infty-gradient} P=\nabla \varphi(\mathbb R^n).
\end{align}

\end{theo}

In order to prove Theorem \ref{limit-solution},  our goal is  to get the following a priori estimate.

\begin{theo}\label{ck-estimate}
For any bounded domain $U\subset \mathbb R^n$,
\begin{align}\label{higher-order}\|\varphi_k\|_{C^l(U)}\le C_l(U).
\end{align}

\end{theo}

\subsection{Estimate of $m_k$}
We begin with   the following $C^0$-estimate.

\begin{prop}\label{co}There exists a uniform constant $m_0$ independent of $k$ such that
\begin{align}\label{c0-estimate}
|m_k|\le m_0.
\end{align}

\end{prop}

\begin{proof}First  we  estimate the lower bound of $m_k$.   Let
$W_{\frac {1}{2}}=\{x\in \mathbb R^n|~\varphi\le \frac{1}{2} \}, W_{\frac{3}{4}}=\{x\in \mathbb R^n| ~\varphi\le\frac{3}{4}\}, W_1=\{x\in \mathbb R^n|~\varphi\le 1\}$
to be three sets in $\mathbb R^n$,  respectively.
We consider a set
$$D_\delta= (W_{\frac{3}{4}}\setminus  W_{\frac{1}{2}})\cap \{(x', t)|~|x'| \le \delta \},$$
where $\delta$ is chosen to be small such that on  $\{(x', 0)|~ |x'| \le \delta  \}$ it holds
$$\varphi(\cdot)\le   \frac{1}{8}.$$

We claim: Set $W=D_\delta\cap \{(x', t)|~t\le 0\}$.  Then  there is a uniform constant $C$ such that 
\begin{align}\label{volume-small}
{\rm vol}(W) \ge C\delta^{n-1}.
\end{align}

Let $x=(x', 0)$ with  $|x'| \le \delta$ and $l_{x'}(t)=(x', t)$ a line staring from $x$, where $t\le 0$.   Then by the convexity,  $\partial_n\varphi (l_{x'}(t))\le 0$ 
 whenever  $l_{x'}(t)\subset  W$, in the  case of  $\partial_n\varphi( l_{x'}{(0)})\le 0$.
It follows that
\begin{align}\label{derivative}|\partial_n\varphi(l_{x'}(t))|\le C_0, ~{\rm on}~ D_\delta.
\end{align}
Let  $p'=(x', t_1)\in   \partial W_{\frac{1}{2}}$  and $p''=(x', t_2)\in   \partial W_{\frac{3}{4}}$ be two points  which lie  in  $l_{x'}$.  Then
$\varphi(p'')-\varphi(p')=\frac{1}{4}$.
Thus by  (\ref{derivative}),  we see that
\begin{align}\label{gap-t}t_2-t_1\ge \frac{1}{4C_0}.
\end{align}

We also need to show that (\ref{gap-t}) holds in  another case of $\partial_n\varphi( l_{x'}{(0)})> 0$. 
Actually,  in this case,  $\frac{1}{8}\ge l_{x'}{(0)}>{\rm {min}} \{l_{x'}(t)\}>0$, and so there exists $t'\le 0 $ such that  $\partial_n\varphi( l_{x'}{(t')})= 0$ and $\varphi( l_{x'}{(t')})<\frac{1}{8}$. Then  as  in the first case,  there are two points $p'\in \partial W_{\frac{1}{2}}$ and $p''\in  \partial W_{\frac{3}{4}}$, respectively, such that
$\partial_n\varphi(p')<0$ and $\partial_n\varphi(p'')<0$. As a consequence, we also get (\ref{gap-t}).  Hence,
$${\rm vol}(W)\ge  \frac{1}{4C_0} {\rm vol}(B_\delta'), $$
where  $ {\rm vol}(B_\delta')$ is the volume of $\delta$-ball  of $(n-1)$-dimension in $\mathbb R^{n-1}$.
(\ref{volume-small}) is proved.

By (\ref{RMA-2}), we see that
\begin{align}
  e^{-m_k-1}{\rm vol}(W)&\le  \int_W e^{-\varphi-m_k}= \int_W  e^{<A_k, \nabla\varphi>} \det(\nabla^2 \varphi) \notag\\
 &\le  \int_{\mathbb R^n}  e^{<A_k, \nabla\varphi>} \det(\nabla^2 \varphi) =\int_{P_k}  e^{<A_k, y>} \notag\\
 &\le A_0.\notag
 \end{align}
 Thus  by (\ref{volume-small}), we get
 \begin{align}\label{lower-mk}
 m_k\ge \ln \frac{{\rm vol}(W)}{A_0}-1\ge (n-1)\ln \delta-C'.
\end{align}

   To get  an upper bound of $m_k$, we follow  the argument  in \cite[Lemma~3.2]{WZ04}. 
     For any nonnegative integer $l$, we denote
$$V_l=\{x\in\mathbb R^n:\  l\le \varphi(x)\le l+1\}. $$
Then  all set $W_l=\bigcup_{i=0}^l V_i=\{\varphi< l+1\}$ are
convex  for any $l\ge 0$.
Note that
\begin{align}\label{barvraphi-inequ}
\det (\nabla^2\varphi) \ge c_0
e^{-\varphi-m_k},
\end{align}
where $c_0>0$ is a constant.  It follows that
$$ \det (\nabla^2\varphi) \ge c_0 e^{-1-m_k}, \ \ \ \text{in}\ \ W_1.$$
Thus the comparison  of  RMAE (cf.  \cite[Lemma~3.2]{WZ04}), we can get
\begin{align}\label{volume-w}
{\rm vol}(W_l)\le C_1e^{\frac{m_k}{2}}.
\end{align}
On the other hand, by the fact
$${\rm vol}(W_l)\le  (2(l+1))^n l^n {\rm vol}(W_1),$$
we have
\begin{align}
\int_{P_k}  e^{<A_k, y>}dy& =\int_{\mathbb R^n} e^{-\varphi-m_k} dx\notag\\
&\le\sum_l e^{-m_k-l}\int_{W_l} dx\le  \sum(2(l+1))^n e^{-m_k-l}\notag\\
&\le  C_2 e^{-m_k} {\rm vol}(W_1)\notag.
\end{align}
Consequently,
\begin{align}\label{volume-w-2}
{\rm vol}(W_1)\ge C_2e^{m_k}.
\end{align}
Hence, combining  (\ref{volume-w}), we obtain
$$m_k\le C'.$$

\end{proof}

\subsection{$C^1$-estimate}
Let $u=u_k$ be the Legendre function of $\varphi_k$.  
We have the following convergence result for $u_k$.

\begin{lem}\label{uk-convergence} $\mathcal L(u_k)$ is uniformly bounded and 
$u_k$ converges locally to a convex function $u_\infty$ on $P$.
\end{lem}

\begin{proof} It suffices  to show that $u_k$ is locally uniformly bounded.  First we bound
$\mathcal L(u_k)$.  Actually, by  a result in \cite[Lemmas~14--15]{CK},  for any convex function $u\in \mathcal C^1(P_k)$ it holds
\begin{align}\label{ck-inequality}
\int_{P_k}u e^{<A_k, y>}dy\ge \frac{c}{2\pi} (\int_{\mathbb R^n}  e^{-\varphi_u} dx)^{\frac{1}{n}}-(n+1),
\end{align}
where $c=\frac{c_1\kappa_n^{\frac{1}{n}}}{4e^{\frac{1}{n}}}, $
$$c_1=\inf_{\theta\in S^{n-1}} 
\{ \int_{P_k} <\theta, y> e^{<A_k, y>}dy\ \}$$
and $\kappa_n$ is a uniform constant depending the dimension $n$.
In particular, for $u=u_k$, we get
$$
\mathcal  D(u)+ \log \int_{\mathbb R^n}  e^{-\varphi_{u}} dx  \ge \frac{c}{2\pi} (\int_{\mathbb R^n}  e^{-\varphi_{u}} dx)^{\frac{1}{n}}-(n+1).$$
Note that $u$ is a minimizer of $\mathcal D(\cdot)$ 
 and so $\mathcal D(u)$ is uniformly bounded above. Thus 
  $\int_{\mathbb R^n}e^{-\varphi_{u}} dx$ is uniformly bounded and so $\mathcal L(u)$ does.

As in \cite{CK}, for any $y_0\in P_k$, we define a function  $f(\cdot)$  on  the unit sphere $S^{n-1}$ in $\mathbb R^n$ by
\begin{align}\label{f-function} f: \theta\in S^{n-1}\to \mu_k (\{y\in  P_k|~<(y-y_0), \theta>0\},
\end{align}
where $ \mu_k(E)=\int_E e^{<A_k, y>}dy$ for any subset $E\subset \mathbb R^n$.
Then it is easy to  see that there is a constant $m_{y_0}>0$ which depends  only on $|y_0|$ and $\delta={\rm dist}(y_0, \partial P)$,  but  is  independent of $k$  and $\theta$ such that
\begin{align}\label{measure-lower} f(\theta)\ge m_{y_0}
\end{align}
as so as $k\gg 1$.

Let $x_0=\nabla u_k(y_0)$. Then  $u_k\ge u_k(y_0) +<x_0, y-y_0>$. Thus
\begin{align}
\mathcal L(u_k)&=\int_{P_k} u_k  e^{<A_k, y>}dy\notag\\
 &\ge \int_{<x_0, y-y_0>0} u_k  e^{<A_k, y>}dy\notag\\
&\ge u_k(y_0) \int_{<x_0, y-y_0>0}  e^{<A_k, y>}dy\ge m_{y_0}u_k(y_0).\notag
\end{align}
 Hence, we get
\begin{align}\label{uk-upper-bound}u_k(y_0)\le  \frac{C}{m_{y_0}}.
\end{align}

\end{proof}

With  help of Lemma \ref{uk-convergence}, we prove

\begin{prop}\label{gradient-estimate}
For any bounded domain $\Omega\subset \mathbb R^n$,
\begin{align}\label{c1-norm}\|\varphi_k\|_{C^1(\Omega)}\le C_1(\Omega).
\end{align}

\end{prop}

\begin{proof}
On contrary, there is a sequence of points $x_k=(x'_k, t_k)\in \mathbb R^n$ such that
$$x_k'\to x_0'\in \mathbb R^{n-1},~t_k\to x_0, ~{\rm and}~\partial_n \varphi(x_k', t_k)\to \infty.$$
 Without loss of generality, we may also assume
\begin{align}\label{c0-blowup}
\varphi_k(x_k)\to\infty,
\end{align}
otherwise, we can take a new sequence  of $(x_k', t_k+1)$ instead of $x_k$. Clearly,
$\partial_n \varphi(x_k', t_k+1)\to \infty.$

Let
$D_1'= \{x\in \mathbb R^{n-1}|~ |x'|<1\}$ be a  ball with radius $1$ and
$$D_k'=\{(x'+x_k', t_k)|~ x'
\in D_1'\}$$
  be a sequence  of sets in $\mathbb R^n$.  Then
$  \varphi_k(x)\to\infty$ uniformly on $D_k'$ by  4) in (\ref{two-condition}).
We claim: for any $ x=(x', t_k) \in D_k'$,
\begin{align}\label{nabla-blowup}
\partial_n \varphi_k( x)\to\infty, ~{\rm as}~k\to\infty.
\end{align}

Let
$$D_0=\{(x'+x_k', 0)|~ x'
\in D_1'\}.$$
Since $\varphi_k(0)=0$,  $\varphi_k$ is uniformly bounded on $D_0$. Consider the restricted function   $l_k(t)=\varphi_k((x'+x_k', t)$  of $\varphi_k(\cdot)$ on an  each line  $(x'+x_k', t)$ with $t\in [0,t_k]$, where $x'\in D_1'$ and $x_k'\in  \mathbb R^{n-1}$ are both fixed.
Then  $l_k(0)$ is a uniformly  bounded number and $l_k(1)$ goes to the infinity as $k\to\infty$.  Thus    by (\ref{c0-blowup})
 there must  exist  $\hat t_k\in [0, t_k]$  such that
$$
l_k'(\hat t_k)\to\infty.$$
By the convexity,
$$
l_k'( t_k)\to\infty.$$
This proves (\ref{nabla-blowup}).

 By (\ref{nabla-blowup}),  as in the  choice of $x_k$ in (\ref{c0-blowup}), we see that  there is a subset
 $$\tilde D_k' =\{(x', t_k'(x'))|~,x' \in D_1', t_k'(x')\in [t_k, t_k+1]\}$$
 such that  for any $ x\in \tilde D_k' $ it holds
   \begin{align}\label{c0-blowup-2}
\varphi_k( x)\to\infty~{\rm and}~\partial_n \varphi_k( x)\to\infty, ~{\rm as}~k\to\infty.
\end{align}
  By  (\ref{c0-blowup-2}), we can perturb  $\tilde D_k' $ to an open  convex  smooth  hypersurface   $\bar D_k' $ in $\mathbb R^n$ such that   for any $x\in \bar D_k' $
   (\ref{c0-blowup-2}) still  holds.   Moreover, we can assume that  $\bar D_k' $ contains a subset which is diffeomorphic to a small $(n-1)$-dimensional hypercube in $\mathbb R^{n-1}$.  
   
   Let
  $V_k$ be a cone with the link $ \bar D_k' $.  We claim: for any integer $p>1$, it holds
  \begin{align}\label{w-estimate}
  \int_{V_k} |\nabla \varphi_k|^p dx\le C_p.
  \end{align}

 Let $\varphi_\infty(x)$ be the Legendre function of $u_\infty$ in Lemma \ref{uk-convergence}. 
Let  $D_\infty$ and $V_\infty$ be  two limit sets of $ \bar D_k' $ and $V_k$, respectively. 
Then for any $y=\nabla\varphi_\infty(x)$,  where $x\in V_\infty$, it holds 
$$|-u_\infty+ <y, \nabla u_\infty>|=o(|y|),$$
since  $\nabla u_\infty$ has a limit as $y\to\infty$. 
This also implies  that
\begin{align}\label{behavior-u}|-u_k+ <y, \nabla u_k>|=o(|y|), 
\end{align}
where $y=\nabla\varphi(x), ~x\in V_k.$

By  (\ref{behavior-u}), we  have
$$-u_k+ <y, \nabla u_k>+<A_k, y>\le \frac{1}{2} <A_k, y>+C_0.$$
Since $u_k$ satisfies
  \begin{align}\label{weak-KRRMA-global-u}
\det(\nabla^2 u_k) =e^{-u_k+ <y, \nabla u_k>+<A_k, y>},
\end{align}
we get
   \begin{align}
   \label{w-estimate-2} &\int_{V_k} |\nabla \varphi_k|^p dx= \int_{\nabla \varphi(V_k)}  |y|^p \det(\nabla^2 u_k) dy\notag\\
   &\le \int_{\nabla \varphi(V_k)}  |y|^p e^{\frac{1}{2} <A_k, y>+C_0} dy\le C,
   \end{align}
   which proves (\ref{w-estimate}).

  By the Sobolev  embedding result together with  (\ref{w-estimate}),  we know
  $$ |\varphi_k( x)|\le C, ~{\rm on}~ V_k.$$
  In particular,  $\varphi_k( x)$ is uniformly bounded on   $\bar D_k' $.  But this is impossible by the first inequality in  (\ref{c0-blowup-2}).  As a consequence,  (\ref{nabla-blowup}) can not be held. Therefore we get a contradiction and so  Lemma \ref{gradient-estimate} is proved.

\end{proof}

For any bounded domain $\Omega\subset \mathbb R^n$,  we define  a quantity
$$N_+(\Omega)=\max \{0, \sup\{t|~ (x', t)\in \Omega\}\}.$$
Then  from  the argument in the proof of Lemma \ref{gradient-estimate}, we see that the constant $C(\Omega)$ in (\ref{c1-norm}) actually depends on $\sup_{\Omega} \varphi_k$ and  $N_+(\Omega)$.    Thus we  prove

\begin{prop}\label{gradient-estimate-2}
For any bounded domain $\Omega\subset \mathbb R^n$,
\begin{align}\label{c1-1}\|\varphi_k\|_{C^1(\Omega)}\le C(\sup_{\Omega} |\varphi_k|, N_+(\Omega)).
\end{align}

\end{prop}

 \begin{defi}\label{good-shape-definition}
 A convex domain $\Omega$ is called to have a good shape if
 \begin{align}\label{radii}\max \{r_i\}/\min\{r_i\}\le c^*,
 \end{align}
 for a constant $c^*$ under controlled, where $r_1, ..., r_n$  are the radii of the minimal ellipsoid $E$ of $\Omega$.
 \end{defi}

\begin{lem}\label{good-shape}Let
$$\Omega_a^k=\{x\in\mathbb R^n|~\varphi_k< a\}$$
 be a level set of $\varphi_k$.  Then there is a uniform $r$-ball centered at $o$ contained in $\Omega_a^k$ and  the diameter of  $\Omega_a^k$ is uniformly bounded.
 As a consequence,  $\Omega_a^k$ has a good shape for any $a>0.$ In the other word,  the constant $c^*$ in
   Definition~\ref{good-shape-definition} associated to $\Omega_a^k$  depends only on $a$.
\end{lem}

 \begin{proof}
Fix $a>0$ and write $u_k=\varphi_k^*$ for the Legendre transform.
By normalization, $\varphi_k(0)=\min\varphi_k=0$ and
$u_k(0)=\min u_k=0$. Set
\[
 Z_k=\int_{P_k} e^{\langle A_k,y\rangle}\,dy,\qquad
 I_k=\int_{\mathbb R^n} e^{-\varphi_k(x)}\,dx,\qquad
 d\nu_k(y)=Z_k^{-1}e^{\langle A_k,y\rangle}\,dy
 \quad(y\in P_k).
\]
Changing variables $y=D\varphi_k(x)$ in \eqref{RMA-2} gives
$I_k=e^{m_k}Z_k$ and the identity of probability measures
\begin{equation}\label{sec6:shape-pushforward}
 (D\varphi_k)_*\bigl(I_k^{-1}e^{-\varphi_k(x)}\,dx\bigr)=\nu_k .
\end{equation}

For each fixed $k$, a small closed ball about $0$ lies in $P_k$.
The Legendre inequality on this ball gives
$\varphi_k(x)\ge b_k|x|-C_k$ with $b_k>0$.
Also, $D\varphi_k$ takes values in the bounded set $P_k$.
We may therefore integrate $\operatorname{div}(x e^{-\varphi_k})$
over $\mathbb R^n$. The Legendre identity then gives
\begin{equation}\label{sec6:shape-dual-average}
 \int_{P_k}u_k\,d\nu_k
 = I_k^{-1}\int_{\mathbb R^n}
       (x\cdot D\varphi_k-\varphi_k)e^{-\varphi_k}\,dx
 = n-I_k^{-1}\int_{\mathbb R^n}\varphi_k e^{-\varphi_k}\,dx
 \le n .
\end{equation}

By \eqref{integral}, $\langle A_0,d\rangle<0$ for every unit
recession direction $d$ of $\overline P$. Indeed, translating a
small ball in $P$ along such a direction would otherwise give an
infinite integral. These directions form a compact set and contain
every limit of $y_j/|y_j|$ with $y_j\in P$ and $|y_j|\to\infty$.
Since $A_k\to A_0$, there are $c>0$ and $C$ such that
$\langle A_k,y\rangle\le C-c|y|$ for $y\in P_k$ and all large $k$.
Choose $b>0$ so that $\overline B_{2b}(0)$ lies in all these $P_k$.
The exponential bound gives a uniform upper bound for $Z_k$, and
integration over this ball gives a positive lower bound. Hence
$\nu_k$ has uniformly bounded moments of every positive order and
its density has a uniform positive lower bound on $B_{2b}(0)$.

For $|z|\le b$, the ball $B_b(z)$ lies in $B_{2b}(0)$.
Jensen's inequality for its normalized Lebesgue measure, the
nonnegativity of $u_k$, and the bound on $\int u_k\,d\nu_k$ give
$u_k(z)\le C_0$ uniformly. Taking the supremum of $x\cdot z-u_k(z)$
over $|z|\le b$ now yields
\begin{equation}\label{sec6:shape-linear-growth}
 \varphi_k(x)\ge b|x|-C_0,\qquad x\in\mathbb R^n .
\end{equation}
Thus $I_k$ is uniformly bounded and
$\Omega_a^k\subset B_{R_a}(0)$ with $R_a=(a+C_0)/b$.

For the inner ball, fix $p>n$ and put $F_k=e^{-\varphi_k/p}$.
The change of variables above and the moment bound give
\begin{equation}\label{sec6:shape-sobolev}
 \int_{\mathbb R^n}|DF_k|^p\,dx
 =p^{-p}I_k\int_{P_k}|y|^p\,d\nu_k(y)\le C_p .
\end{equation}
Since $\|F_k\|_{L^p}^p=I_k$ and $F_k(0)=1$, Morrey's inequality
gives $1-F_k(x)\le C'_p|x|^{1-n/p}$. Choose $r_a>0$ so that
$C'_p r_a^{1-n/p}<1-e^{-a/p}$. Then $F_k>e^{-a/p}$ on $B_{r_a}(0)$,
and consequently
\begin{equation}\label{sec6:shape-inclusions}
 B_{r_a}(0)\subset\Omega_a^k\subset B_{R_a}(0).
\end{equation}
Shrinking $r_a$ and enlarging $R_a$ covers the finitely many
remaining indices, by continuity at $0$ and the fixed-$k$ growth
bound above.

The minimum-volume ellipsoid containing $\Omega_a^k$ has every
semiaxis at least $r_a$, by comparison of widths. Its volume is at
most that of $B_{R_a}(0)$, so its largest semiaxis is at most
$R_a^n/r_a^{\,n-1}$. The ratio in Definition~\ref{good-shape-definition}
is therefore at most $(R_a/r_a)^n$, uniformly in $k$.
\end{proof}

\subsection{$C^2$-estimate}

To get $C^2$-estimate, we consider a  more general  MA-type equation with Dirichlet problem:
\begin{equation}\label{MA-general}
\begin{cases}
\ \log\det u_{ij} = -u-X(u)-f &  \ \mbox{in $\Omega$}, \\[1mm]
\ u = 0 & \ \mbox{on $\partial\Omega$},
\end{cases}
\end{equation}
where  $\Omega\subset\mathbb{R}^{n}$ is  a bounded domain,  $f\in C^{\infty}(\overline{\Omega})$, $X$ is a  smooth vector field on $\overline{\Omega}$, and $u$ is a smooth convex  solution.
We establish the following  Pogorelov's estimate (cf. \cite{GT}).

\begin{lem}\label{c2-more}\footnote{This lemma was from a draft  written  by  Jianchun Chu last year. The authors would like to thank him for his contribution. } Let $u$ be a smooth convex  solution of (\ref{MA-general}). Then 
	there exist constants $\beta$ and $C$ depending only on $\|u\|_{C^0}, \|\nabla u\|_{C^{0}}$, $\|\nabla f\|_{C^0}$,$\|\nabla^2 f\|_{C^0}$, $\|X\|_{C^{2}}$ and $n$ such that
	\[
	\sup_{\Omega}\Big((-u)^{\beta}|D^{2}u|\Big) \leq C.
	\]
\end{lem}

\begin{proof}
	Since $u$ is convex, it then suffices to establish the upper bound of $D^{2}u$. For $(x,\xi)\in\Omega\times\mathbb{S}^{n}$, we consider the quantity
	\[
	Q(x,\xi) = \log u_{\xi\xi}+\frac{\alpha}{2}|Du|^{2}+\beta\log(-u),
	\]
	where $\alpha$ and $\beta$ are constants to be determined later. Suppose that $(x_{0},\xi_{0})$ is the maximum point of $Q$. Using $u|_{\partial\Omega}\equiv0$, $x_{0}$ must be an interior point of $\Omega$. We choose the coordinate system such that
	\[
	u_{ij}(x_{0}) = u_{ii}\delta_{ij}, \ \ u_{11} \geq u_{22} \geq \cdots \geq u_{nn}.
	\]
	Then $\xi_{0}=(1,0,\cdots,0)$ and the new quantity
	\[
	\hat{Q} = \log u_{11}+\frac{\alpha}{2}|Du|^{2}+\beta\log(-u)
	\]
	still achieves its maximum at $x_{0}$.
	
	\smallskip
	
	Define the elliptic operator
	\[
	L = u^{ij}\partial_{i}\partial_{j}+X^{i}\partial_{i}.
	\]
	The maximum principle shows $\hat{Q}_{i}=0$ and $L\hat{Q}\leq0$ at $x_{0}$. Then we obtain
	\begin{equation}\label{hat Q i}
		\frac{u_{11i}}{u_{11}}+\alpha u_{i}u_{ii}+\frac{\beta u_{i}}{u} = 0
	\end{equation}
	and
	\begin{equation}\label{L hat Q}
		\begin{split}
			0 \geq L\hat{Q}
			= {} & \frac{u^{ii}u_{11ii}+X^{i}u_{11i}}{u_{11}}-\frac{u^{ii}u_{11i}^{2}}{u_{11}^{2}}+\alpha u_{k}(u^{ii}u_{kii}+X^{i}u_{ki}) \\
			& +\alpha u^{ii}u_{ii}^{2}+\frac{\beta (u^{ii}u_{ii}+X^{i}u_{i})}{u}-\frac{\beta u^{ii}u_{i}^{2}}{u^{2}} \\[2mm]
			= {} & \frac{u^{ii}u_{ii11}+X^{i}u_{i11}}{u_{11}}-\frac{u^{ii}u_{11i}^{2}}{u_{11}^{2}}
			+\alpha u_{k}(u^{ii}u_{iik}+X^{i}u_{ik}) \\[1mm]
			& +\alpha\Delta u+\frac{\beta (n+X^{i}u_{i})}{u}-\frac{\beta u^{ii}u_{i}^{2}}{u^{2}}.
		\end{split}
	\end{equation}
	Applying $\partial_{k}$ and $\partial_{1}\partial_{1}$ on both sides of  the first relation in \eqref{MA-general} respectively, we obtain
	\[
	u^{ii}u_{iik}+X^{i}u_{ik} = -u_{k}-(\partial_{k}X^{i})u_{i}-f_{k}
	\]
	and
	\[
	u^{ii}u_{ii11}+X^{i}u_{i11} = u^{ii}u^{jj}u_{ij1}^{2}-u_{11}-2(\partial_{1}X^{i})u_{i1}-(\partial_{1}\partial_{1}X^{i})u_{i}-f_{11}.
	\]
	Then we see that
	\begin{equation}\label{eqn 1}
		\alpha u_{k}(u^{ii}u_{iik}+X^{i}u_{ik}) \geq -C\alpha
	\end{equation}
	and
	\begin{equation}\label{eqn 2}
		u^{ii}u_{ii11}+X^{i}u_{i11} \geq u^{ii}u^{jj}u_{ij1}^{2}-Cu_{11}-C.
	\end{equation}
	Substituting \eqref{eqn 1} and \eqref{eqn 2} into \eqref{L hat Q},
	\begin{equation}\label{eqn 3}
		\begin{split}
			0 \geq {} & \frac{u^{ii}u^{jj}u_{ij1}^{2}-Cu_{11}-C}{u_{11}}-\frac{u^{ii}u_{11i}^{2}}{u_{11}^{2}}
			-C\alpha+\alpha\Delta u+\frac{C\beta}{u}-\frac{\beta u^{ii}u_{i}^{2}}{u^{2}} \\
			\geq {} & \frac{u^{ii}u^{jj}u_{ij1}^{2}}{u_{11}}-\frac{u^{ii}u_{11i}^{2}}{u_{11}^{2}}-\frac{\beta u^{ii}u_{i}^{2}}{u^{2}}
			+\alpha u_{11}+\frac{C\beta}{u}-C\alpha-C-\frac{C}{u_{11}}.
		\end{split}
	\end{equation}
	For the third order term, using
	\[
	\frac{u^{ii}u^{jj}u_{ij1}^{2}}{u_{11}} \geq \sum_{i=1 \, \mathrm{or} \, j=1}\frac{u^{ii}u^{jj}u_{ij1}^{2}}{u_{11}}
	= \frac{u^{jj}u_{11j}^{2}}{u_{11}^{2}}+\sum_{i>1}\frac{u^{ii}u_{11i}^{2}}{u_{11}^{2}},
	\]
	then \eqref{eqn 3} becomes
	\begin{equation}\label{eqn 4}
		0 \geq \sum_{i>1}\frac{u^{ii}u_{11i}^{2}}{u_{11}^{2}}-\frac{\beta u^{ii}u_{i}^{2}}{u^{2}}
		+\alpha u_{11}+\frac{C\beta}{u}-C\alpha-C-\frac{C}{u_{11}}.
	\end{equation}
	By \eqref{hat Q i},
	\[
	\begin{split}
		\frac{\beta u^{ii}u_{i}^{2}}{u^{2}}
		= {} & \frac{\beta u^{11}u_{1}^{2}}{u^{2}}+\sum_{i>1}\frac{\beta u^{ii}u_{i}^{2}}{u^{2}} \\
		= {} & \frac{\beta u_{1}^{2}}{u^{2}u_{11}}+\sum_{i>1}\frac{u^{ii}}{\beta}\left(\frac{u_{11i}}{u_{11}}+u_{i}u_{ii}\right)^{2} \\
		\leq {} & \frac{\beta u_{1}^{2}}{u^{2}u_{11}}+\sum_{i>1}\frac{2u^{ii}u_{11i}^{2}}{\beta u_{11}^{2}}+\sum_{i>1}\frac{2u_{i}^{2}u_{ii}}{\beta} \\
		\leq {} & \frac{C\beta}{u^{2}u_{11}}+\frac{2}{\beta}\sum_{i>1}\frac{u^{ii}u_{11i}^{2}}{u_{11}^{2}}+\frac{Cu_{11}}{\beta}.
	\end{split}
	\]
	Substituting this into \eqref{eqn 4},
	\[
	0 \geq \left(1-\frac{2}{\beta}\right)\sum_{i>1}\frac{u^{ii}u_{11i}^{2}}{u_{11}^{2}}+\left(\alpha-\frac{C}{\beta}\right)u_{11}
	-\frac{C\beta}{u^{2}u_{11}}+\frac{C\beta}{u}-C\alpha-C-\frac{C}{u_{11}}.
	\]
	Choosing $\alpha=2$ and $\beta\gg1$, we obtain
	\[
	0 \geq u_{11}-\frac{C}{u^{2}u_{11}}+\frac{C}{u}-C-\frac{C}{u_{11}},
	\]
	and so
	\begin{equation}\label{eqn 5}
		0 \geq u^{2}u_{11}^{2}-C+uu_{11}-Cu^{2}u_{11}-Cu^{2}.
	\end{equation}
	Recalling $\hat{Q} = \log u_{11}+|Du|^{2}+\beta\log(-u)$, we have $e^{\hat{Q}}=(-u)^{\beta}u_{11}e^{|Du|^{2}}$. It then follows from \eqref{eqn 5} that
	\[
	0 \geq e^{2\hat{Q}}-Ce^{\hat{Q}}-C.
	\]
	By the Cauchy-Schwarz inequality, we obtain $\hat{Q}\leq C$, as required.
\end{proof}

\begin{proof}[Proof of  Theorem  \ref{ck-estimate}]
	Let  $\varphi$ be a solution of (\ref{RMA-2}) and $u=\varphi-m-2a$, where $m=\inf_{\mathbb R^n}\varphi$. Then  $u$ satisfies
	\begin{equation}\label{RMA-normalized} 
		\begin{cases} &\det(\nabla^2 u) =e^{-2a-m-u-<A_k, \nabla u>}, ~{\rm in}~ \Omega_{2a},\\
			&u=0, ~\partial  \Omega_{2a}.
		\end{cases}
	\end{equation}
	Then by Proposition  \ref{gradient-estimate} and Lemma \ref{c2-more},  we see
	$$\|\varphi\|_{C^2(\Omega_{a})}\le C(m, a).$$
	By the higher order estimate for the MA equation, it follows that
	\begin{align}\label{cl-u}\|\varphi\|_{C^l(V)}\le C(l, m, a, \delta), 
	\end{align}
	where
	$V=\{x\in \Omega_{a})|~{\rm dist}(x, \Omega_{a})\ge \delta\}.$
	On the other hand,  by  Lemma  \ref{gradient-estimate},  for any bounded domain  $U$ in $\mathbb R^n$, there is a large $a$ such that
	$$U\subset \Omega_{a}^k, ~\forall~k\gg 1.$$ 
	Hence, we get  (\ref{higher-order}) from (\ref{cl-u})  immediately. 
	
	\end{proof}

\begin{proof}[Another  proof of  Theorem   \ref{ck-estimate}]We can also give another   proof of Theorem \ref{gradient-estimate} by  Caffarelli's $C^{2,\alpha}$-estimate for the MA equation (cf.  \cite[Theorem~3.2(ii)]{WT}). Actually, 
	by  Lemma \ref{good-shape},    $\Omega_{2a}=\Omega_{2a}^k$  is uniformly bounded for any $a>0$. Then by  Lemma \ref{gradient-estimate},  the Monge-Amp\`ere equation  (\ref{RMA-2}) is locally strictly non-degenerate.   It follows that
	\begin{align}\label{c1+}\|\varphi_k\|_{C^{1,\alpha}(\bar\Omega_{2a})}\le C(a, m).\end{align}
	Thus as in (\ref{RMA-normalized} ), $u=\varphi-m-2a$ is a  convex solution of Monge-Amp\`ere equation on $\Omega_{2a}$, 
	\begin{align}\label{MA-initial}\det(\nabla^2 u) =f, 
	\end{align}
	where $\|f\|_{C^{\alpha}(\bar\Omega_{2a})}\le C(a, m). $
	
	On the other hand,  since $\Omega_{2a}$ has a good shade,  we may assume  that $\Omega_{2a}$ is 
	a normalized convex set. Then by an interior $C^{2,\alpha}$-estimate for (\ref{MA-initial}) in \cite[Theorem~4.42]{Fig}, 
	\begin{align}\label{c2-alpha}\|u\|_{C^{2,\alpha}(\Omega_{a})}\le C(a, m).
	\end{align}
	Hence, by the higher order estimate for the MA equation,   we also  get  (\ref{higher-order}).

\end{proof}

\begin{proof}[Proof of Theorem \ref{limit-solution}]
	By Theorem \ref{ck-estimate}, for any  $R$-ball  $B_R(0)$ centered at the original,
	there is a subsequence of $\{\varphi_k\}$ which converges to a limit solution of 
	(\ref{RMA-limit}) on  $B_R(0)$. By the  diagonal method, we will get a global solution 
	$\varphi_\infty$ of  (\ref{RMA-limit}) on  $\mathbb R^n$. It remains to check 
	(\ref{infty-gradient})  is satisfied.
	
	Let $u_k$ be the Legendre function of $\varphi_k$ on $P_k$
	as in Lemma \ref{uk-convergence}. By the local bound of $u_k$,
	for any bounded  open set $U\subset P_k$ $(k>>1)$, $\nabla u_k$ is uniformly bounded on $U$. This implies that  there exist an integer  $k_0$ and  a bounded  open set $K\subset \mathbb R^n$ such that   $U\subset \nabla\varphi_k( K)$ for any $k\ge k_0$.  By the above convergence of $\varphi_k$, 
	$U\subset \nabla\varphi_\infty( K)$. Thus we get (\ref{infty-gradient}).
\end{proof}

\begin{cor}\label{legendre-pair} Let  $u_\infty$ be the convex function on $P$ in Lemma \ref{uk-convergence}. Then  
	$u_\infty$ 
	is a Legendre function of $\varphi_\infty$,  which satisfies 
	\begin{align}\label{limit-legendre} u_\infty(y(x))=<x, \nabla\varphi_\infty(x)>-\varphi_\infty(x),
	\end{align}
	where $y=\nabla \varphi_\infty(x)$.
\end{cor}

\begin{proof} By both convergences of  $\varphi_k$ and $u_k$,  it easy to see that 
	the relation (\ref{limit-legendre}) is just the limit of 
	$$ u_k(\nabla\varphi_k(x))=<x, \nabla\varphi_k(x)>-\varphi_k(x)$$
	as $k\to\infty$.
	On the other hand, by Theorem \ref{limit-solution}, $\varphi_\infty$ is a strictly convex function on $\mathbb R^n$.  Then for any $y\in P$,  convex function
	$$ h_y(x)=<x,y>-\varphi_\infty(x)$$ 
	can attain the maximum at some point $x_y\in \mathbb R^n$. It follows that
	$y=\nabla \varphi_\infty(x_y)$. Namely,
	$$\sup\{x\in \mathbb R^n |~<x,y>-\varphi_\infty(x)\}=<x_y,  \nabla \varphi_\infty(x_y)>-\varphi_\infty(x_y)=u_\infty(\nabla \varphi_\infty(x_y)).$$
	Since   $\nabla \varphi_\infty$ gives a diffeomorphism from $\mathbb R^n$ to $P$, we prove that
	$u_\infty(\cdot)$ is the Legendre function of $\varphi_\infty$.

\end{proof}

By Lemma \ref{c2-more}, we also get  a uniform upper bound of  second-order estimate for $\varphi_\infty$.

\begin{prop}[]\label{sec7:pogorelov}
	For any  positive number $M_0$  there exists a constant $C(M_0)$ such that  for any $x\in \mathbb R^n$ with $|\nabla\varphi_\infty (x)|\le M_0$ it holds
	\begin{align}\label{c2upper-bound-uniform}
		D^2\varphi_\infty(x)\le C(M_0).
	\end{align}
	
\end{prop}

\begin{proof}
Fix $M_0>0$ and choose $R>|a_0|(M_0+1)$ sufficiently large that
the compact truncation
$K_R=\overline P\cap\{y:\langle a_0,y\rangle\le R\}$ has nonempty
interior. We use the gradient-image estimate of
Lemma~\ref{sec7:section-confinement}. Its proof uses
Lemma~\ref{bound} and the Legendre inequality
\eqref{sec7:superlinear}, both of which follow independently of this
proposition. Let $B$ and $\delta=\delta_R>0$ be the constants of
that lemma.

Fix $x_0$ with $|p_0|\le M_0$, where $p_0=D\varphi_\infty(x_0)$,
and set
\[
 h(x)=\varphi_\infty(x)-\varphi_\infty(x_0)-p_0\cdot(x-x_0),
 \qquad S=\{x:h(x)<\delta\}.
\]
Define $p_{0,k}=D\varphi_k(x_0)$ and
$h_k(x)=\varphi_k(x)-\varphi_k(x_0)-p_{0,k}\cdot(x-x_0)$.
Along the subsequence of Theorem~\ref{limit-solution},
$p_{0,k}\to p_0$, so $p_{0,k}\in K_R$ for all sufficiently large $k$.
For each fixed $x\in S$, local smooth convergence gives
$h_k(x)<\delta$ for all sufficiently large $k$.
Lemma~\ref{sec7:section-confinement} then places $D\varphi_k(x)$ in
the fixed compact set $K_B$. Passing to the limit gives
$D\varphi_\infty(S)\subset K_B$. In particular,
$|Dh|\le M_0+\max_{K_B}|y|$ throughout $S$.

Each such $S$ is a bounded domain. Indeed, $p_0$ is an interior
point of $P$, so choose $\varepsilon>0$ with
$\overline B_\varepsilon(p_0)\subset P^\circ$. The Legendre
inequality on this ball gives
$h(x)\ge\varepsilon|x|-C_{x_0}$, since $u_\infty$ is bounded on
the ball. Strict convexity gives $h(x_0)=0$ as its unique minimum,
so $\{h=\delta\}$ is a smooth level set.

Set $v=h-\delta$ on $S$. Then $-\delta\le v\le0$ and
$v|_{\partial S}=0$. Equation~\eqref{RMA-limit} gives
\begin{equation}\label{equ-h}
 \log\det D^2v=-v-A_0\cdot Dv-f,\qquad
 f(x)=p_0\cdot(x-x_0)+\varphi_\infty(x_0)+m_0+A_0\cdot p_0+\delta .
\end{equation}
Here $|Df|=|p_0|\le M_0$ and $D^2f=0$.
The bounds on $v$ and $Dv$ therefore allow us to apply
Lemma~\ref{c2-more} with uniform constants. At $x_0$ we have
$-v(x_0)=\delta$, and hence
$D^2\varphi_\infty(x_0)=D^2v(x_0)\le C\delta^{-\beta}I$.
Since $R$ and $\delta$ were fixed in terms of $M_0$, this gives
\eqref{c2upper-bound-uniform}.
\end{proof}

Both of Corollary \ref {legendre-pair} and Proposition \ref{sec7:pogorelov} will be used in next section.

\section{Global smoothness of the limiting soliton}
\label{sec:boundary-estimates}

In this section,  we prove the regularity of $\omega_\infty=-\sqrt{-1}\partial\bar\partial\varphi_\infty$ on 
$\bar M\setminus  U_{D,\epsilon}$, where $U_{D,\epsilon}$ is any  $\epsilon$-neighborhood of $D$.  We retain the parabolic assumption and the normalization of Section~6.
In this section $u_k$ always denotes the Legendre transform of the
solution $\varphi_k$, rather than the reference potential in
\eqref{G-convex-F}. We work along the convergent subsequence furnished
by Theorem~\ref{limit-solution}, and continue to denote its parameter
by $k$. Thus $\varphi_k(0)=\min\varphi_k=0$, $A_k\to A_0$, and
$m_k\to m_0$. Taking the Legendre transform of \eqref{RMA-2} gives
\begin{equation}\label{sec7:dual-equation}
 \log\det D^2u_k
 =\langle y,Du_k\rangle-u_k+\langle A_k,y\rangle+m_k.
\end{equation}
In particular, the additive normalization contributes $+m_k$ to the
dual equation. We write $P^\circ$ for the interior of the polyhedron
and $\overline P$ for its closure, including all its finite faces.
The notation $\omega_k$ below refers to the normalized solution
$\omega_{\phi_k}$, which is smooth and positive on
$M=\bar M\setminus D$ for each fixed $k$ by
Theorem~\ref{partial-smooth}.

Recall the reference function from \eqref{u0-convex},
\[
 u_0(y)=\sum_{\alpha=1}^d l_\alpha(y)\log l_\alpha(y),
 \qquad l_\alpha(y)=1-\langle a_\alpha,y\rangle.
\]
These are
precisely the retained facets. The moving facet
$l_{0,k}=k-\langle a_0,y\rangle=0$ is disjoint from every fixed
compact subset of $\overline P$ for sufficiently large $k$.

We will prove

\begin{theo}\label{smoothness}
The form
$\omega_\infty=\sqrt{-1}\partial\bar\partial\varphi_\infty$ on the
dense orbit extends to a smooth positive K\"ahler metric on all of
$M=\bar M\setminus D$. Along the sequence fixed in Section~6,
\begin{equation}\label{sec7:smooth-metric-convergence}
 \omega_k\longrightarrow\omega_\infty
 \quad\hbox{in }C^\infty_{\mathrm{loc}}(M).
\end{equation}
The extended metric satisfies the KR soliton equation
on $M$ with soliton vector field determined by $A_0$. Moreover, the metric $\omega_\infty$ is complete.

\end{theo}

Theorem \ref{smoothness} implies 
Theorem \ref{main-resullt}.

\subsection{Uniform $C^0$-estimate}

For $R$ sufficiently large, put
\[
K_R=\overline P\cap\{\langle a_0,y\rangle\le R\},
\qquad \Omega_R=\operatorname{Int}K_R.
\]
The cone condition implies that $K_R$ is a compact polytope with
nonempty interior, and these sets exhaust $\overline P$. For $k>R$,
the moving facet of $P_k$ does not meet $K_R$. All estimates below
are independent of sufficiently large $k$; their constants may depend
on the fixed region $K_R$.

\begin{lem}\label{bound}
For every $R$, every $1<p<\infty$, and every $0<\alpha<1$, there are
constants $C_{R,p}$ and $C_{R,\alpha}$ such that
\begin{equation}\label{sec7:holder}
\|u_k\|_{W^{1,p}(\Omega_R)}\le C_{R,p},
\qquad
\|u_k\|_{C^\alpha(K_R)}\le C_{R,\alpha}.
\end{equation}
Along the sequence fixed in Section~6, $u_k\to u_\infty$ uniformly
on every $K_R$. In particular, $u_\infty$ extends continuously to
$\overline P$ and belongs to $C^\alpha_{\mathrm{loc}}(\overline P)$
for every $\alpha<1$.
\end{lem}

\begin{proof}
The uniform bound for the diameter of the sublevel set
$\{\varphi_k<1\}$ in Lemma~\ref{good-shape}, together with
$\varphi_k(0)=0$, gives constants $c>0$ and $C$ such that
\begin{equation}\label{sec7:linear-growth}
\varphi_k(x)\ge c|x|-C\qquad(x\in\mathbb R^n).
\end{equation}
Indeed, choose $R_*>0$ with $\{\varphi_k<1\}\subset B_{R_*}(0)$
for every $k$. Convexity along each ray from the origin gives
$\varphi_k(x)\ge |x|/R_*$ for $|x|\ge R_*$, after increasing $R_*$
if necessary.

The Legendre relations, \eqref{RMA-2}, and the uniform bounds for
$m_k$ and $A_k$ imply
\begin{align}\label{sec7:gradient-integral}
\int_{\Omega_R}|Du_k|^p\,dy
&=\int_{Du_k(\Omega_R)}
 |x|^p\det D^2\varphi_k(x)\,dx\notag\\
&=\int_{Du_k(\Omega_R)}
 |x|^p e^{-\varphi_k(x)-m_k-\langle A_k,D\varphi_k(x)\rangle}\,dx
 \notag\\
&\le C_R\int_{\mathbb R^n}|x|^p e^{-c|x|}\,dx\le C_{R,p}.
\end{align}
Here $D\varphi_k(x)\in\Omega_R$ on the domain of integration,
which bounds the factor involving $A_k$. Choose a fixed ball
$B\Subset\Omega_R$. Lemma~\ref{uk-convergence} bounds $u_k$ on $B$.
The Poincar\'e inequality with the average over $B$ as normalization
therefore gives the full $W^{1,p}(\Omega_R)$ bound.
Since $\Omega_R$ is a fixed bounded Lipschitz domain, Morrey's
inequality, with $p$ chosen so that $1-n/p>\alpha$, proves
\eqref{sec7:holder}.

The equicontinuity on $K_R$ and the interior convergence from
Lemma~\ref{uk-convergence} identify every uniformly convergent
subsequence with the same continuous extension of $u_\infty$.
This proves convergence on $K_R$ and compatibility of the
extensions as $R$ increases.
\end{proof}

We shall also need the following consequence of the entire Legendre
potential. For every $L>0$,
\begin{equation}\label{sec7:superlinear}
u_\infty(y)\ge L|y|-\sup_{|x|\le L}\varphi_\infty(x),
\qquad y\in\overline P.
\end{equation}
For $y\in P^\circ$, this is the defining inequality for the Legendre
transform, evaluated at $x=Ly/|y|$ when $y\ne0$; at $y=0$ it follows
by taking the supremum over $|x|\le L$. Lemma~\ref{bound} extends the
inequality to $\partial P$. Thus $u_\infty$ grows faster than every
fixed linear function along every unbounded ray in $\overline P$.

\begin{lem}\label{uniformly-c0} Given any  number $M_0>0$. 
	Then there exists an integer $k_0>0$
	depending  only on $|y|=M_0$ such that for any $x\in \mathbb R^n$ with  $|\nabla\varphi_\infty(x)|\le M_0$ it holds
	\begin{align}\label{c0-order-limit}|\bar\phi(x)| \le C(M_0).\end{align}
	where $\bar\phi(x)=\varphi_\infty-\psi_{k_0}$ and $\psi_{k_0}$  is the convex function  on $\mathbb R^n$  defined as in Corollary \ref{metric-comparison}.
	associated to the polytope $P_{k_0}$.
\end{lem}

\begin{proof}
	Suppose that  $|y_0=\nabla\varphi_\infty(x_0)|\le M_0$.  Then we choose $P_{k_0}$ such that $y_0\in P_{k_0}$ for some $k_0$.    Recall the support function  $ v_{P_{k_0}}(x)$ on $\mathbb R^n$ defined  by
	(\ref{support-function}), 
	$$v_{P_{k_0}}(x)=\max_{p\in P_{k_0}}\{<x, p>\}.$$
	It is known that
	\begin{align}\label{two-equivalent} |\psi_{k_0}-v_{k_0}|\le  b_1(k_0).
	\end{align}
	Thus,  by Corollary \ref{legendre-pair},  for any $x\in \Omega= \nabla u_\infty(P_{k_0})$  we get
	\begin{align}\label{upper-bound-phi} \bar\phi(x)&=\varphi_\infty(x)-\psi_{k_0}(x)\notag\\
		&\le  \varphi_\infty(x) -v_{k_0}(x)+ b_1(M_0)\le   <x, \nabla u_\infty(y)> -u_\infty(y) -v_{k_0}(x)+ b_1(M_0)\notag\\
		&\le  b_1(M_0).
	\end{align}
	In particular, 
	$$\bar\phi(x_0)\le  b_1(M_0).$$
	
	We extend the above argument to prove  (\ref{c0-order-limit}).   Let $\tilde u_\infty$ be a restricted convex function  of $u_\infty$ on $P_{k_0}$ and
	$\tilde \varphi_\infty(x)$ its  Legendre dual  function on $\mathbb R^n$. Then as in (\ref{two-equivalent}), we have
	$$|\psi_{k_0}-  \tilde \varphi_\infty(x)|\le  b_1'(k_0), ~\forall~x\in \mathbb R^n.$$
	It remains  to show that  for any $x_0\in \mathbb R^n$ with  $|y_0=\nabla\varphi_\infty(x_0)|\le M_0$ it holds
	\begin{align}\label{different-u}
		|\tilde \varphi_\infty(x_0)-\varphi_\infty(x_0)|
		\le   b_2(M_0).
	\end{align}
	Actually, we claim:
	For any $x\in  \Omega$ it holds
	\begin{align}\label{equivalent}
		\tilde \varphi_\infty(x)=\varphi_\infty(x).
	\end{align}
	Hence,  (\ref{c0-order-limit}) is obtained.

	For any  $x=\nabla u_\infty(y)\in \Omega$,   we consider a convex function  $F(y')=  <x, y'>- u_\infty(y')$ on  $P_{k_0}$.  Clearly, $y$ is a unique critical point of  $F(y')$.  It follows that
	$$
	F(y)=\max_{y'\in P_{k_0}} ( <x, y'>- u_\infty(y')).$$
	By Corollary \ref{legendre-pair}, we  see that
	\begin{align}\label{legendre-2}\varphi_\infty (x)
		=\max_{y'\in P_{k_0}} ( <x, y'>- u_\infty(y')).
	\end{align}
	Recall that
	\begin{align}\tilde \varphi_\infty (x)
		&=\max_{y'\in P_{k_0}} ( <x, y'>-\tilde u_\infty(y'))\notag\\
		&=\max_{y'\in P_{k_0}} ( <x, y'>- u_\infty(y'))\notag.
	\end{align}
	Thus (\ref{equivalent}) is true.

	\end{proof}

The argument in  the proof of Lemma \ref{uniformly-c0} also implies

\begin{lem}\label{sec7:section-confinement}
For every fixed $R$, there are $B>R$, $\delta_R>0$, and $k_R$ with
the following property. If $k\ge k_R$ and
$p_0=D\varphi_k(x_0)\in K_R\cap P_k^\circ$, then
\begin{equation}\label{sec7:confined-section}
D\varphi_k\bigl(\{x:
\varphi_k(x)-\varphi_k(x_0)-p_0\cdot(x-x_0)<\delta_R\}\bigr)
\subset K_B.
\end{equation}
\end{lem}

\begin{proof}
Otherwise choose $k_j\to\infty$, $p_j\in K_R\cap P_{k_j}^\circ$,
and $q_j\in P_{k_j}^\circ$ such that $|q_j|\to\infty$ and
\begin{equation}\label{sec7:bregman}
0\le u_{k_j}(p_j)-u_{k_j}(q_j)
-Du_{k_j}(q_j)\cdot(p_j-q_j)\longrightarrow0.
\end{equation}
This is precisely the height of the section in
\eqref{sec7:confined-section}, expressed by Legendre duality.
Pass to a subsequence such that $p_j\to p_*$ and
$d_j=(q_j-p_j)/|q_j-p_j|\to d$, where $|d|=1$.
Let $\lambda_j$ be the supporting affine function of $u_{k_j}$ at
$q_j$, and write $e_j=u_{k_j}-\lambda_j$. Convexity gives
\[
0\le e_j\le e_j(p_j)\quad\hbox{on }[p_j,q_j],
\qquad e_j(p_j)\longrightarrow0.
\]
Fix any three ordered points on the ray $p_*+\mathbb R_{\ge0}d$,
and approximate them by points at the same fixed distances along
$[p_j,q_j]$. Their three-point Jensen differences for $u_{k_j}$
tend to zero, because subtracting $\lambda_j$ leaves these
differences unchanged. These points lie in a fixed compact
truncation, so Lemma~\ref{bound} passes the equality to the limit.
The ray is contained in $\overline P$, and $u_\infty$ is affine on
it. This contradicts \eqref{sec7:superlinear}.
\end{proof}

Together with Lemma \ref{sec7:section-confinement},   as in the proof of Proposition \ref{sec7:pogorelov}  in Section 6, we  get   

\begin{prop}\label{sec7:dual-lower-prop}
For every fixed $R$, there are $c_R>0$ and $k_R$ such that
\begin{equation}\label{sec7:dual-lower}
D^2u_k(y)\ge c_R I,
\qquad y\in K_R\cap P^\circ,\quad k\ge k_R.
\end{equation}
\end{prop}

\begin{proof}
Fix $p_0=D\varphi_k(x_0)\in K_R\cap P^\circ$ and set
\[
h(x)=\varphi_k(x)-\varphi_k(x_0)-p_0\cdot(x-x_0).
\]
The section $\{h<\delta_R\}$ is relatively compact: since
$p_0\in P_k^\circ$, the Legendre inequality on a small closed ball
around $p_0$ gives $h(x)\ge\varepsilon|x|-C$, with constants that
may depend on this section. Lemma~\ref{sec7:section-confinement} gives the uniform bound $|Dh|\le C_R$ on the section. Moreover,
\[
\log\det D^2h
=-h-p_0\cdot(x-x_0)-A_k\cdot Dh+d_{k,x_0},
\]
where $d_{k,x_0}$ is some constant.  Thus we can apply Lemma \ref{c2-more} (also see Proposition \ref{sec7:pogorelov} ),  with $f=-p_0\cdot(x-x_0)+d_{k,x_0}$ and $X=A_k$,
to obtain $D^2\varphi_k(x_0)\le C_RI$. The constants do not depend
on $d_{k,x_0}$ or on $x_0$. Taking inverse matrices in the Legendre
identity proves \eqref{sec7:dual-lower}.
\end{proof}

 Proposition \ref{sec7:dual-lower-prop}  will be used to construct a barrier function in deriving   a uniform $C^2$-estimate below.

\subsection{A proper toric coordinate and uniform localization}

We use the parabolic assumption of Section~6 to relate bounded moment
regions to fixed compact subsets of the complex manifold.  Throughout
this subsection, $Z=(Z_1,\ldots,Z_n)$ denotes multiplicative coordinates
on the dense torus, with
\[
 x_i=\log|Z_i|^2,\qquad \mu_k=D\varphi_k.
\]
This convention corresponds to coefficient $1$ in the retained-facet
terms $l_\alpha\log l_\alpha$.  For each fixed $k$,
Theorem~\ref{partial-smooth} gives a smooth toric K\"ahler metric on
$M$ and a smooth moment map $\mu_k$ there.  Only this individual
smoothness is used until uniform estimates are established.

\begin{lem}\label{sec7:geo-proper-character}
After an integral change of coordinates, and reversing the last
coordinate if necessary, the recession cone of $\overline P$ is
$\mathbb R_{\geq0}e_n$.  There are a compact polytope
$Q\subset\mathbb R^{n-1}$ and a number $T_0$ such that
\begin{equation}\label{sec7:geo-product-tail}
 \overline P\cap\{y_n\geq T_0\}=Q\times[T_0,\infty).
\end{equation}
The toric character $w=Z_n$ extends to a proper holomorphic map
\begin{equation}\label{sec7:geo-proper-map}
 w:M\longrightarrow\mathbb C.
\end{equation}
In particular, $\{|w|^2\leq e^b\}$ is compact in $M$ for every real $b$.
\end{lem}

\begin{proof}
Boundedness of the projection of $P$ to the first $n-1$ coordinates
forces its recession cone to lie in the last coordinate line.
The cone is nonzero because $P$ is an unbounded polyhedron, and it
cannot contain both directions of that line because the truncations
$P_k$ are bounded.  Thus it is a single ray.  If a preliminary change
of coordinates is needed, this ray has a primitive integral generator:
the defining facet normals are integral, so the recession cone is
rational.  A primitive integral vector can be completed to a lattice
basis.  We therefore make an integral change sending the ray to
$\mathbb R_{\geq0}e_n$, and express all moment, logarithmic, and soliton
data in the corresponding dual coordinates.

Write the retained primitive inward normals as
$\nu_\alpha=-a_\alpha$.  They satisfy $(\nu_\alpha)_n\geq0$.
The vertical facets, for which $(\nu_\alpha)_n=0$, define
\[
 Q=\{y'\in\mathbb R^{n-1}:
       1+\langle\nu_\alpha',y'\rangle\geq0
       \text{ whenever }(\nu_\alpha)_n=0\}.
\]
This is exactly the projection of $\overline P$: for a fixed $y'\in Q$,
all remaining inequalities hold when $y_n$ is sufficiently large.
Consequently $Q$ is bounded.  Since $Q$ is compact and there are only
finitely many remaining inequalities, a single $T_0$ makes all of
them strictly positive on $Q\times[T_0,\infty)$.  This proves
\eqref{sec7:geo-product-tail}.  It also shows that $y_n$ is bounded
below on $\overline P$; equivalently, every set
$\overline P\cap\{y_n\leq R\}$ is compact.

Let $\nu_D$ be the primitive inward normal of the moving facet.
Boundedness of $P_k$ implies $(\nu_D)_n<0$.  On a sufficiently distant
moving facet, $y_n>T_0$ everywhere, since $y'$ remains in the fixed
bounded polytope $Q$.  Thus the only retained facets meeting it are
vertical.  The common fan of the truncations describes the fixed
compactification $\bar M$; hence every cone containing the ray of $D$
has all its other rays generated by normals with last component zero.

The order of the character $w=Z_n$ along a toric prime divisor of
primitive inward normal $\nu$ is $\langle e_n,\nu\rangle$.
On any toric chart not meeting $D$, all these orders are nonnegative,
so $w$ is holomorphic.  On a chart meeting $D$, its only negative order
is along $D$, and all other divisor orders on that chart vanish.
In local toric coordinates, $1/w$ is therefore a positive integral
power of a defining coordinate of $D$, multiplied by a nonvanishing
holomorphic function.  There is no simultaneous zero and pole giving
an indeterminacy.  These descriptions extend $w$ to a holomorphic map
$\bar M\to\mathbb P^1$ whose fiber over infinity is exactly $D$.
Compactness of $\bar M$ proves properness of its restriction
\eqref{sec7:geo-proper-map}.
\end{proof}

Set $t_k=\mu_k^n$, and, for sufficiently large $R$, write
\[
 \mathcal C_R=\overline P\cap\{y_n\leq R\}.
\]
The next lemma gives both directions of the localization that will
be needed for the maximum principle.

\begin{lem}\label{sec7:geo-localization}
For each fixed $R$, there is a compact set $E_R\Subset M$ such that
\begin{equation}\label{sec7:geo-moment-to-compact}
 \{p\in M:t_k(p)\leq R\}\subset E_R
 \qquad\text{for all sufficiently large }k.
\end{equation}
Conversely, for each compact set $E\Subset M$, there are $R_E$ and
$k_E$ such that
\begin{equation}\label{sec7:geo-compact-to-moment}
 E\subset\{p\in M:t_k(p)\leq R_E\}
 \qquad(k\geq k_E).
\end{equation}
\end{lem}

\begin{proof}
First consider a point of the dense torus with
$y=\mu_k(p)\in\mathcal C_R$.  Since $e_n$ is a recession direction,
both $y$ and $y+e_n$ belong to the fixed compact moment region
$\mathcal C_{R+1}$.  They belong to $P_k$ for all sufficiently large
$k$.  Convexity and Lemma~\ref{bound} give
\begin{equation}\label{sec7:geo-forward-secant}
 x_n=u_{k,n}(y)
 \leq u_k(y+e_n)-u_k(y)\leq C_R.
\end{equation}
Therefore $|w(p)|^2\leq e^{C_R}$.  The moment map and $|w|^2$ are
continuous on $M$, so the same conclusion holds at retained toric
strata.  To include a point with $t_k=R$ directly, approximate it by
dense-torus points with $t_k<R+1$ and use the preceding bound for
$R+1$.  Properness of $w$ now proves
\eqref{sec7:geo-moment-to-compact}, after enlarging $C_R$ if necessary.

For the reverse inclusion, it is enough to bound $t_k$ uniformly on
$\{|w|^2<e^b\}$ for an arbitrary fixed $b$; this open set contains
any prescribed compact set when $b$ is chosen sufficiently large.
Take $b\geq0$ and fix $L>b$.  Theorem~\ref{ck-estimate} bounds
$\varphi_k(Le_n)$ independently of $k$.  The defining Legendre
inequality therefore gives
\begin{equation}\label{sec7:geo-fenchel-line}
 u_k(y',t)\geq Lt-\varphi_k(Le_n)\geq Lt-C_L.
\end{equation}
Also Lemma~\ref{bound} bounds $u_k(y',T_0)$ uniformly for $y'\in Q$.
If $t>T_0$ and $y=(y',t)$ is an interior moment point, the vertical
segment from $(y',T_0)$ to $y$ lies in $P_k$.  Convexity along this
segment yields
\[
 u_{k,n}(y',t)
 \geq\frac{u_k(y',t)-u_k(y',T_0)}{t-T_0}
 \geq\frac{Lt-C_L-C_0}{t-T_0}.
\]
At a dense-torus point with $|w|^2<e^b$, the left side is
$x_n<b$.  Consequently
\[
 (L-b)t\leq C_L+C_0-bT_0.
\]
This bounds $t=t_k$ independently of $k$.  Points with $t\leq T_0$
already satisfy a fixed upper bound.  Finally, density of the torus
and continuity of $t_k$ extend the estimate throughout the open set
$\{|w|^2<e^b\}$, including the zero set of $w$.
This proves \eqref{sec7:geo-compact-to-moment}.
\end{proof}

\subsection{Uniform $C^2$-estimate} Now we begin to do  the uniform $C^2$-estimate on     any $\epsilon$-neighborhood $U_{D,\epsilon}$  of $D$.
Fix a sufficiently large moment level $R$ above.  Choose $k_0$ so large
that $\mathcal C_{R+1}\subset\overline{P_{k_0}}$ and the moving facet
of $P_{k_0}$ is disjoint from $\mathcal C_{R+1}$.  This $k_0$ may
depend on $R$, but is held fixed as the approximation parameter
$k$ tends to infinity.  Let $\omega_0$ be a smooth torus-invariant
K\"ahler metric on $\bar M$ with moment polytope $P_{k_0}$, and let
$\psi_0$ be its convex logarithmic potential.  Write
\[
 u_{\mathrm{ref}}=\psi_0^*,\qquad
 \mu_0=D\psi_0,\qquad v_k=\varphi_k-\psi_0.
\]
The reference dual potential $u_{\mathrm{ref}}$ is continuous on
$\overline{P_{k_0}}$.  Each $v_k$ extends to a smooth function on
$M$.  Indeed, the retained facets have the same primitive normals
and the same support constants $1$ for both metrics.  In a toric
chart along such a face, the same affine logarithmic term converts
$\varphi_k$ and $\psi_0$ into their smooth local K\"ahler potentials.
Subtracting cancels that term.  These extensions agree on their
dense overlaps, so
\[
 \omega_k=\omega_0+\sqrt{-1}\partial\bar\partial v_k
 \qquad\text{on }M.
\]
No uniform boundary regularity is assumed in this individual extension.

\begin{lem}\label{sec7:geo-background-bounds}
On $\{t_k\leq R\}$ one has uniform bounds
\begin{equation}\label{sec7:geo-potential-hamiltonian-bounds}
 |v_k|+|X_kv_k|\leq C_R,
\end{equation}
where
$X_k=\sum_i(A_k)_iZ_i\partial/\partial Z_i$.
There are smooth functions $h_k$ on $M$ satisfying
\begin{equation}\label{sec7:geo-background-equation}
 \omega_k^n=e^{h_k-v_k-X_kv_k}\omega_0^n.
\end{equation}
For every fixed compact subset of $M$, the functions $h_k$ and
the vector fields $X_k$ have uniform bounds for all derivatives.
In particular,
\begin{equation}\label{sec7:geo-volume-ratio}
 C_R^{-1}\leq\frac{\omega_k^n}{\omega_0^n}\leq C_R
 \qquad\text{on }\{t_k\leq R\}.
\end{equation}
\end{lem}

\begin{proof}
At a dense-torus point with $t_k\leq R$, the maximizer in the
Legendre formula for $\varphi_k$ lies in
$\mathcal C_R\subset\overline{P_{k_0}}$.  Thus
\[
 \varphi_k(x)=\sup_{y\in\overline{P_{k_0}}}
                      (x\cdot y-u_k(y)),
 \qquad
 \psi_0(x)=\sup_{y\in\overline{P_{k_0}}}
                      (x\cdot y-u_{\mathrm{ref}}(y)).
\]
Consequently,
\[
 |v_k(x)|\leq
 \sup_{\overline{P_{k_0}}}|u_k-u_{\mathrm{ref}}|\leq C_R
\]
by Lemma~\ref{bound}.  The identity
\[
 X_kv_k=\langle A_k,\mu_k-\mu_0\rangle
\]
bounds the other term because $A_k$ is bounded, $\mu_k$ lies in
$\mathcal C_R$, and $\mu_0$ lies in the fixed compact polytope
$\overline{P_{k_0}}$.  Both bounds extend to the retained strata by
continuity (or by using a slightly larger fixed level $R$).

On the dense torus define
\[
 h_{\mathrm{ref}}=-\psi_0-\log\det D^2\psi_0,
 \qquad
 h_k=h_{\mathrm{ref}}-m_k-\langle A_k,\mu_0\rangle.
\]
The function $h_{\mathrm{ref}}$ extends smoothly to $M$.
To see this, choose integral toric coordinates along a retained
face.  If its normal logarithmic variables are $\tau_1,\ldots,\tau_r$,
then $\psi_0+\sum_{i=1}^r\tau_i$ is a smooth local potential.
The logarithm of its complex Hessian determinant differs from
$\log\det D^2\psi_0$ by $-\sum_{i=1}^r\tau_i$.
These logarithmic terms cancel in $h_{\mathrm{ref}}$; any additional
nonvanishing coordinate Jacobian contributes only a fixed smooth
term.  This uses exactly the retained support constants $1$.

Equation \eqref{RMA-2} now gives
\eqref{sec7:geo-background-equation}, since the logarithmic
coordinate Jacobians cancel in the ratio of the two volume forms.
The reference data are fixed, while $A_k$ and $m_k$ are bounded,
so $h_k$ and $X_k$ have uniform smooth bounds on every fixed compact
set.  By \eqref{sec7:geo-moment-to-compact}, the sets $\{t_k\leq R\}$
lie in one such compact set.  Combining this observation with
\eqref{sec7:geo-potential-hamiltonian-bounds} proves
\eqref{sec7:geo-volume-ratio}.
\end{proof}


We use the compactness and potential bounds established above to estimate
the original smooth approximating metrics directly. Write
\[
 t_k=\frac{\partial\varphi_k}{\partial x_n}=\mu_k^n,
 \qquad \Omega_{k,R}=\{p\in M:t_k(p)\le R\}.
\]
Fix two sufficiently large levels $R_1<R_2$. Apply the preceding
subsection with $R=R_2$. The preceding comparison
between geometric and moment regions places $\Omega_{k,R_2}$ in a fixed
compact subset of $M$, for all sufficiently large $k$. On that subset
choose the fixed background $\omega_0$ and potential $\psi_0$ used above,
and put $v_k=\varphi_k-\psi_0$. Although the logarithmic torus potentials
need not extend individually, $v_k$ is a smooth function on $M$ and
\[
 \omega_k=\omega_0+\sqrt{-1}\partial\bar\partial v_k.
\]
Here each $\omega_k$ is smooth and positive on $M$ by
Theorem~\ref{partial-smooth}. No uniform boundary estimate is being
assumed at this stage.

Let
\[
 X_k=\sum_{j=1}^n(A_k)_j Z_j\frac{\partial}{\partial Z_j},
 \qquad \mathcal L_k=\Delta_{\omega_k}+X_k,
 \qquad \Delta_{\omega_k}=g_k^{i\bar j}\partial_i\partial_{\bar j}.
\]
The vector field $X_k$ extends holomorphically to $M$. On torus-invariant
real functions its action is real, so $\mathcal L_k$ has the usual
maximum principle on the functions considered below. The relative
Monge--Amp\`ere equation is
\begin{equation}\label{sec7:yau-relative-equation}
 \mathcal F_k:=\log\frac{\omega_k^n}{\omega_0^n}
       =h_k-v_k-X_kv_k.
\end{equation}
The previously obtained bounds give, on $\Omega_{k,R_2}$,
\begin{equation}\label{sec7:yau-input-bounds}
 |v_k|+|X_kv_k|\le C,
 \qquad C^{-1}\le e^{\mathcal F_k}\le C.
\end{equation}
The functions $h_k$, the background metric, and the coefficients of
$X_k$ have uniformly bounded derivatives on the fixed geometric compact
set containing these regions. On the dense orbit one may write
\[
 h_k=-\log\det D^2\psi_0-\psi_0
       -\langle A_k,D\psi_0\rangle-m_k.
\]
Its smooth extension and its uniform bounds here refer only to a fixed
region away from $D$.

As Lemma \ref{C^2}, we have
\begin{lem}\label{sec7:yau-drift-trace}
Let
\[
 T_k=\operatorname{tr}_{\omega_0}\omega_k,
 \qquad S_k=\operatorname{tr}_{\omega_k}\omega_0.
\]
On $\Omega_{k,R_2}$ there is a constant independent of $k$ such that
\begin{equation}\label{sec7:yau-trace-inequality}
 \mathcal L_k\log T_k\ge-CS_k-C.
\end{equation}
\end{lem}

\begin{proof}
Suppress the subscript $k$ and write $\Delta_0=\Delta_{\omega_0}$. From (\ref{Yau's C^2}), we have
\begin{equation}\label{sec7:yau-basic-trace}
 \Delta_{\omega}\log T
 \ge\frac{\Delta_0\mathcal F}{T}-CS-\frac CT.
\end{equation}
The constants in this inequality depend only on the background geometry
on the fixed compact set; see \cite{Yau} and
\cite[Appendix~B, (B.3)]{BBEGZ}. We record the commutator with
$X$ explicitly in order to control the drift without assuming any
higher derivatives of $v$. In holomorphic coordinates, write
$X=X^a\partial_a$. Since $\partial_{\bar j}X^a=0$,
\begin{equation}\label{sec7:yau-holomorphic-commutator}
 \Delta_0(Xv)-X(\Delta_0v)
 =g_0^{i\bar j}(\partial_iX^a)v_{a\bar j}
     -X(g_0^{i\bar j})v_{i\bar j}.
\end{equation}
There is no term involving an uncontrolled first or third derivative of
$v$ in this difference. Since $(v_{i\bar j})=(g_{i\bar j})-
(g_{0,i\bar j})$, positivity of $\omega$ implies
\[
 |\Delta_0(Xv)-X(\Delta_0v)|\le C(T+1).
\]
Indeed, in an $\omega_0$-unitary frame diagonalizing $\omega$, the norm
of this Hermitian tensor is bounded by $T+n$; all the coefficients in
\eqref{sec7:yau-holomorphic-commutator} are uniformly bounded.

Using $\Delta_0v=T-n$ in \eqref{sec7:yau-relative-equation} now gives
\[
 \Delta_0\mathcal F
 =\Delta_0h-(T-n)-X(T-n)
   -\bigl(\Delta_0(Xv)-X(\Delta_0v)\bigr)
 \ge-XT-C(T+1).
\]
Substitution in \eqref{sec7:yau-basic-trace} yields
\[
 (\Delta_{\omega}+X)\log T
 \ge-CS-C-C/T.
\]
The determinant lower bound in \eqref{sec7:yau-input-bounds} implies
$T\ge n e^{\mathcal F/n}\ge c>0$. This proves
\eqref{sec7:yau-trace-inequality}.
\end{proof}

The identities needed for the cutoff also contain only quantities
already controlled. On the dense orbit, differentiating \eqref{RMA-2}
in $x_n$ gives
\begin{align}\label{sec7:yau-moment-identities}
 \mathcal L_k t_k&=-t_k,\notag\\
 |\partial t_k|_{\omega_k}^2
 &=\sum_{i,j}\varphi_k^{ij}\varphi_{k,ni}\varphi_{k,nj}
   =\varphi_{k,nn}
   =u_k^{nn}\bigl(D\varphi_k\bigr)\le C_{R_2}.
\end{align}
The last bound follows from \eqref{sec7:dual-lower}, since
$\mathcal C_{R_2}$ is contained in a fixed truncation $K_R$.
Each quantity in
these identities extends smoothly across the retained divisors for a
fixed $k$, and the inequalities therefore hold there by continuity.

Choose a smooth function $\chi:\mathbb R\to[0,1]$ with
$\chi=1$ on $(-\infty,R_1]$, $\chi=0$ on $[R_2,\infty)$, and
$\chi>0$ on $(-\infty,R_2)$. Write $\chi_k=\chi(t_k)$.
Its first two derivatives are bounded by constants depending only on
the chosen levels. The moment coordinates on $\Omega_{k,R_2}$ are
uniformly bounded, so \eqref{sec7:yau-moment-identities} gives
\begin{equation}\label{sec7:yau-cutoff}
 \mathcal L_k\log\chi_k
 =-\frac{\chi'(t_k)}{\chi_k}t_k
   +\left(\frac{\chi''(t_k)}{\chi_k}
          -\frac{\chi'(t_k)^2}{\chi_k^2}\right)
       |\partial t_k|_{\omega_k}^2
 \ge-C\chi_k^{-2}
\end{equation}
where $\chi_k>0$.

\begin{prop}\label{sec7:yau-local-equivalence}
For every fixed compact set $K\subset U_{D,\epsilon}$,  there is a constant $C_K$
such that the approximating metrics satisfy
\begin{equation}\label{sec7:yau-metric-equivalence}
 C_K^{-1}\omega_0\le\omega_k\le C_K\omega_0
 \qquad\hbox{on }K
\end{equation}
for all sufficiently large $k$. The background in this statement is
fixed once a larger compact neighborhood of $K$ has been chosen.
\end{prop}

\begin{proof}
First establish the estimate on $\Omega_{k,R_1}$. Since
\[
 \mathcal L_kv_k=n-S_k+X_kv_k,
\]
Lemma~\ref{sec7:yau-drift-trace} and
\eqref{sec7:yau-input-bounds} allow us to choose a fixed $A>0$ such that
\[
 \mathcal L_k(\log T_k-Av_k)\ge S_k-C.
\]
Consider
\[
 G_k=\log T_k-Av_k+2n\log\chi_k
\]
on $\{\chi_k>0\}$. The closure of this domain is contained in the fixed
geometric compact set described above. For every fixed $k$, both metrics and
$v_k$ are smooth there, including at all retained divisors. Thus
$G_k\to-\infty$ as $\chi_k\to0$, and its maximum is attained at a
point of $M$ where $\chi_k>0$. In particular, a point on a retained
divisor is an interior point for this maximum principle.

At a maximum point, \eqref{sec7:yau-cutoff} gives
\[
 0\ge\mathcal L_kG_k\ge S_k-C-C\chi_k^{-2},
 \qquad S_k\le C\chi_k^{-2}.
\]
For positive eigenvalues of $\omega_k$ relative to $\omega_0$,
the elementary trace--determinant inequality gives
\[
 T_k\le C_n e^{\mathcal F_k}S_k^{n-1}.
\]
When $n=1$ this means simply $T_k=e^{\mathcal F_k}$. Consequently, at the
maximum point,
\[
 T_k\le C\chi_k^{-2(n-1)},\qquad
 G_k\le C-Av_k+2\log\chi_k\le C.
\]
On $\Omega_{k,R_1}$ we have $\chi_k=1$, so the uniform bound for
$v_k$ gives $T_k\le C$. The determinant lower bound then bounds
each eigenvalue away from zero and proves the two-sided metric
estimate on that region.

The preceding geometric--moment comparison also shows that each fixed
compact $K\subset M$ is contained in $\Omega_{k,R_1}$ after choosing
$R_1$ sufficiently large, uniformly for all large $k$. The estimate
therefore holds on $K$, as claimed.
\end{proof}

\begin{proof}[Proof of Theorem \ref{smoothness}]
By Proposition \ref{sec7:yau-local-equivalence}, from the   Evans--Krylov estimate  and bootstraping,  we get
\begin{equation}\label{sec7:boot-all-orders}
 \|v_k\|_{C^N(K)}\le C_{K,N}
 \qquad\hbox{for every compact }K\subset M\hbox{ and every }N\ge0.
\end{equation}

By Arzel\`a--Ascoli and a diagonal argument, every subsequence has a
further subsequence converging smoothly on compact subsets of $M$.
On the dense orbit its limit is already determined by the convergence
$\varphi_k\to\varphi_\infty$ in Section~6. Smooth limits agree on the
dense orbit and hence agree everywhere. Thus the whole sequence under
consideration has the same smooth limit. This proves
\eqref{sec7:smooth-metric-convergence}.
The lower bound in \eqref{sec7:yau-metric-equivalence} shows that the
limit is a positive K\"ahler metric. On the dense orbit it is
$\sqrt{-1}\partial\bar\partial\varphi_\infty$, so this constructs
its smooth extension across all retained divisors and their
intersections at the same time.

The moment maps also converge smoothly in these charts: their
components are the corresponding fixed background moment components
plus the action of the torus holomorphic generators on $v_k$.
Denote the resulting smooth moment map by $\mu_\infty$, and set
\[
 \theta_\infty=\langle A_0,\mu_\infty\rangle,
 \qquad f=-\theta_\infty,
 \qquad X_0=\sum_j(A_0)_jZ_j\frac{\partial}{\partial Z_j}.
\]
Taking $\sqrt{-1}\partial\bar\partial$ of \eqref{RMA-limit} on the
dense orbit gives
\begin{equation}\label{sec7:extended-soliton}
 \operatorname{Ric}(\omega_\infty)-\omega_\infty
 =\sqrt{-1}\partial\bar\partial\theta_\infty
 =L_{X_0}\omega_\infty.
\end{equation}
The additive constant $m_0$ and the logarithmic Jacobian terms have
zero $\partial\bar\partial$ there. All terms in
\eqref{sec7:extended-soliton} are now smooth on $M$, so the identity
extends to the whole manifold. With the displayed convention for
logarithmic coordinates, the Legendre identities give
$\nabla^{1,0}f=-X_0$. Thus
\[
 \operatorname{Ric}(\omega_\infty)
       +\sqrt{-1}\partial\bar\partial f=\omega_\infty,
\]
and the limiting metric is a gradient KR soliton. To prove the completeness, we will improve (\ref{sec7:yau-moment-identities}) to \begin{equation}\label{sec7:global-height-gradient}
 |\partial t_k|_{\omega_k}^2\le a(t_k+a)\qquad\text{on }M.
\end{equation} 
 We will prove this using maximal principle (see also \cite[Proposition~4.4]{CHTX}).  Fix $a>0$ such that
$y_n+a\ge1$ on $\overline P$, hence also on every $P_k$.

Fix $k$ and write $\varphi=\varphi_k$, $t=\varphi_n$, $B=\varphi_{nn}$,
and $L=\varphi^{ij}\partial_i\partial_j+A_k\cdot D$ on $\mathbb R^n$.
The identities in \eqref{sec7:yau-moment-identities} give
\[
 Lt=-t,\qquad \varphi^{ij}t_i t_j=B.
\]
Set $H=D^2\varphi$ and $S=\partial_n H$. Differentiating the logarithm
of \eqref{RMA-2} twice in the $x_n$ direction gives
\[
 LB=-B+\operatorname{tr}\bigl((H^{-1}S)^2\bigr).
\]
Since $DB=Se_n$ and $B=e_n^T H e_n$, the Cauchy--Schwarz inequality
gives $(DB)^T H^{-1}DB\le B\,\operatorname{tr}((H^{-1}S)^2)$.
Consequently $L\log B\ge-1$, and hence
\[
 L\log\frac{B}{t+a}
 \ge -1+\frac{t}{t+a}+\frac{B}{(t+a)^2}
 =\frac{1}{t+a}\left(\frac{B}{t+a}-a\right).
\]

For this fixed $k$, $B$ is bounded on $\mathbb R^n$. Indeed, at any
$x_0$ subtract the tangent plane to $\varphi$ at $x_0$ and denote the
difference by $h$. The section $\{h<1\}$ is bounded because
$D\varphi(x_0)\in P_k^\circ$. Apply Lemma~\ref{c2-more} to $h-1$ on
this section, as in the proof of Proposition~\ref{sec7:dual-lower-prop}.
Here $-1\le h-1\le0$, $|Dh|\le\operatorname{diam}P_k$, and the affine
term in the equation has gradient of norm $|D\varphi(x_0)|$, bounded
independently of $x_0$ since $P_k$ is bounded. At $x_0$ we have $h=0$,
so the lemma gives $B(x_0)\le C_k$.

By \eqref{sec7:linear-growth}, $\varphi(x)\to+\infty$ as
$|x|\to\infty$. Thus for every $\varepsilon>0$ the function
\[
 G_\varepsilon=\log\frac{B}{t+a}-\varepsilon\varphi
\]
attains its maximum at some $x_\varepsilon\in\mathbb R^n$. Since
$L\varphi=n+\langle A_k,D\varphi\rangle$, the preceding differential
inequality and $LG_\varepsilon(x_\varepsilon)\le0$ give, at this point,
\[
 \frac{B}{t+a}\le
 a+\varepsilon(t+a)\bigl(n+\langle A_k,D\varphi\rangle\bigr)
 \le a+\varepsilon C_k,
\]
where $C_k$ is finite because $D\varphi\in P_k$. The maximum property
and the normalization $\varphi\ge0$ now yield, for every fixed $x$,
\[
 \frac{B(x)}{t(x)+a}
 \le (a+\varepsilon C_k)e^{\varepsilon\varphi(x)}.
\]
Letting $\varepsilon\downarrow0$ with $k$ fixed gives $B\le a(t+a)$.
The identity $B=|\partial t_k|_{\omega_k}^2$ on the dense orbit and
smoothness across the retained divisors prove
\eqref{sec7:global-height-gradient} with the same $a$ for all $k$.

Let $t_\infty=\mu_\infty^n$. By local smooth convergence of the metrics
and moment maps, \eqref{sec7:global-height-gradient} gives
\[
 \left|\partial\sqrt{t_\infty+a}\right|_{\omega_\infty}^2
 =\frac{|\partial t_\infty|_{\omega_\infty}^2}{4(t_\infty+a)}
 \le \frac a4.
\] The completeness follows easily.

\end{proof}

\section{Complete KR solitons on  $\mathbf G$-manifolds}

In this section, we generalize Theorem \ref{main-resullt} to $\mathbf G$-manifolds.  Unlike toric manifolds,  there is no any KR soliton on a $\mathbf G$-manifold, otherwise it  should satisfy some obstruction condition \cite{Del, LZZ, LTZ}.
Thus there  need more condition to guarantee the existence of complete KR solitons on  $\mathbf G$-manifolds.

Given   a complex reductive Lie group  $\mathbf G$, i.e., a complexification of compact group $\mathbf K$.  By a $\mathbf G$-manifold $(M, J)$, we mean a {\it (bi-equivariant) compactification of $\mathbf G$} which admits a holomorphic $\mathbf G\times \mathbf G$-action and has an open and dense orbit isomorphic to $\mathbf G$ as a $\mathbf G\times \mathbf G$-homogeneous space (cf. \cite{Del1, LZZ}).  $(M, L)$ is called a {\it polarized compactification} of $\mathbf G$   if  $L$ is a $\mathbf G\times \mathbf G$-linearized ample line bundle on $M$ \cite{AK}.
 A toric manifold is a special case of $\mathbf G$-manifold when $\mathbf G$ is  torus.

By \cite{AB1, AB2},  for the $\mathbf G$-manifold $(\bar M, J)$ there is an $n$-dimensional toric complex submanifold $(Z,J|_Z)$ associated to a maximal torus $T^\mathbb C$ of $\mathbf G$.
Let $\bar P$ be a moment polytope of $(Z,J|_Z)$ associated to  a {\it polarized compactification}  $(\bar M, L)$.    Assume that $\Phi$ is the root system of $(\mathbf G, \mathbf T^\mathbb C)$ in $\mathfrak  a^*$ and $\Phi_+=\{\alpha_1,..., \alpha_d\} $ is a chosen set of positive roots.  
Denote  $\bar P_+$  to be   the part of $\bar P$ defined by $\Phi_+$, which 
is a quotient space of $\bar P$ by the  Weyl  subgroup $W$ of $\mathbf T^\mathbb C$.   We call  $P$  wonderful  if  $P_i\cap  (\partial \bar P\setminus  (\partial \bar P)^o)=\emptyset$ for any $i=1,...d$,  where  $P_i=\{y\in \mathfrak a^* |~ \langle\alpha_i, y\rangle_{killing}=0\}$ is a hyperplane in $\mathfrak   a^*$.  In this section, we always assume $\bar P$  is wonderful.

 Let $\{F_A\}_{A=0,1,...,d_0}$ be the facets of $\bar P$ and $\{F_A\}_{A=0,1,...,d_+}$ be those whose relative interior intersects $\mathfrak a_+^*$. Suppose that
\begin{align}\label{polytope-eq}
\bar P=\cap_{A=0}^{d_0}\{l_A=\lambda_{A}-u_A(y)>0\}
\end{align}
for some primitive vector $u_A\in\mathfrak N$ and the facet
$$F_A\subseteq\{l_{A}=0\},A=0,1,...,d_0.$$
By  the $W$-invariance, for each $A\in\{1,...,d_0\}$, there is some $w_A\in W$ such that $w_A(F_A)\in\{F_B\}_{B=0,1,...,d_+}$, where $F_B$ is a  facet of $\bar P_+$. Moreover, in the case of Fano $\mathbf G$-manifold, each number $\lambda_{A}$ 
satisfies the following quantity (cf. \cite{LTZ}),
\begin{align}\label{relation-lambda}\lambda_{A}\,=\,1+2\rho_A(u_A),
\end{align}
where $\rho_A\,=\,w_A^{-1}(\rho)$ and  $\rho$ is given by
\begin{align}\label{rho}
\rho\,=\,\frac12\sum_{\alpha_i\in\Phi_+}\alpha_i.
\end{align}
Note that $\rho_A(u_A)$ is independent of the choice of $w_A\in W$ and hence it is well-defined.

Let $D$ be a smooth semi-ample $\mathbf G\times \mathbf G$-invariant prime divisor.
We assume that its restriction to $Z$ is the reduced union of the toric
divisors indexed by a $W$-orbit $S$.  For a real parameter $k\geq1$, let
\[
 \Omega_k=2\pi\bigl(c_1(\bar M)+(k-1)c_1(\mathcal O(D))\bigr).
\]
The compatible moment polytope is
\begin{equation}\label{pk-u}
 \begin{split}
 l_{A,k}(y)&=\lambda_A
                  +(k-1)\mathbf1_S(A)-u_A(y),\\
 P_k&=\bigcap_{A=0}^{d_0}\{l_{A,k}\geq0\},
 \qquad P_k^+=P_k\cap\overline{\mathfrak a_+^*}.
 \end{split}
\end{equation}
The entire orbit $S$ moves together, so $P_k$ is $W$-invariant.

Set
\begin{equation}\label{g:noncompact-polyhedron}
 P_\infty
 =\bigcap_{A\notin S}
   \{y\in\mathfrak a^*:u_A(y)\leq 1+2\rho_A(u_A)\},
 \qquad
 P_\infty^+=P_\infty\cap\overline{\mathfrak a_+^*}.
\end{equation}
The family $P_k$ is increasing and exhausts a noncompact polytope $P_\infty$.  The latter is
$W$-invariant and contains the origin in its interior.  For simplicity, we call $P_\infty$ toric moment polytope associated to $\mathbf G$-manifold $\bar M\setminus D$.
   Its recession cone is
\[
 C_\infty:=\operatorname{rec}(P_\infty)
 =\{d\in\mathfrak a^*:u_A(d)\leq0\text{ for every }A\notin S\}.
\]

By modifying  Definition \ref{cone-polytope},   we introduce

\begin{defi}\label{Fano-cone-G}
The polyhedron $P_\infty$ is called a Fano cone  if  all $P_k$ in (\ref{pk-u}) satisfy  the Delzant  condition  in addition that $C_\infty$ is strictly convex, i.e., 
\begin{equation}\label{cone-condition}
C_\infty\cap(-C_\infty)=\{0\}.
\end{equation}
\end{defi}

\subsection{Singular KR solitons} In this subsection, we  generalize Theorem \ref{existence-1} to $\mathbf G$-manifold.
 
First we determine the KR soliton VF $X$ in  K\"ahler class $(\bar M, \Omega_k)$ associated to $\bar P$.  Suppose that the center $\mathbf T_0\subset 
\mathbf T^\mathbb C$  of $\mathbf G$  is not trivial.  Let $\{E_i\}$ is a basis of Lie algebra of $\mathbf T_0$. 
Following the $H$-invariant construction in \cite{TZ-survey} and the weighted moment formula in \cite[Section~5.1]{LZZ},    
         for any torus HVF $v=(v_i)\in \mathbb R^r$ induced by $ \sum_{i=1}^r a_iE_i$, where $a=(a_i)$ is a  real vector in $\mathbb R^r$ , we introduce a restricted H-invariant for $v$ 
          by
            \begin{align}\label{G-case-1}H(v)=\int_{ P_k^+} e^{<v, y>} \pi(y) dy,
  \end{align}
 where $<v, y>=\sum_{i=1}^r v_iy_i$ and
\begin{align}\label{p-function} \pi(y) =
  \prod_{\alpha_i\in\Phi_+}\langle\alpha_i,y\rangle^2.
  \end{align}
By the convexity of $H(\cdot)$, there exist a unique critical  $v_0$ which attains the minimum of $H(\cdot)$.  Without confusion,   as the case of toric manifolds in Section 3, we call   the corresponding HVF  $X$ of $v_0$ a KR soliton VF in  K\"ahler class $(\bar M, \Omega_k)$. 
  
  Next we shall construct an admissible  K\"ahler metric $\omega_k$ associated to $D$  in  any K\"ahler class $(\bar M, \Omega_k)$ as in the case of toric manifolds via the Gulliemin metric.  Actually,
  \[
  u_k=\sum_A c_{A,k}l_{A,k}\log l_{A,k},\qquad
  c_{A,k}=\begin{cases}1,&A\notin S,\\1/k,&A\in S.\end{cases}
  \]
  is a $W$-invariant convex function on $\mathfrak   a^*$.
  Then  its Legendre function $\varphi_{k}$  is also  $W$-invariant, where
\begin{align}\label{ Legendre-function} \varphi_{k}\,=\,\sup_{y\in P_k} (\langle x,y\rangle- u_{k}(y)), ~\forall~x\in \mathfrak{a}.
\end{align}
Thus, by \cite[Theorem 2]{Ab-JDG} and \cite{AL} (also see \cite{LZZ, LTZ}),  we can extend
$$\omega_{k}\,=\, \sqrt{-1}\partial\bar\partial \varphi_{k},~{\rm on}~\mathfrak a^\mathbb C, $$
to a $K\times K$-invariant metric on $\bar M\setminus D$. Since  the Ricci potential $h_0$ of  $\omega_{k}$ is given by (cf.  \cite{LZZ, LTZ}), 
\begin{align}\label{g:h-potential}
h_0&=-\log\det(\varphi_{k,ij})-\varphi_{k}+\log \mathbf J(x)-\log\prod_{\alpha_i\in\Phi_+}\alpha_i^2(\nabla\varphi_{k})\\
&=\log\det(u_{{k},ij})-y_iu_{{k},i}+u_{k}+\log \mathbf J(\nabla u_{k})-\log\pi(y),
\end{align}
where 
\begin{align}\label{j-convex}
\mathbf J(x)=\prod_{\alpha\in\Phi_+}
             \sinh^2\!\left(\frac{\alpha(x)}2\right),
\qquad x\in\mathfrak a_+.
\end{align}
Here we use the Cartan parametrization $g=k_1\exp(x/2)k_2$,
so the Legendre gradient image is $P_k$.  In the convention
$g=k_1\exp(x)k_2$, the gradient image is $2P_k$ and the root
Jacobian is $\prod_{\alpha>0}\sinh^2\alpha(x)$.
The latter root Jacobian also occurs in the proof of Theorem~\ref{g-para}.

Write $\lambda_{A,k}=\lambda_A+(k-1)\mathbf1_S(A)$.
Since $\sum_Ac_{A,k}u_A$ is $W$-invariant and hence central,
direct differentiation gives the following identity on the positive
chamber, up to an additive normalization constant:
\begin{align}\label{h0-sum}
h_0={}&\log\det\left(\sum_A
             \frac{c_{A,k}u_A\otimes u_A}{l_{A,k}}\right)
 +\sum_A c_{A,k}\lambda_{A,k}\log l_{A,k}
 +\sum_A c_{A,k}u_A(y)\notag\\
&+2\sum_{\alpha\in\Phi_+}
 \log\frac{\sinh\!\left(-\frac12\sum_A
                     c_{A,k}\alpha(u_A)\log l_{A,k}\right)}
               {\alpha(y)}.
\end{align}
The quotients along a Weyl wall are interpreted by continuity.
For each fixed $k$, the determinant and root factors have the
boundary cancellations used in \cite[Proposition 5.1]{LTZ};
the smoothness and Delzant hypotheses give a bounded Ricci potential
on $\mathbf G$.  The bound here is allowed to depend on $k$.

 By the first relation in (\ref{g:h-potential}),   $\varphi_k$ satisfies the following MA equation, 
 $$\pi (\nabla\varphi){\rm det}(\nabla^2 \varphi)= J(x) e^{-h-\varphi -<v_0, \nabla\varphi>},
 $$
 where $h=h_0-<v_0, \nabla\varphi>.$
 Thus,  $h$ becomes zero when    $\omega=\sqrt{-1}\partial\bar\partial\varphi$ is a singular  KR soliton  in $\Omega_k $   on $\bar M$ induced by a $K\times K$-invariant function $\varphi$ on $\mathbf G$. Namely, 
 $\varphi$ satisfies a  reduced KR soliton equation on $\mathfrak a_+$,
 \begin{align}\label{KR-soliton-G}\pi (\nabla\varphi){\rm det}(\nabla^2 \varphi)= J(x) e^{-\varphi -<v_0, \nabla\varphi>}.
 \end{align}
 
 As in the proof of Theorem \ref{existence-1}  for toric manifolds, we use the variation method to solve (\ref{KR-soliton-G}). Actually,  similar to the canonical class on Fano $\mathbf G$-manifolds
 (cf. \cite[Section~5]{LZZ} and, for the reduced Ding functional, \cite[Section~6]{LTZ}), we introduce a reduced  Ding-energy on $P_k^+$ by 
  $$\mathcal D^X(u)=\mathcal N^X(u)+\mathcal L^X(u), $$
  where $X$ is the holomorphic VF associated to $v_0$, and 
\begin{align} 
\mathcal L^X(u)= \int_{P_k^+}\langle y-2\rho,\nabla u\rangle e^{<v_0, y>  }\pi\,dy, \nonumber
\end{align}
\begin{align}
\mathcal N^X(u)=-\int_{P_k^+}\left(\log\det\left(u_{ij}\right)+ \ln J(\nabla u)-\langle 2\rho,\nabla u\rangle\right)e^{ <v_0, y>}\pi\,dy. \label{MN}
\end{align}
Here $u$ is chosen in the following space, 
\begin{align}\mathcal E^1_{K\times K}(P_k)\,=\,\{u| ~&\text{$u$ is convex, $W$-invariant on $P_k$ which satisfies $u(o)=\inf_{P_k} u=0$
  and}
\notag\\
&\int_{P_k^+} u \pi\,dy<+\infty\}.\notag
\end{align}
Then by extending \cite[Theorem 5.1]{LZZ} to the K\"ahler class $\Omega_k$ with the help of background $\omega_k$  constructed above, we have

\begin{prop}\label{thm40}
Let
\begin{align}\label{baracenter}{\rm bar}_X( P_k^+)={\frac{\int_{ P_k^+}y e^{\theta_X(y)}\pi(y)\,dy}{\int_{ P_k^+} e^{\theta_X(y)}\pi(y) \,dy}}
\end{align}
be  a  restricted barycenter  on  $ P_k^+$,
where $\theta_X(y)=<v_0, y>. $ 
Suppose 
 \begin{align}\label{bary}
{\rm bar}_X( P_k^+)\,\in\,  2\rho+\Xi.
\end{align}
Then $D^X(\cdot)$ is proper on $\mathcal E^1_{K\times K}(P_k)$.
\end{prop}
  
 The following is a generalization of Theorem \ref{existence-1}  to   $\mathbf G$-manifolds. 

\begin{theo}\label{existence-G}  Let   $(\bar M, L)$ be a  polarized $\mathbf G$-manifold with  wonderful   moment polytope $\bar P$ given by (\ref{polytope-eq}).  Then on  any  K\"ahler class $(\bar M, \Omega_k)$ ($k\ge \lambda_0$) 
 there exists a singular KR soliton $\omega_{\phi_k}\in \Omega_k$, if (\ref{bary}) is satisfied.
 
\end{theo}

\begin{proof}By Proposition \ref{thm40},  there exists a minimizer $u$ of $D^X(\cdot)$ in $\mathcal E^1_{K\times K}(P_k)$, whose $K\times K$-invariant  Legendre function $\varphi_u$   satisfies  (\ref{KR-soliton-G})  almost everywhere  on $\mathfrak a_+$ (cf. \cite{LTZ}).  Choose a $K\times K$-invariant convex function $\psi_0$ on $\mathfrak a$  such that the extension of 
$\omega_0=\sqrt{-1}\partial\bar\partial \psi_0$ is a smooth K\"ahler metric in $\Omega_k$. Then    $\phi=\varphi_u-\psi_0$ can be extended to a K\"ahler potential  in  $\mathcal E_{K\times K}(\bar M, \Omega_k)$ associated to $\omega_0$ (cf. \cite[Theorem 4.2]{LTZ}). Thus by \cite[Lemma~4.5]{LTZ},  $ \omega_{\phi}^n={\rm det}(\nabla^2 \varphi)$.   As a consequence,    $\phi$ is a weak solution of following CMA equation,
$$\omega_{\phi}^n= e^{h-\phi-X(\phi)} \omega_0^n |s|^{2(k-1)},$$
where $h$ is a uniform bounded function on $\bar M$ and $s$ is a defining holomorphic section of $D$.  Thus by the Kolodziej's result, $\phi$ is uniformly bounded.  This proves that
$\sqrt{-1}\partial\bar\partial \phi_u$ can be extended as a singular KR soliton in $ \Omega_k$.

\end{proof}

It seems the condition  (\ref{bary}) is also necessary as same as in the case of $\mathbb Q$-Fano varieties \cite{LTZ, TZ24}. In 
Subsection~\ref{sec9:sl2-blowup} and \ref{sec9:gl2}, we will  give  some lower dimensional 
$\mathbf G$-manifolds,
for which \eqref{bary} holds for every $k\ge1$.  Thus there always  exist a  singular KR soliton in $\Omega_k$ for any integer $k\ge 1$ on those $GL(2,\mathbb C)$-manifolds.

With respect to the existence of complete KR solitons on $\bar M\setminus D$, we propose   the following conjecture.

\begin{conj}\label{conj-g}
	Let   $(\bar M, L)$ be a  polarized $\mathbf G$-manifold with  wonderful   moment polytope $\bar P$ given by (\ref{polytope-eq}). 
	Suppose  that $P_\infty$ is  a Fano cone  associated to   $\bar M\setminus  D$ and there exists a large integer $k$ such that 
	(\ref{bary}) is satisfied.
	Then there exists a complete GSKRS on $\bar M\setminus D$ with moment polyhedron $P_\infty$. 
\end{conj}

In the following,  we  verify the above conjecture on $\bar M\setminus D$  with   parabolic Fano cone  $P_\infty$ as a generalization of  Theorem \ref{main-resullt}  for toric manifolds.

\subsection{ $\mathbf G$-manifolds  $\bar M\setminus D$  with parabolic  Fano cone}


\begin{defi}\label{Fano-cone-G-parabolic}
We call $P_\infty$ \emph{parabolic} if $C_\infty$ is a single nonzero ray.
\end{defi}

The central decay condition is equivalent to $C_\infty$ containing no
line.  Indeed, the condition excludes both $d$ and $-d$ from $C_\infty$
when $d\neq0$.  Conversely, a pointed polyhedral cone admits a linear
functional strictly negative on its nonzero elements.  Averaging this
functional over $W$ preserves strict negativity, since $C_\infty$ is
$W$-invariant, and produces an element of $\mathfrak a^W=\mathfrak z$.

For $v\in\mathfrak z$ define the extended-valued function
\begin{equation}\label{g:limiting-H}
 H_\infty(v)
 =\int_{P_\infty^+}e^{\langle v,y\rangle}\pi(y)\,dy
 \in(0,+\infty].
\end{equation}
All Lebesgue measures below use a fixed $W$-invariant Euclidean structure
compatible with the central and semisimple decomposition.

\begin{lem}\label{g:limiting-central-vector}
 $H_\infty(v)<\infty$ if and only if $v\in\mathcal V(P_\infty)$.  On this open set $H_\infty$ is smooth
and strictly convex and has a unique minimizer $v_\infty$.  It satisfies
\begin{equation}\label{g:limiting-central-moment}
 \int_{P_\infty^+}\operatorname{pr}_{\mathfrak z^*}(y)
        e^{\langle v_\infty,y\rangle}\pi(y)\,dy=0.
\end{equation}
Moreover, for every compact subset $K\Subset\mathcal V(P_\infty)$ and
every nonnegative integer $m$,
\begin{equation}\label{g:uniform-weighted-moments}
 \sup_{v\in K}\int_{P_\infty^+}(1+|y|)^m
       e^{\langle v,y\rangle}\pi(y)\,dy<\infty.
\end{equation}
\end{lem}

\begin{proof}
The Minkowski--Weyl decomposition writes $P_\infty=Q_0+C_\infty$ with
$Q_0$ compact.  Strict negativity on the unit section of $C_\infty$ gives
\[
 \langle v,y\rangle\leq C_K-c_K|y|,
 \qquad y\in P_\infty,\quad v\in K,
\]
where $c_K>0$ for every compact $K\Subset\mathcal V(P_\infty)$.  Since
$\pi$ is a polynomial, this proves finiteness and
\eqref{g:uniform-weighted-moments}, and justifies differentiation under
the integral.

If $v\notin\mathcal V(P_\infty)$, choose $0\neq d\in C_\infty$ with
$\langle v,d\rangle\geq0$.  After applying an element of $W$, we may
assume $d\in\overline{\mathfrak a_+^*}$; the pairing is unchanged because
$v$ is central.  Choose a small ball $B$ strictly inside $P_\infty^+$.
Then $B+td\subset P_\infty^+$ for every $t\geq0$.  On this tube $\pi$
has a positive lower bound, because every positive root is positive on
$B$ and nonnegative on $d$.  The exponential factor also has a positive
lower bound.  Infinitely many disjoint translates of $B$ therefore show
that $H_\infty(v)=+\infty$.

For $0\neq a\in\mathfrak z$, differentiation gives
\[
 D^2H_\infty(v)[a,a]
 =\int_{P_\infty^+}\langle a,y\rangle^2
       e^{\langle v,y\rangle}\pi(y)\,dy>0.
\]
To prove existence of a minimizer, use $0\in\operatorname{int}P_1$.
A product of a small central ball about the origin and a small ball
strictly inside the semisimple chamber lies in $P_1^+$.  Integrating over
this product shows that
\[
 H_1(v):=\int_{P_1^+}e^{\langle v,y\rangle}\pi(y)\,dy
 \geq c\exp(c'|v|),\qquad v\in\mathfrak z,
\]
for some $c,c'>0$; when the root system is empty, only the central ball
is needed.  Hence every finite sublevel set of $H_\infty\geq H_1$ is
bounded.  Fatou's lemma makes it closed in $\mathfrak z$, and the
finiteness criterion just proved places it inside
$\mathcal V(P_\infty)$.  A minimum is attained, and strict convexity
makes it unique.  Its first variation in central directions gives
\eqref{g:limiting-central-moment}.
\end{proof}

For each $k$, let $v_k$ minimize
\[
 H_k(v)=\int_{P_k^+}e^{\langle v,y\rangle}\pi(y)\,dy,
 \qquad v\in\mathfrak z,
\]
as in the finite-parameter construction.  Define
\begin{equation}\label{g:limiting-barycenter}
 b_\infty
 =\frac{\displaystyle\int_{P_\infty^+}y
             e^{\langle v_\infty,y\rangle}\pi(y)\,dy}
        {H_\infty(v_\infty)},
 \qquad
 b_k
 =\frac{\displaystyle\int_{P_k^+}y
             e^{\langle v_k,y\rangle}\pi(y)\,dy}
        {H_k(v_k)}.
\end{equation}
Their central components vanish by the respective critical-point
equations.

\begin{prop}\label{g:central-data-convergence}
\begin{equation}\label{g:central-data-limits}
 v_k\longrightarrow v_\infty,
 \qquad H_k(v_k)\longrightarrow H_\infty(v_\infty),
 \qquad b_k\longrightarrow b_\infty.
\end{equation}
More generally, for every polynomial $p$ on $\mathfrak a^*$,
\begin{equation}\label{g:polynomial-moment-convergence}
 \int_{P_k^+}p(y)e^{\langle v_k,y\rangle}\pi(y)\,dy
 \longrightarrow
 \int_{P_\infty^+}p(y)e^{\langle v_\infty,y\rangle}\pi(y)\,dy.
\end{equation}
If $b_\infty\in2\rho+\Xi$, then $b_k\in2\rho+\Xi$ for all sufficiently
large $k$.  In particular, Theorem~\ref{existence-G} supplies a
finite-parameter singular KR soliton for every sufficiently large $k$.
\end{prop}

\begin{proof}
Fix $v_*\in\mathcal V(P_\infty)$.  Minimality and inclusion give
\[
 H_1(v_k)\leq H_k(v_k)\leq H_k(v_*)\leq H_\infty(v_*).
\]
The coercivity of $H_1$ therefore bounds $v_k$ in $\mathfrak z$.
Take any subsequence $v_{k_j}\to\bar v$.  Since
$P_{k_j}^+\uparrow P_\infty^+$, Fatou's lemma yields
\[
 H_\infty(\bar v)
 \leq\liminf_{j\to\infty}H_{k_j}(v_{k_j})
 \leq H_\infty(v_*)<\infty.
\]
Thus $\bar v\in\mathcal V(P_\infty)$.  For any
$v\in\mathcal V(P_\infty)$, the inequality
$H_{k_j}(v_{k_j})\leq H_{k_j}(v)$ and monotone convergence give
\[
 H_\infty(\bar v)\leq H_\infty(v).
\]
Uniqueness of the minimizer implies $\bar v=v_\infty$, proving convergence
of the full sequence of central vectors.  These vectors consequently
lie, for all large $k$, in a fixed compact subset of
$\mathcal V(P_\infty)$.  The uniform exponential bound in the preceding
proof now permits dominated convergence with any polynomial factor.
This proves \eqref{g:polynomial-moment-convergence} and
\eqref{g:central-data-limits}.

Finally, $b_k-2\rho$ and $b_\infty-2\rho$ belong to the semisimple root
space, since all central components vanish.  The set $\Xi$ is open in
that space.  Its strict membership for the limit therefore persists
for all sufficiently large $k$.
\end{proof}

\begin{lem}\label{g:parabolic-central-ray}
Suppose that $P_\infty$ is parabolic, with
$C_\infty=\mathbb R_{\geq0}e$ and $|e|=1$.  Then $e\in\mathfrak z^*$,
the projection of $P_\infty$ to $e^\perp$ is bounded. More precisely, there are a
compact polytope $Q\subset e^\perp$ and a continuous convex piecewise
affine function $h$ on $Q$ such that
\begin{equation}\label{g:parabolic-epigraph}
 P_\infty=\{q+te:q\in Q,\ t\geq h(q)\}.
\end{equation}
Thus $P_\infty$ has a product tail $Q\times[T_0,\infty)$, and
$\pi(q+te)=\pi(q)$.
\end{lem}

\begin{proof}
Every element of $W$ preserves the recession ray.  It sends $e$ to a
positive multiple of $e$, which must equal $e$ because $W$ is finite.
Hence $e$ is central.  The Minkowski--Weyl decomposition projects
$P_\infty$ to a compact polytope $Q$ in $e^\perp$.  In the coordinates
$y=q+te$, each retained inequality has $u_A(e)\leq0$.  Those with
$u_A(e)=0$ constrain only $q$; those with $u_A(e)<0$ give affine lower
bounds for $t$.  Their maximum is $h$.  There is at least one such lower
bound, since otherwise $-e$ would also belong to $C_\infty$.
This proves \eqref{g:parabolic-epigraph}.  Take
$T_0=\max_Qh$.  Every root vanishes on $e$, so $\pi$ is independent of
$t$.  Any central vector pairing negatively with $e$ belongs to
$\mathcal V(P_\infty)$.
\end{proof}

The following result generalizes Theorem \ref{main-resullt} to $\mathbf G$-manifolds, which confirms 
Conjecture \ref{conj-g} in the  parabolic case.

\begin{theo}\label{g-para}
Assume the fixed compactification and divisor hypotheses of
Theorem~\ref{existence-G}, and suppose that the noncompact polyhedron
$P_\infty$ is parabolic.  Let
$v_\infty$ and $b_\infty$ be defined by
Lemma~\ref{g:limiting-central-vector} and
\eqref{g:limiting-barycenter}.  If
\begin{equation}\label{g:limiting-strict-barycenter}
 b_\infty\in2\rho+\Xi,
\end{equation}
then $M=\bar M\setminus D$ admits a $K\times K$-invariant complete
gradient shrinking KR soliton $\omega_\infty$ with moment
polyhedron $P_\infty$ and central soliton vector field $X_\infty$
induced by $v_\infty$.  
Moreover, one can choose the finite-parameter solitons $\omega_k$ for all
sufficiently large $k$,
a sequence $k_j\to\infty$, and automorphisms $\sigma_j$ induced by
$Z(\mathbf G)^0$ such that
\begin{equation}\label{g:conjectural-smooth-limit}
 \sigma_j^*\omega_{k_j}\longrightarrow\omega_\infty
 \quad\text{in }C^\infty_{\mathrm{loc}}(M).
\end{equation}
\end{theo}

\begin{proof}
Proposition~\ref{g:central-data-convergence} and
Theorem~\ref{existence-G} give the finite-parameter solitons for all
sufficiently large $k$. For the convergence assertion we take integer $k$.
Write $r=\dim\mathfrak a$, $N=\dim_{\mathbb C}\mathbf G$ and
$\beta=2\sum_{\alpha>0}\alpha=2\rho$ in the convention of
\eqref{rho} and \eqref{j-convex}. The minimum of the strictly convex,
$W$-invariant potential lies in $\mathfrak a^W$. By a central
translation and addition of a constant we may therefore assume
$\varphi_k(0)=\min\varphi_k=0$. Write $u_k=\varphi_k^*$, so that
$u_k(0)=\min u_k=0$. Denote the constant introduced on the right of
the logarithmic equation by $-m_k$.

\smallskip
\noindent\emph{The central Hamiltonian.}
Put $\xi_k=-v_k$, $f_k=\partial_{\xi_k}\varphi_k$ and
$B_k=\partial_{\xi_k}^2\varphi_k$. If $e$ is the unit generator
in Lemma~\ref{g:parabolic-central-ray}, then
$\langle\xi_k,e\rangle\ge a>0$ for large $k$. The bounded transverse
projection of $P_\infty$ consequently gives a fixed $A>0$ such that
$f_k+A\ge1$. We claim that
\[
                 B_k\le A(f_k+A).
\]
We verify the maximum principle of \cite[Proposition 4.4]{CHTX} in the
presence of roots. Lift $\varphi_k$ to the invariant function $\Phi_k$
on the Euclidean space $V=i\mathfrak k$ by the adjoint action of $K$.
It is smooth by pulling back the smooth potential on $\mathbf G$ by the
Cartan exponential. At a regular Cartan point the flat Hessian on a real
root plane has the repeated eigenvalue
$\alpha(D\varphi_k)/\alpha(x)$, as follows by differentiating the
equivariance of $D\Phi_k$ along the adjoint orbit. The resulting determinant identity extends by continuity to $V$.
Its right-hand side is positive, so the Hessian is positive also at
the walls and we may write
\[
 \log\det D^2\Phi_k=-\Phi_k+\partial_{\xi_k}\Phi_k+F-m_k,
 \qquad
 F|_{\mathfrak a}=2\sum_{\alpha>0}
            \log\frac{\sinh\alpha(x)}{\alpha(x)} .
\]
Here the quotient has value $1$ at zero and the Euclidean volume
normalization can be absorbed into $m_k$. The function $F$ is constant along central directions and
has bounded first and second derivatives. Indeed
$\log(\sinh s/s)$ has bounded first and second derivatives, and on
a root plane the bound for $\alpha(DF)/\alpha(x)$ follows by reflection
across the wall and the Hessian bound on $\mathfrak a$.

Set $L_k=\Phi_k^{ij}\partial_i\partial_j-\partial_{\xi_k}$.
Differentiating in the central direction gives
\[
 L_k f_k=-f_k,\qquad L_k\log B_k\ge-1,\qquad L_k\Phi_k=N-f_k .
\]
The middle inequality follows by Cauchy--Schwarz on the quadratic
third-derivative term, as in \cite[(4.4)]{CHTX}. For fixed $k$,
$B_k$ is bounded on $V$: apply Lemma~\ref{c2-more}, as in
Proposition~\ref{sec7:dual-lower-prop}, to each height-one section
above a tangent plane. Each section is bounded since its slope is in
the interior of the gradient image, the gradients are bounded by
compactness of $P_k$, and the additional term $F$ has the derivative
bounds just proved. The resulting constant may depend on $k$.
Also $\Phi_k$ is proper since $0\in P_k^\circ$. At a maximum of
$\log(B_k/(f_k+A))-\varepsilon\Phi_k$ we consequently have
\[
 0\ge-\frac{A}{f_k+A}+\frac{B_k}{(f_k+A)^2}
               +\varepsilon(f_k-N).
\]
Since $(N-f_k)_+(f_k+A)$ is uniformly bounded and $\Phi_k\ge0$,
letting $\varepsilon\downarrow0$ proves the claim. On $M$ the function
$f_k$ is the central Hamiltonian and
$|\partial f_k|_{\omega_k}^2=B_k$ on $\mathbf G$, hence the same bound holds
on $M$ by continuity.

The claim also controls $\varphi_k$ on bounded sets of $\mathfrak a$.
Along $s\xi_k$, one has $f_k(0)=0$ and
$0\le d f_k/ds=B_k\le A(f_k+A)$. Thus $f_k$ and $\varphi_k$ are
uniformly bounded on bounded segments of that line. Every bounded $x$
can be written $x=s\xi_k+x'$, where $\langle e,x'\rangle=0$ and
$s,x'$ are uniformly bounded. The transverse derivatives of $\varphi_k$
are bounded by Lemma~\ref{g:parabolic-central-ray}, proving the assertion.

\smallskip
\noindent\emph{The weighted estimate on closed moment truncations.}
Let $\nu_k=H_k(v_k)^{-1}e^{\langle v_k,y\rangle}\pi(y)\,dy$ on
$P_k^+$, and put $I_k=\int u_k\,d\nu_k$. It is finite for fixed $k$:
the bounded relative potential of the finite soliton makes $u_k$
bounded on $P_k$. The mean of $\nu_k$ is $b_k$.
Choose $s<1$, close enough to $1$ that $s b_k-\beta$ belongs to the
closed positive root cone for every sufficiently large $k$.
This uses the strict limiting barycenter assumption and vanishing of
the central components. Choose $\eta>0$ with $c_0=s(1+\eta)<1$.
For some $d>0$ independent of $k$, every $b_k+z$, $|z|\le d$, is the
mean of a probability measure dominated by $(1+\eta)\nu_k$. Indeed,
choose a smooth bump $\zeta$ of integral one in a fixed ball in
$(P_1^+)^\circ$ and use $d\nu_k-z\cdot D\zeta\,dy$. On the ball the
densities of $\nu_k$ have a uniform positive lower bound.
Jensen's inequality and $u_k(0)=0$ give
$u_k(s(b_k+z))\le c_0I_k$. Maximizing the corresponding Legendre
lower bound over $z$ gives, for $x\in\overline{\mathfrak a_+}$,
\[
 \varphi_k(x)\ge\langle\beta,x\rangle+sd|x|-c_0 I_k .
\]

Let $d\mu_k=J e^{-\varphi_k}dx/\int J e^{-\varphi_k}dx$ on the
positive chamber. By \eqref{KR-soliton-G},
$(D\varphi_k)_*\mu_k=\nu_k$. Legendre duality and integration of
$\operatorname{div}(xJ e^{-\varphi_k})$ give
\[
 I_k=r+\int(x\cdot D\log J-\varphi_k)\,d\mu_k
     \le N+\int(\langle\beta,x\rangle-\varphi_k)\,d\mu_k
     \le N+c_0 I_k .
\]
Here $t\coth t\le t+1$, $N=r+2|\Phi_+|$, and $J$ vanishes on the
chamber walls. The preceding lower bound for $\varphi_k$ gives
exponential decay for each fixed $k$ and justifies integration at
infinity. It now gives the uniform estimate
$J e^{-\varphi_k}\le C e^{-sd|x|}$. Since $\varphi_k$ is uniformly
bounded on one fixed interior ball, $\int J e^{-\varphi_k}$ is
bounded above and below. The mass identity and
Proposition~\ref{g:central-data-convergence} therefore give $|m_k|\le C$.

For a fixed closed truncation $K_R=P_\infty\cap\{t\le R\}$, with
$y=q+te$ as in Lemma~\ref{g:parabolic-central-ray}, the change of
variables in Lemma~\ref{bound} now gives
$\int_{K_R}|Du_k|^p\pi\,dy\le C_{p,R}$ for every finite $p$.
We used $W$-symmetry to pass to the full truncation.
A sufficiently small negative power of $\pi$ is locally integrable
even at intersections of root walls, by H\"older's inequality applied
to its finitely many linear factors. Hence a second application of
H\"older, using the preceding estimate at a higher exponent, gives
$\int_{K_R}|Du_k|^pdy\le C'_{p,R}$. Morrey's inequality and $u_k(0)=0$
give the bounds of Lemma~\ref{bound} on every closed truncation.

\smallskip
\noindent\emph{Localization and passage to the limit.}
On $Z\cap M$ the proof of Lemma~\ref{sec7:geo-proper-character}
applies to the primitive character on the recession ray. All moving
inward normals pair negatively with it, since they form one $W$-orbit,
and the only retained facets meeting a sufficiently distant moving
facet are vertical. Thus the character has only poles near the removed
divisors, including their intersections, and is proper on $Z\cap M$.
The two secant inequalities of Lemma~\ref{sec7:geo-localization}
apply using the closed-truncation bound and the bounds for $\varphi_k$
on fixed balls. By the Cartan decomposition and compactness,
$\bar M=(K\times K)Z$, so saturation gives both localization inclusions
on $M$. They also hold with $f_k$ as height, since
$f_k=\langle\xi_k,e\rangle t+\langle\xi_k,q\rangle$,
$\langle\xi_k,e\rangle\ge a$ and $q$ is bounded.

Choose the fixed invariant reference metric as in
Lemma~\ref{sec7:geo-background-bounds}. Its Legendre comparison and
the boundedness of $m_k$ give the same relative potential, Hamiltonian
and volume-ratio bounds. The relative potential is smooth on $M$:
the defining section of $D$ identifies the polarizations there with
$-K_M$. The reference volume induced by the resulting Hermitian metric
on $-K_M$ gives the smooth reference term in the relative equation.
Thus its derivatives are uniformly bounded on fixed compact sets.
For the complex drift Laplacian on $M$ one has
$(\Delta_{\omega_k}+X_k) f_k=-f_k$. On the regular chamber this is the
central derivative of \eqref{KR-soliton-G}: the term
$D\log\pi(D\varphi_k)\cdot D f_k$ is the root contribution to
$\Delta_{\omega_k} f_k$. Smoothness extends the identity throughout $M$.
Our bound on $|\partial f_k|_{\omega_k}^2$ now gives exactly the
cutoff estimate of Section~\ref{sec:boundary-estimates} on each
$\{f_k\le R\}$. Lemma~\ref{sec7:yau-drift-trace},
the proof of Proposition~\ref{sec7:yau-local-equivalence}, and the
bootstrap following it give smooth local convergence along a subsequence
to a positive metric $\omega_\infty$. The pullbacks used here are
precisely the central translations chosen at the start of the proof.
Passing to the equation yields
$\operatorname{Ric}(\omega_\infty)-\omega_\infty
=\mathcal L_{X_\infty}\omega_\infty$ and
$\nabla^{1,0} f_\infty=-X_\infty$.

The normalized moment maps on $Z\cap M$ converge locally smoothly,
by the reference Hamiltonian formula. For each $y\in P_\infty$,
including its finite faces, choose preimages for the finite moment maps
for all large $k$. These lie in a fixed compact subset of $Z\cap M$
by localization, so their limit along a subsequence realizes $y$.
The reverse inclusion follows from $P_k\subset P_\infty$ and
closedness. Finally $f_\infty$ is proper: for each $R$, pointwise
convergence places $\{f_\infty\le R\}$ in the common compact set
containing $\{f_k\le R+1\}$ for all large $k$. The central estimate
passes to the limit and gives
\[
 |\partial\sqrt{f_\infty+A}|_{\omega_\infty}^2\le A/4 .
\]
This proper function with bounded gradient proves completeness, as in
\cite[Proposition 4.1]{CHTX}.
\end{proof}

\section{Appendix: Examples}
\label{sec9:examples}

In this appendix, we discuss several examples to show a relationship between  singular KR 
solitons  and complete gradient  shrinking  KR solitons. 

\subsection{One dimensional toric case}
\label{sec9:dimension-one}

In complex dimension one, after reversing the lattice coordinate
if necessary, the noncompact Fano moment interval is
$P^\circ=(-1,\infty)$. The compact approximations are
$P_k^\circ=(-1,k)$, with $\bar M=\mathbb P^1$ and
$D=\{\infty\}$; in particular, $\mathcal O(D)=\mathcal O(1)$
is ample. For $a<0$, the weighted volume is
\[
 H_0(a)=\int_{-1}^{\infty}e^{ay}\,dy=-\frac{e^{-a}}a,
 \qquad H_0'(a)=\frac{e^{-a}(a+1)}{a^2}.
\]
Consequently the vanishing weighted barycenter selects $A=-1$.
The normalized solution and its Legendre transform are
\begin{equation}\label{sec9:gaussian-potentials}
 \varphi(x)=e^x-x-1,\qquad m=0,\qquad
 u(y)=(y+1)\log(y+1)-y.
\end{equation}
Indeed, $y=\varphi'(x)=e^x-1$, $u'(y)=\log(y+1)$, and
\[
 \varphi''(x)=e^x=e^{-\varphi(x)+\varphi'(x)},
 \qquad u''(y)=\frac1{y+1}.
\]
Thus the retained endpoint has the required Guillemin singularity.
On $\mathbb C^*$, the term $-\log|Z|^2-1$ is pluriharmonic,
and the K\"ahler form extends to
\[
 \omega=\sqrt{-1}\,dZ\wedge d\overline Z
 \quad\hbox{on }\mathbb C.
\]
This is the complete Gaussian shrinking soliton, with
$X=-Z\partial_Z$ and $f=|Z|^2-1$.

The compact approximation can also be written explicitly in
inverse-Hessian coordinates. Let $a_k$ be the unique number with
\begin{equation}\label{sec9:interval-barycenter}
 \int_{-1}^{k}y e^{a_ky}\,dy=0.
\end{equation}
Strict monotonicity of this integral as a function of $a_k$
gives uniqueness. Moreover, $a_1=0$, $-1<a_k<0$ for $k>1$,
and $a_k\downarrow-1$ as $k\to\infty$.
For $q_k(y)=1/u_k''(y)$, differentiating the one-dimensional
equation \eqref{RMA-2} gives
\begin{equation}\label{sec9:interval-ode}
 q_k'(y)+a_kq_k(y)=-y,\qquad q_k(-1)=q_k(k)=0.
\end{equation}
The endpoint conditions follow from $u_k'((-1,k))=\mathbb R$:
the linear equation extends $q_k$ smoothly to each finite endpoint,
and a positive endpoint value would make the integral of $q_k^{-1}$
finite there, contradicting this gradient image.
Its solution is
\begin{align}\label{sec9:interval-inverse-hessian}
 q_k(y)
 &=-e^{-a_ky}\int_{-1}^{y}t e^{a_kt}\,dt\notag\\
 &=\frac{1-a_ky-(1+a_k)e^{-a_k(y+1)}}{a_k^2}
       \quad (k>1),
 \qquad q_1(y)=\frac{1-y^2}{2}.
\end{align}
The first expression and \eqref{sec9:interval-barycenter}
show that $q_k>0$ on $(-1,k)$. Integrating $u_k''=q_k^{-1}$
with $u_k(0)=u_k'(0)=0$ therefore defines a strictly convex
potential whose derivative maps $(-1,k)$ onto $\mathbb R$.
Its Legendre transform is normalized as above and satisfies
\eqref{RMA-2} with
\[
 m_k=-\log q_k(0).
\]

The endpoint slopes follow directly from
\eqref{sec9:interval-ode}:
\[
 q_k'(-1)=1,\qquad q_k'(k)=-k.
\]
In particular,
\begin{equation}\label{sec9:interval-boundary}
 u_k(y)=(y+1)\log(y+1)
       +\frac1k(k-y)\log(k-y)+r_k(y),
 \qquad r_k\in C^\infty([-1,k]).
\end{equation}
The coefficient at the moving endpoint is $1/k$. More explicitly,
with $w=Z^{-1}$ a coordinate at $D$, one has
$k-y\sim c_k|w|^{2k}$ for a positive constant $c_k$.
The associated Riemannian metric consequently has leading form
$c_k'|w|^{2k-2}|dw|^2$, with $c_k'>0$, and cone angle $2\pi k$.
Thus these are the large-angle compact singular solitons used
in the approximation. Finally, \eqref{sec9:interval-inverse-hessian}
gives $q_k\to y+1$ smoothly on every fixed compact subinterval
of $[-1,\infty)$ and $m_k\to0$. The normalized potentials and
metrics converge to \eqref{sec9:gaussian-potentials} and the
Euclidean metric, respectively.

\subsection{Two dimensional parabolic toric case}
\label{sec9:polygon-classification}

In dimension two the relevant polyhedra admit an elementary
classification. This is a classification of the parabolic
Fano polyhedra with primitive normals and smooth finite
vertices in the present construction. It concerns the noncompact
polyhedron $P$ and its toric manifold, rather than the pair
$(P,a_0)$: different choices of the moving facet can give
different compactifications. In the two examples below we
choose $a_0=e_2$.

\begin{prop}\label{sec9:two-polygons}
Let $\overline P\subset\mathbb R^2$ be a two-dimensional
parabolic Fano polyhedron in the sense used above, with
primitive integral facet normals. Up to a change of lattice
coordinates in $\operatorname{GL}(2,\mathbb Z)$, it is one of
\begin{align}\label{sec9:classified-polygons}
 \overline P_{\mathrm{prod}}
   &=\{(r,s):-1\le r\le1,\ s\ge-1\},\notag\\
 \overline P_{\mathrm{bl}}
   &=\{(r,s):-1\le r\le1,\ s\ge-1,\ r+s\ge-1\}.
\end{align}
\end{prop}

\begin{proof}
Choose lattice coordinates in which the recession cone is
$\mathbb R_{\ge0}e_2$. Every primitive inward normal then has
the form $(a,b)$ with $b\ge0$. Boundedness of the horizontal
projection requires vertical supporting lines on both sides.
Their primitive inward normals are $(1,0)$ and $(-1,0)$;
the Fano normalization therefore makes these lines $r=-1$
and $r=1$.

The remaining facets form a finite chain of compact lower
edges from the left supporting line to the right one.
Write their primitive inward normals as $(a_i,b_i)$,
where $b_i\in\mathbb Z_{>0}$. Every finite vertex is a
lattice point: its two adjacent primitive normals form a
lattice basis by the Delzant condition, and both defining
constants are $-1$. Hence the lattice length $\ell_i$ of
each lower edge is a positive integer. When an edge is
oriented from left to right, its displacement vector is
$\ell_i(b_i,-a_i)$. Adding horizontal displacements gives
\begin{equation}\label{sec9:edge-width}
 \sum_i b_i\ell_i=2.
\end{equation}
There are therefore at most two lower edges.

If there is only one lower edge, its normal is $(a,b)$.
The Delzant condition at its intersection with a vertical
facet gives $b=1$, and \eqref{sec9:edge-width} gives
$\ell=2$. The integral shear $(r,s)\mapsto(r,s+ar)$
then gives $\overline P_{\mathrm{prod}}$.
If there are two lower edges, \eqref{sec9:edge-width}
forces $b_1=b_2=\ell_1=\ell_2=1$. Denote their first
normal components by $a_L,a_R$, in left-to-right order.
Convexity and the Delzant condition at their common
vertex give $a_L-a_R=1$. After the shear
$(r,s)\mapsto(r,s+a_Rr)$, their inequalities become
$r+s\ge-1$ and $s\ge-1$. This is
$\overline P_{\mathrm{bl}}$.
\end{proof}

\begin{figure}[htbp]
\centering
\begingroup
\setlength{\unitlength}{1pt}
\begin{picture}(170,175)
 \thinlines
 \put(15,80){\vector(1,0){140}}
 \put(85,15){\vector(0,1){140}}
 \put(157,77){\makebox(0,0)[l]{$r$}}
 \put(89,157){\makebox(0,0)[l]{$s$}}
 \thicklines
 \put(45,40){\line(1,0){80}}
 \put(45,40){\vector(0,1){100}}
 \put(125,40){\vector(0,1){100}}
 \put(45,40){\circle*{3}}
 \put(125,40){\circle*{3}}
 \put(43,31){\makebox(0,0)[r]{$(-1,-1)$}}
 \put(127,31){\makebox(0,0)[l]{$(1,-1)$}}
 \put(85,5){\makebox(0,0){$\overline P_{\mathrm{prod}}$}}
\end{picture}
\hspace{28pt}
\begin{picture}(170,175)
 \thinlines
 \put(15,80){\vector(1,0){140}}
 \put(85,15){\vector(0,1){140}}
 \put(157,77){\makebox(0,0)[l]{$r$}}
 \put(89,157){\makebox(0,0)[l]{$s$}}
 \multiput(45,40)(4,0){10}{\line(1,0){2}}
 \multiput(45,40)(0,4){10}{\line(0,1){2}}
 \thicklines
 \put(45,80){\line(1,-1){40}}
 \put(85,40){\line(1,0){40}}
 \put(45,80){\vector(0,1){60}}
 \put(125,40){\vector(0,1){100}}
 \put(45,80){\circle*{3}}
 \put(85,40){\circle*{3}}
 \put(125,40){\circle*{3}}
 \put(41,91){\makebox(0,0)[r]{$(-1,0)$}}
 \put(82,29){\makebox(0,0)[r]{$(0,-1)$}}
 \put(128,29){\makebox(0,0)[l]{$(1,-1)$}}
 \put(85,5){\makebox(0,0){$\overline P_{\mathrm{bl}}$}}
\end{picture}
\endgroup
\caption{The two parabolic Fano polygons in complex dimension two.
The regions lie above the lower boundary and between $r=-1$ and
$r=1$; arrows indicate the unbounded direction. On the right,
the added facet $r+s=-1$ cuts off the lower-left corner of the
product polygon, shown by dashed segments.}
\label{sec9:polygon-figure}
\end{figure}
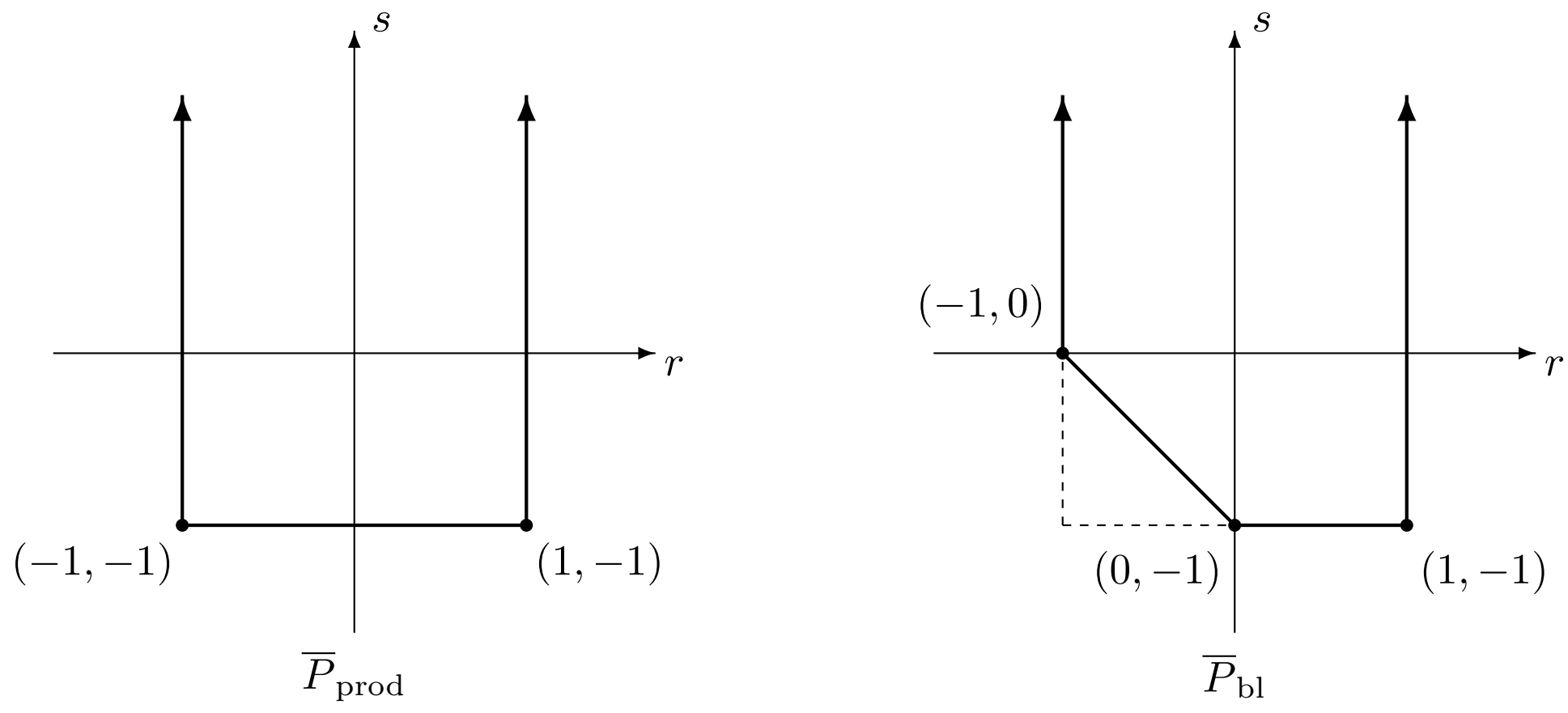

The two canonical cuts at the two lower corners of the
product strip would instead impose both $r+s\ge-1$
and $-r+s\ge-1$. The original horizontal lower edge
then disappears, and the two remaining lower normals
$(1,1)$ and $(-1,1)$ have determinant of absolute value
$2$ at $(0,-1)$. Thus this two-point blowup does not give
another Delzant polyhedron with all facet constants equal
to $1$.

\subsubsection{The product on $\mathbb P^1\times\mathbb C$}
\label{sec9:product}

For $\overline P_{\mathrm{prod}}$, take
\[
 P_k^\circ=(-1,1)\times(-1,k),\qquad k\ge1.
\]
The compactification is $\bar M=\mathbb P^1\times\mathbb P^1$,
and $D=\mathbb P^1\times\{\infty\}$. Its line bundle is
$\mathcal O(D)=\operatorname{pr}_2^*\mathcal O(1)$,
which is generated by global sections and hence semiample.
The weighted integral factors into the two one-dimensional
integrals, so its critical vector is $A=(0,-1)$.

The normalized potential is
\begin{equation}\label{sec9:product-potential}
 \varphi(x_1,x_2)
   =2\log\cosh(x_1/2)+e^{x_2}-x_2-1,
 \qquad m=\log2.
\end{equation}
Indeed,
\[
 D\varphi=\bigl(\tanh(x_1/2),e^{x_2}-1\bigr),
 \qquad
 \det D^2\varphi
   =\tfrac12\operatorname{sech}^2(x_1/2)e^{x_2}
   =e^{-\varphi-\log2+\varphi_{x_2}}.
\]
The Legendre potential is
\begin{equation}\label{sec9:product-dual}
 u(r,s)=(1+r)\log(1+r)+(1-r)\log(1-r)
           +(1+s)\log(1+s)-s.
\end{equation}
In particular, the Guillemin remainder is the smooth
affine function $-s$. On the affine chart of the first
factor, the metric is
\[
 \omega
  =2\sqrt{-1}\,\partial\bar\partial\log(1+|Z_1|^2)
      +\sqrt{-1}\,dZ_2\wedge d\overline Z_2.
\]
It is the complete product of the K\"ahler--Einstein metric
on $\mathbb P^1$ with $\operatorname{Ric}(\omega_{\mathbb P^1})
=\omega_{\mathbb P^1}$ and the Gaussian metric on $\mathbb C$.
Here $X=-Z_2\partial_{Z_2}$ and $f=|Z_2|^2-1$.
The compact approximations also split: if
$\varphi_k^{(1)}$ and $m_k^{(1)}$ are the one-dimensional
objects above, then
\[
 \varphi_k=2\log\cosh(x_1/2)+\varphi_k^{(1)}(x_2),
 \quad A_k=(0,a_k),\quad m_k=\log2+m_k^{(1)}.
\]
The uniqueness theorem for complete toric shrinking solitons
\cite[Theorem A]{Ci} identifies the limit furnished by
Theorem~\ref{smoothness} and the completeness argument following
\eqref{sec7:extended-soliton} with this product, up to biholomorphism.

\subsubsection{The blowup of $\mathbb P^1\times\mathbb C$ at one point}
\label{sec9:bccd}

For $\overline P_{\mathrm{bl}}$, choose
\begin{equation}\label{sec9:bccd-truncations}
 P_k^\circ
   =\{(r,s):-1<r<1,\ -1<s<k,\ r+s>-1\},
 \qquad k\ge1.
\end{equation}
The vertices of its closure, in boundary order, are
\[
 (-1,0),\quad(-1,k),\quad(1,k),\quad(1,-1),\quad(0,-1).
\]
The adjacent primitive normals have determinant of
absolute value $1$, so every $P_k$ is Delzant. The new
normal $(1,1)$ is the sum of the normals $(1,0)$ and
$(0,1)$ at the lower left corner of the product polygon.
Thus
\[
 \bar M=\operatorname{Bl}_{p}(\mathbb P^1\times\mathbb P^1),
 \qquad
 M=\operatorname{Bl}_{p}(\mathbb P^1\times\mathbb C),
\]
where $p$ is the torus fixed point corresponding to
$(-1,-1)$ before the blowup. The divisor $D$ is the
strict transform of $\mathbb P^1\times\{\infty\}$.
Since the blowup center lies on the zero fiber, it is
disjoint from this divisor. For the blowdown map $\pi$,
\[
 \mathcal O(D)=\pi^*\operatorname{pr}_2^*\mathcal O(1),
\]
so $D$ is again semiample. This is the toric model of
\cite[Example 2.33]{CCD-type1}, with the two moment
coordinates interchanged.

The weighted barycenter vector has a particularly simple
description in these coordinates. The integral linear
involution
\[
 (r,s)\longmapsto(-r,s+r)
\]
preserves $\overline P_{\mathrm{bl}}$, Lebesgue measure,
and the function $r+2s$. Consequently
\[
 \int_{P_{\mathrm{bl}}}r e^{-\lambda(r+2s)}\,dr\,ds=0
 \qquad(\lambda>0).
\]
Integrating first in $s$ gives
\begin{align}\label{sec9:bccd-volume}
 \mathcal H(\lambda)
  &:=\int_{P_{\mathrm{bl}}}e^{-\lambda(r+2s)}\,dr\,ds\notag\\
  &=\frac{e^{2\lambda}}{2\lambda}
       \left(\int_{-1}^0e^{\lambda r}\,dr
              +\int_0^1e^{-\lambda r}\,dr\right)
    =\frac{e^{2\lambda}-e^\lambda}{\lambda^2}.
\end{align}
At a critical point of $\mathcal H$, the vanishing of
the $r$-moment also forces the $s$-moment to vanish.
The function $\mathcal H$ is strictly convex and tends
to infinity at both ends of $(0,\infty)$. Since
$\mathcal H'(1)=e>0$, its unique critical point is the unique solution
$\lambda\in(0,1)$ of
\begin{equation}\label{sec9:bccd-root}
 2e^\lambda(1-\lambda)=2-\lambda,
 \qquad \lambda\simeq0.6437977582.
\end{equation}
The formal root $\lambda=0$ of the displayed equation
is excluded because the weighted volume
is infinite there. Strict convexity of the full weighted
volume functional now gives
\begin{equation}\label{sec9:bccd-vector}
 A_0=(-\lambda,-2\lambda),\qquad
 X_0=-\lambda\bigl(Z_1\partial_{Z_1}
                     +2Z_2\partial_{Z_2}\bigr),
 \qquad f=\lambda(\mu^1+2\mu^2)+\text{constant}.
\end{equation}
This calculation agrees with the soliton vector determined
in \cite[Example 2.33]{CCD-type1}, after accounting for
the coordinate order and the sign convention for $X$.

Bamler--Cifarelli--Conlon--Deruelle proved the existence
of a complete torus-invariant shrinking gradient
KR soliton on this blowup
\cite[Theorem A]{BCCD}. Our limiting metric is smooth
and torus-invariant by Theorem~\ref{smoothness} and
complete by the argument following \eqref{sec7:extended-soliton}.
Since this
toric surface has finitely many torus fixed points,
the uniqueness theorem \cite[Theorem A]{Ci} identifies
it, up to biholomorphism, with that soliton. This use
of toric uniqueness requires no additional curvature
bound on the metric constructed here.

\begin{rem}\label{sec9:scope}
Proposition~\ref{sec9:two-polygons} concerns the
one-dimensional recession cone imposed by the
parabolic hypothesis. For comparison, the canonical
polyhedra of $\mathbb C^2$ and
$\operatorname{Bl}_0\mathbb C^2$ are respectively
\[
 \{(r,s):r\ge-1,\ s\ge-1\},\qquad
 \{(r,s):r\ge-1,\ s\ge-1,\ r+s\ge-1\}.
\]
Both have recession cone $\mathbb R_{\ge0}^2$.
The Gaussian soliton on the former and the
Feldman--Ilmanen--Knopf soliton on the latter
\cite{FIK} therefore lie outside the parabolic
classification of this section.
\end{rem}

\subsection{A non-product parabolic compactification of
\texorpdfstring{$SL(2,\mathbb C)\times\mathbb C^*$}{SL(2,C) times C*}}
\label{sec9:sl2-blowup}

Let $Y=Q^3=\{[A:\tau]\in
\mathbb P(\operatorname{Mat}_{2\times2}(\mathbb C)\oplus\mathbb C):
\det A=\tau^2\}$. The left--right action on $A$ fixes $\tau$.
The chart $\tau=1$ is $SL(2,\mathbb C)$, and
$B=\{\tau=0\}\simeq\mathbb P^1\times\mathbb P^1$ is its closed orbit.
For $\mathbf G=SL(2,\mathbb C)\times\mathbb C^*$, consider the equivariant blowup
\[
 \beta:\bar M=\operatorname{Bl}_{B\times\{0\}}(Y\times\mathbb P^1)
       \longrightarrow Y\times\mathbb P^1,\qquad
 f=\operatorname{pr}_2\circ\beta,\qquad D=f^{-1}(\infty).
\]
The central factor acts by multiplication on the affine coordinate $Z$
of $\mathbb P^1$. The open orbit is $\mathbf G$. In the classification of
Delcroix and Montagard this blowup is the embedding $4$-$2$-$3$
\cite[Example~6.4]{DM}.

Write $H$ and $F$ for the pullbacks of the hyperplane class of $Y$
and the point class of $\mathbb P^1$, and let $E$ be the exceptional divisor.
Adjunction gives $-K_Y=\mathcal O_Y(3)$, and the center has
codimension two. Thus
\[
 -K_{\bar M}=3H+2F-E,\qquad
 [D]=F,\qquad
 \Omega_k/(2\pi)=3H+(k+1)F-E .
\]
To see ampleness, the center is the complete intersection of divisors
of classes $H$ and $F$. Its ideal sheaf twisted by $H+F$ is globally
generated. The resulting sections of $H+F-E$ realize the blowup as a
closed graph in $(Y\times\mathbb P^1)\times\mathbb P^N$.
Hence $(2H+F)+(H+F-E)=-K_{\bar M}$ is the restriction of an ample
external product. The smooth divisor $D\simeq Q^3$ is semiample since
$\mathcal O(D)=f^*\mathcal O_{\mathbb P^1}(1)$.

\smallskip
\noindent\emph{Non-product geometry.}
The general fiber of $f$ is $Q^3$, whereas
$f^{-1}(0)=F_0\cup E$, where $F_0$ is the strict transform of
$Y\times\{0\}$. Blowing up the Cartier divisor $B\subset Y$ leaves
$F_0\simeq Q^3$, and $F_0\cap E\simeq B$. Since $f|_M$ is proper with
connected fibers over the normal affine line, $\Gamma(M,\mathcal O_M)
=\mathbb C[Z]$. Its affinization thus has a reducible fiber, which
excludes an algebraic product $Y'\times\mathbb A^1$ for a connected
smooth projective $Y'$.

The compactification itself is also indecomposable. If
$\bar M=U\times V$ with both factors positive-dimensional, both are
Fano, and $\operatorname{Pic}^0(U)=\operatorname{Pic}^0(V)=0$.
Thus $f^*\mathcal O(1)$ is an external product of nef line bundles.
Its square is zero. Intersecting with products of curves shows that
one factor is numerically trivial, so $f$ is constant on the
corresponding slices and factors through the other factor, say $U$.
A general fiber $Q^3$ has Picard number one and hence cannot split as
a product of two positive-dimensional projective varieties. Thus $U$
is a curve. Connectedness of the fibers implies that $U\to\mathbb P^1$
has degree one. This would make every fiber isomorphic to $V$,
contrary to the reducibility of $f^{-1}(0)$.

\smallskip
\noindent\emph{The moment polygons.}
Parametrize a maximal torus of $SL(2,\mathbb C)$ by
$\operatorname{diag}(z,z^{-1})$. Its closure in $Y$ is the conic
$ab=\tau^2$ of diagonal matrices. It has hyperplane degree two, so
$-K_Y$ has degree six on it. In the weight coordinate $y$ for the
character $z\mapsto z$, Weyl symmetry centers its moment interval at
zero, giving $[-3,3]$. Write $t$ for the central moment coordinate,
oriented so that $D$ is the upper facet. Before the blowup the polygon
is $[-3,3]\times[-1,k]$.

\begin{figure}[htbp]
\centering
\begin{tikzpicture}[
  x=.76cm,y=1.00cm,font=\small,
  axis/.style={->,>=stealth,draw=black!52,line width=.4pt},
  fixed/.style={draw=black!82,line width=.85pt},
  moving/.style={draw=blue!65!black,line width=1.6pt},
  former/.style={draw=blue!65!black,dashed,line width=.65pt},
  exceptional/.style={draw=orange!70!black,line width=1.25pt},
  wall/.style={draw=black!56,densely dotted,line width=.65pt},
  ticklabel/.style={font=\scriptsize,inner sep=1pt}]
  \foreach \panelshift/\islimit in {0/0,6.55/1} {
    \begin{scope}[xshift=\panelshift cm]
      \ifnum\islimit=0
        \def\topheight{2}
      \else
        \def\topheight{2.34}
      \fi
      \fill[black!4] (-3,\topheight) -- (-3,0) -- (-2,-1)
         -- (0,-1) -- (0,\topheight) -- cycle;
      \fill[blue!10] (0,\topheight) -- (0,-1) -- (2,-1)
         -- (3,0) -- (3,\topheight) -- cycle;
      \draw[axis] (-3.56,-1.39) -- (3.49,-1.39) node[right,black] {$y$};
      \draw[axis] (-3.56,-1.39) -- (-3.56,2.51) node[above,black] {$t$};
      \foreach \yv in {-3,-2,0,2,3}{
        \draw[black!55] (\yv,-1.39) -- +(0,-.055)
          node[ticklabel,below=2pt] {$\yv$};
      }
      \foreach \tv in {-1,0,1,2}{
        \draw[black!55] (-3.56,\tv) -- +(-.07,0)
          node[ticklabel,left=2pt] {$\tv$};
      }
      \draw[wall] (0,-1) -- (0,\topheight);
      \draw[fixed] (-2,-1) -- (2,-1);
      \draw[exceptional] (-3,0) -- (-2,-1);
      \draw[exceptional] (2,-1) -- (3,0);
      \node[ticklabel,text=orange!65!black] at (-2.12,-.51) {$E$};
      \node[ticklabel,text=orange!65!black] at (2.12,-.51) {$E$};
      \ifnum\islimit=0
        \draw[fixed] (-3,0) -- (-3,2) (3,0) -- (3,2);
        \draw[former] (-3,1) -- (3,1);
        \draw[moving] (-3,2) -- (3,2)
          node[midway,above=2pt,text=blue!65!black] {$D$};
        \draw[->,>=stealth,blue!65!black,line width=.8pt]
          (-1.48,1.13) -- (-1.48,1.87);
        \node at (1.48,.40) {$P_2^+$};
        \node at (0,2.83) {$P_2$};
      \else
        \draw[fixed,->,>=stealth] (-3,0) -- (-3,2.52);
        \draw[fixed,->,>=stealth] (3,0) -- (3,2.52);
        \draw[->,>=stealth,blue!65!black,line width=.85pt]
          (1.48,.95) -- (1.48,2.13);
        \node at (1.48,.40) {$P_\infty^+$};
        \node at (0,2.83) {$P_\infty$};
      \fi
    \end{scope}
  }
\end{tikzpicture}
\caption{Moving the divisor $D=f^{-1}(\infty)$ on the blowup of
$Q^3\times\mathbb P^1$. The upper edge moves from $t=1$ (dashed)
to $t=2$ (left) and escapes to infinity (right). The two sloping
edges marked $E$ are exchanged by the Weyl group and represent the
exceptional divisor on the torus closure. The bottom edge represents
the strict transform of $Q^3\times\{0\}$, and the vertical sides
represent the strict transform of $Q^2\times\mathbb P^1$.
The dotted line is the Weyl wall $y=0$ and the blue shading marks
the positive chamber.}
\label{sec9:sl2-blowup-polytopes-figure}
\end{figure}
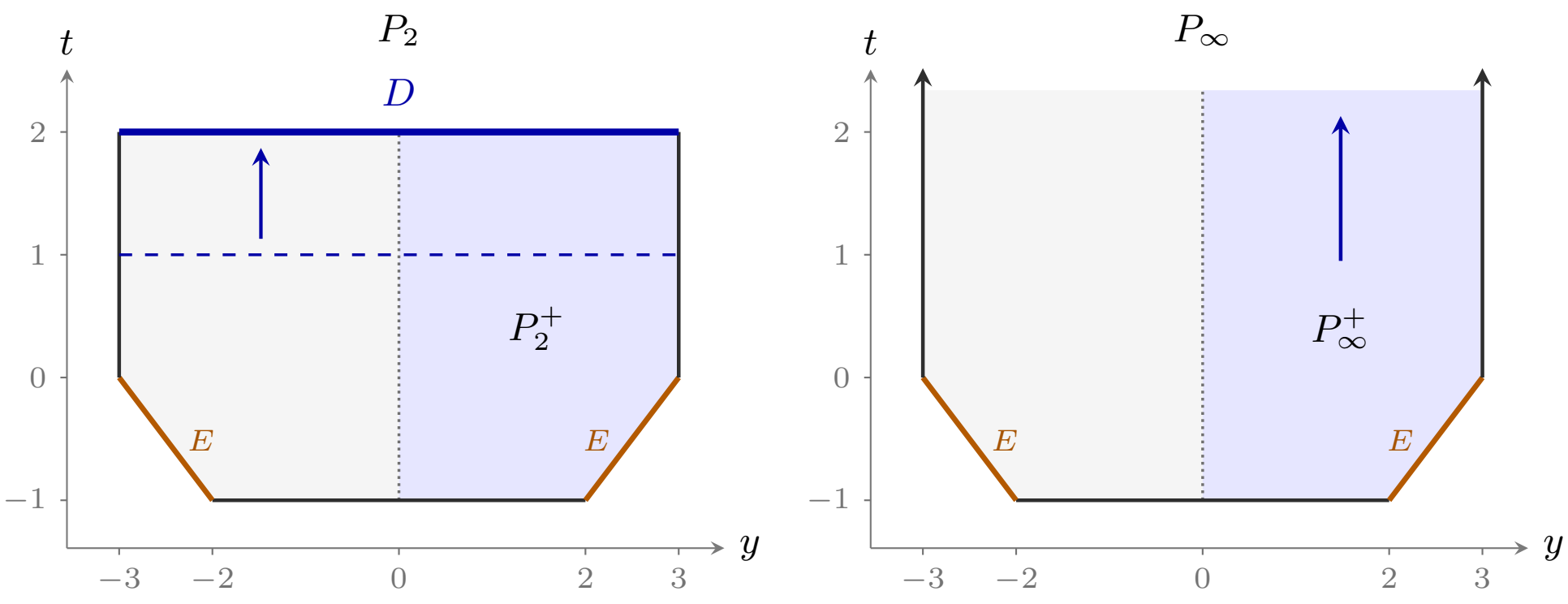

The conic meets $B$ transversely at its two endpoints $\tau=0$.
Thus the maximal torus closure in $\bar M$ is the blowup of the torus
surface at its two lower fixed points. At each point the exceptional
normal is the sum of the adjacent primitive normals, namely
$(1,-1)$ or $(-1,-1)$. The coefficient $-1$ of $E$ lowers the new
support number from $3+1$ to $3$. Consequently, with
$R(t)=3+\min(t,0)$,
\begin{equation}\label{sec9:sl2-blowup-polygons}
 \begin{aligned}
 P_k&=\{(y,t):-1\le t\le k,\ |y|\le3,\ |y|-t\le3\},\\
 P_\infty&=\{(y,t):t\ge-1,\ |y|\le R(t)\},\qquad
 \operatorname{rec}(P_\infty)=\mathbb R_{\ge0}(0,1).
 \end{aligned}
\end{equation}
The vertices of $P_k$ are $(\pm2,-1)$, $(\pm3,0)$ and $(\pm3,k)$.
The Weyl wall $y=0$ therefore meets the boundary only in the relative
interiors of facets, as required by the wonderful hypothesis.
Only the facet $t=k$ moves. Its removal leaves a parabolic polyhedron
whose sections widen from $[-2,2]$ at $t=-1$ to $[-3,3]$ at $t=0$.
For $t\ge0$ the polyhedron is exactly $[-3,3]\times[0,\infty)$.
Figure~\ref{sec9:sl2-blowup-polytopes-figure} shows the two cuts and
the moving facet.

\smallskip
\noindent\emph{The weighted barycenter.}
In these coordinates $2\rho=(2,0)$, $\Xi=\mathbb R_{>0}(1,0)$,
and the root density on the positive chamber is a positive constant
times $y^2$. For $1\le k\le\infty$, write $v_k=(0,-a_k)$.
Integration over $0\le y\le R(t)$ gives the central critical-point
equation and the semisimple moment formula
\begin{equation}\label{sec9:sl2-blowup-moments}
 \int_{-1}^k t R(t)^3 e^{-a_k t}\,dt=0,\qquad
 b_{k,y}=\frac34\,
 \frac{\int_{-1}^k R(t)^4 e^{-a_k t}\,dt}
      {\int_{-1}^k R(t)^3 e^{-a_k t}\,dt}.
\end{equation}
The central component is zero by the first identity. To compute the
limiting vector, put $a=a_\infty>0$. Splitting the first integral at
zero and integrating on the positive half-line gives the scalar equation
\[
 \int_0^1 s(3-s)^3e^{as}\,ds=\frac{27}{a^2}.
\]
Its left side increases strictly and its right side decreases strictly.
At $a=1$, the function $e^s(3-s)^3$ decreases on $[0,1]$, so the left
side is at most $27/2$. Hence the unique solution satisfies $a>1$.

Normalize $e^{-at}R(t)^3\,dt$ to a probability measure and let $T$ have
this law. Set $T_-=\max(-T,0)$ and $m=\mathbb E[-T\mid T<0]$.
On $T>0$ the conditional law is exponential of mean $1/a$.
As $\mathbb ET=0$, it follows that
$\mathbb E T_-=m/(1+am)$. On $T<0$, the variable $s=-T$ has
density $w(s)$ proportional to $e^{as}(3-s)^3$ on $[0,1]$.
Integrating $(s(1-s)w(s))'$ and using $a>1$ gives
\[
 2m-1
 =\int_0^1 s(1-s)\left(a-\frac3{3-s}\right)w(s)\,ds
 <(a-1)m.
\]
Indeed $3/(3-s)\ge1$ and $s(1-s)<s$ for $0<s<1$.
Thus $3m<1+am$, and \eqref{sec9:sl2-blowup-moments} yields
\begin{equation}\label{sec9:sl2-blowup-strict}
 \mathbb E T_-<\tfrac13,\qquad
 b_{\infty,y}=\tfrac34(3-\mathbb E T_-)>2,\qquad
 b_\infty\in2\rho+\Xi .
\end{equation}
Numerically, $a_\infty\simeq1.332168$ and
$b_{\infty,y}\simeq2.022231$.

In fact the strict inequality holds for every finite $k\ge1$.
Let $p_k$ and $p_\infty$ denote the normalized marginal densities of $t$
in \eqref{sec9:sl2-blowup-moments}. Truncating $p_\infty$ at $k$ removes
a positive first moment. On $[-1,k]$ the mean for density proportional to
$e^{-a t}R(t)^3$ decreases strictly with $a$ (its derivative is minus
the variance), so restoring mean zero gives $a_k<a_\infty$.
On $[-1,k]$ the ratio $p_k/p_\infty$ is therefore strictly increasing,
and beyond $k$ it is zero. Equality of the masses and first moments
forces a crossing $c\in(-1,k)$. Indeed, $p_k$ must exceed $p_\infty$
somewhere to compensate for the discarded tail. If $p_k\ge p_\infty$
throughout $[-1,k]$, then $\int(t-k)(p_k-p_\infty)\,dt<0$.
Thus $p_k-p_\infty$ has successive signs $-,+,-$ on the intervals
cut out by $c$ and $k$. For the convex function $g(t)=\max(-t,0)$,
let $\ell$ be the affine function agreeing with $g$ at $c$ and $k$.
Then $(g-\ell)(p_k-p_\infty)\le0$ everywhere and
$\int\ell(p_k-p_\infty)\,dt=0$. Hence
$\mathbb E_k g\le\mathbb E_\infty g$, which proves
$b_{k,y}\ge b_{\infty,y}>2$.

\begin{prop}\label{sec9:sl2-blowup-soliton}
For each $k\ge1$ the pair $(\bar M,D)$ admits a singular
$K\times K$-invariant KR soliton in $\Omega_k$.
The complement
\[
 M=\operatorname{Bl}_{Q^2\times\{0\}}(Q^3\times\mathbb C)
\]
admits a complete $K\times K$-invariant gradient shrinking
KR soliton with moment polyhedron
$P_\infty$ given by \eqref{sec9:sl2-blowup-polygons} and central soliton vector
$v_\infty=(0,-a_\infty)$. Along a sequence $k_j\to\infty$ the
finite-parameter solitons converge to it in $C^\infty_{\mathrm{loc}}(M)$
after pullback by central automorphisms.
\end{prop}
\begin{proof}
Here $K=SU(2)\times S^1$. Apply Theorem~\ref{existence-G} using
$b_{k,y}>2$ and Theorem~\ref{g-para} using
\eqref{sec9:sl2-blowup-polygons} and \eqref{sec9:sl2-blowup-strict}.
\end{proof}

\subsection{Three single-divisor deformations of $GL(2,\mathbb C)$-compactifications}
\label{sec9:gl2}

In this subsection, we discuss three  $\bar M=GL(2,\mathbb C)$-manifolds described in \cite{Zhu},  each of which moment polytope   corresponds  to a Fano cone after moving a smooth divisor of  $\bar M$.   Precisely, we move the infinity section in two projective bundles $Y_m$ and its
strict transform in a blowup $X_8$. In each case a single Weyl-fixed
facet recedes to infinity, leaving a two-dimensional recession cone.
The line bundles admit a common radial construction. For the blowup,
we use its proper map to matrix space to localize the metric estimates.

\subsubsection{The pairs \texorpdfstring{$(Y_m,I_m)$, $m=1,3$}{(Ym, Im), m = 1, 3}}
\label{sec9:gl2-y1}\label{sec9:gl2-y3}

Write $V=\operatorname{Mat}_{2\times2}(\mathbb C)$ and
$\mathbb P^3=\mathbb P(V)$.  For $m=1,3$, set
\[
 Y_m=\mathbb P_{\mathbb P^3}
       \bigl(\mathcal O\oplus\mathcal O(m)\bigr),\qquad
 \pi_m:Y_m\longrightarrow\mathbb P^3,
\]
and let $I_m$ be the section defined by the quotient onto
$\mathcal O(m)$.  If $\xi_m=c_1(\mathcal O_{Y_m}(1))$ and
$H=\pi_m^*c_1(\mathcal O_{\mathbb P^3}(1))$, then
\begin{equation}\label{sec9:gl2-ym-geometry}
 \begin{gathered}
 N_{I_m/Y_m}\simeq\mathcal O_{\mathbb P^3}(m),\qquad
 -K_{Y_m}=2\xi_m+(4-m)H,\\
 Y_m\setminus I_m\simeq
 \operatorname{Tot}\mathcal O_{\mathbb P^3}(-m).
 \end{gathered}
\end{equation}
The bundle $\mathcal O\oplus\mathcal O(m)$ is globally generated,
so $I_m$ is a smooth semiample divisor.  The displayed anticanonical
class is ample for $m=1,3$.

In the weight coordinates $y=(y_1,y_2)$ of $g$, put $s=y_1+y_2$
and $t=y_1-y_2$. The Weyl group exchanges $y_1$ and $y_2$,
the positive chamber is $t\ge0$, and $2\rho=(1,-1)$.
On the diagonal torus the orders of the two characters of $g$
along the two coordinate hyperplanes of $A$ give outward normals
$(d+1,d)$ and $(d,d+1)$. Their anticanonical support numbers are
both $2$ by \eqref{relation-lambda}. The two sections have outward
normals $\pm(1,1)$ and support number $1$ at $k=1$.
Moving $I_m$ therefore gives

\begin{equation}\label{sec9:gl2-ym-polytopes}
 P_{m,k}
 =\{-k\le s\le1,\quad |t|+ms\le4\},\qquad k\ge1.
\end{equation}

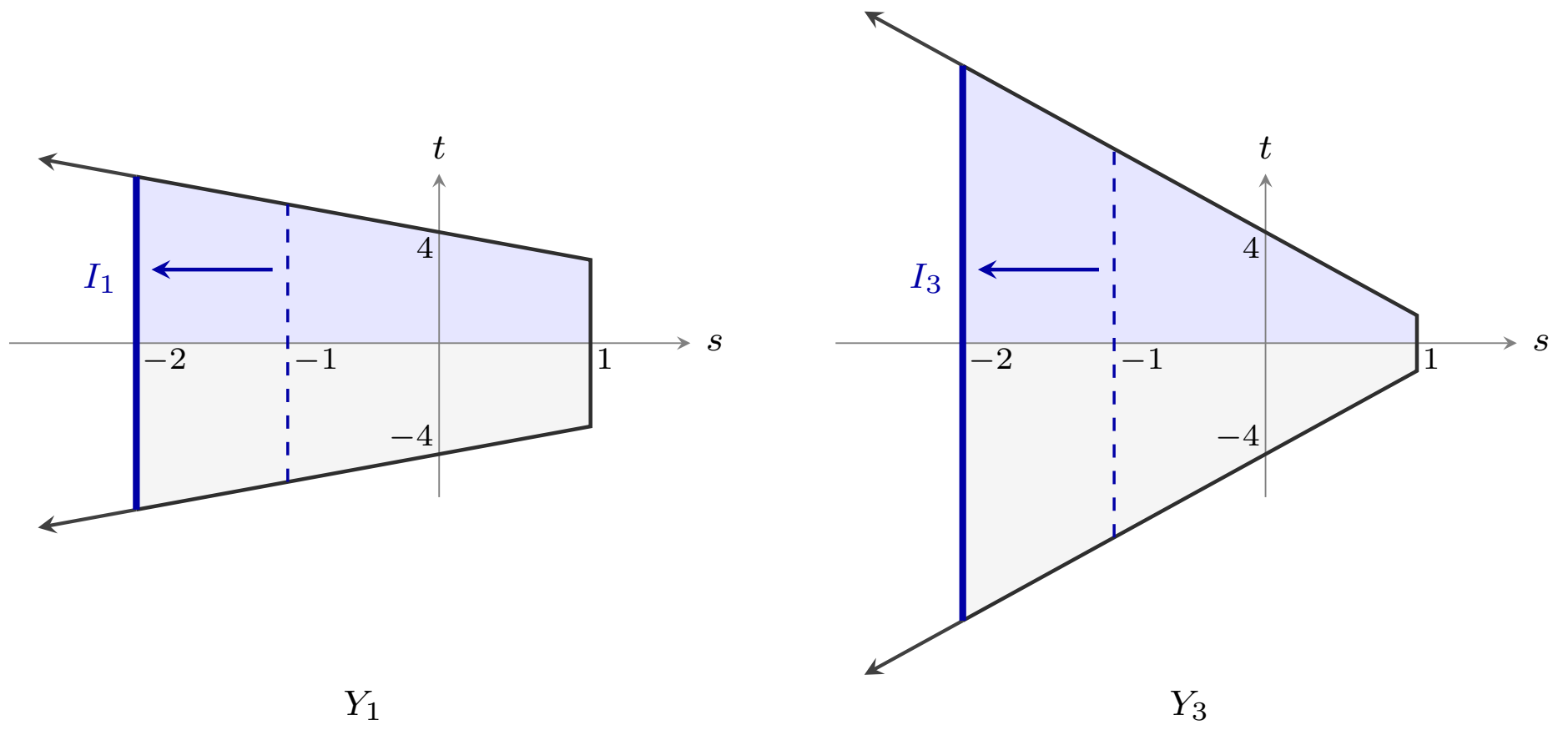
\begin{figure}[htbp]
\centering
\begin{tikzpicture}[
  x=1.2cm,y=.22cm,font=\footnotesize,
  axis/.style={->,>=stealth,draw=black!50,line width=.35pt},
  fixed/.style={draw=black!82,line width=.85pt},
  moving/.style={draw=blue!65!black,line width=1.55pt},
  former/.style={draw=blue!65!black,dashed,line width=.65pt},
  ray/.style={->,>=stealth,draw=black!75,line width=.8pt},
  ticktext/.style={font=\scriptsize,inner sep=1pt}]
  \foreach \panelshift/\ymindex in {0/1,6.55/3}{
    \begin{scope}[xshift=\panelshift cm]
      \pgfmathsetmacro{\topnew}{4+2*\ymindex}
      \pgfmathsetmacro{\topold}{4+\ymindex}
      \pgfmathsetmacro{\topright}{4-\ymindex}
      \pgfmathsetmacro{\topext}{4+2.65*\ymindex}
      \fill[black!4] (-2,0) -- (-2,-\topnew) -- (1,-\topright) -- (1,0) -- cycle;
      \fill[blue!10] (-2,0) -- (-2,\topnew) -- (1,\topright) -- (1,0) -- cycle;
      \draw[axis] (-2.84,0) -- (1.66,0) node[right,black] {$s$};
      \draw[axis] (0,-5.55) -- (0,6.10) node[above,black] {$t$};
      \draw[fixed] (-2,\topnew) -- (1,\topright)
        -- (1,-\topright) -- (-2,-\topnew);
      \draw[ray] (-2,\topnew) -- (-2.65,\topext);
      \draw[ray] (-2,-\topnew) -- (-2.65,-\topext);
      \draw[former] (-1,-\topold) -- (-1,\topold);
      \draw[moving] (-2,-\topnew) -- (-2,\topnew);
      \draw[->,>=stealth,blue!65!black,line width=.85pt] (-1.10,2.65) -- (-1.90,2.65);
      \node[anchor=east,text=blue!65!black] at (-2.02,2.35) {$I_{\ymindex}$};
      \foreach \sv in {-2,-1,1}{
        \node[ticktext,below right] at (\sv,0) {$\sv$};
      }
      \node[ticktext,below left] at (0,4) {$4$};
      \node[ticktext,above left] at (0,-4) {$-4$};
      \node at (-.50,-13.05) {$Y_{\ymindex}$};
    \end{scope}
  }
\end{tikzpicture}
\caption{The moment polygons of $Y_1$ and $Y_3$ at $k=2$, drawn at the
same scale. The blue shading marks the chamber $t\ge0$. The facet of
$I_m$ moves from $s=-1$ (dashed) to $s=-2$ (solid blue). The arrows
on the sloping boundaries $t=\pm(4-ms)$ show their extensions in
$P_{m,\infty}$. The right edge $s=1$ represents the negative section.}
\label{sec9:gl2-ym-figure}
\end{figure}

The Duistermaat--Heckman density is a positive constant multiple
of $t^2$.  The central soliton vector is therefore
$v_{m,k}=a_{m,k}(1,1)$, which is determined by the modified Futaki-invariant vanishing, i.e., 
\begin{equation}\label{sec9:gl2-ym-central-moment}
 \int_{-k}^{1}s\,e^{a_{m,k}s}(4-ms)^3\,ds=0.
\end{equation}
The left-hand side is strictly increasing in $a$ and positive for
large $a$. At $a=0$ the contributions at $s$ and $-s$ for $0<s<1$
have negative sum, and the interval $[-k,-1]$ contributes a
nonpositive term. Hence there is a unique solution $a_{m,k}>0$.  If $\mathbb E_{m,k}$ denotes
expectation for the probability density proportional to
$e^{a_{m,k}s}(4-ms)^3$ on $[-k,1]$, then the conditional mean of $t$ at fixed $s$ is
$\int_0^{4-ms}t^3\,dt/\int_0^{4-ms}t^2\,dt=3(4-ms)/4$. Hence
\begin{align}\label{sec9:gl2-ym-barycenter}
 \mathbb E_{m,k}(s)&=0,\notag\\
 b_{t,m,k}
 &=\frac34\mathbb E_{m,k}(4-ms)=3,\qquad
 \operatorname{bar}_{v_{m,k}}(P_{m,k}^+)
   =\left(\frac32,-\frac32\right)
   \in 2\rho+\Xi.
\end{align}
In particular, the finite-parameter barycenter condition holds
for every $k\ge1$, with a fixed strict margin.

The limiting polyhedron and its recession cone are
\begin{equation}\label{sec9:gl2-ym-limit}
 P_{m,\infty}=\{s\le1,\quad |t|+ms\le4\},\qquad
 C_{m,\infty}
   =\{(\delta s,\delta t):m\delta s+|\delta t|\le0\}.
\end{equation}
The cone is pointed and two-dimensional, and its central decay
domain is $a>0$.  Let $a_{m,\infty}$ be the
unique positive solution of
\begin{equation}\label{sec9:gl2-ym-limiting-moment}
 \int_{-\infty}^{1}s\,e^{a_{m,\infty}s}(4-ms)^3\,ds=0.
\end{equation}
Thus
$0<a_{m,1}\le a_{m,k}<a_{m,\infty}$, since increasing $k$ adds
only a negative contribution to \eqref{sec9:gl2-ym-central-moment}
at fixed $a$.  These bounds give uniform
weighted polynomial moments, and dominated convergence identifies
the increasing limit as $a_{m,\infty}$.  Together with
\eqref{sec9:gl2-ym-barycenter}, it also gives
$\operatorname{bar}_{v_{m,\infty}}(P_{m,\infty}^+)
=(3/2,-3/2)$. Following the conjecture (\ref{conj-g}), we make the following claim: 

\begin{claim}\label{sec9:gl2-ym-smooth-limit}
For $m=1,3$ and every integer $k\ge1$, one can choose a
$U(4)$-invariant singular KR soliton $\omega_{m,k}\in\Omega_k$ on
$Y_m$ with cone angle $2\pi k$ along $I_m$. After normalizing by
fiber dilations, the whole sequence converges as $k\to\infty$ in
$C^\infty_{\mathrm{loc}}(Y_m\setminus I_m)$ to the complete
FIK shrinker on
$\operatorname{Tot}\mathcal O_{\mathbb P^3}(-m)$.
\end{claim}

\subsubsection{The pair $(X_8,I_1')$}
\label{sec9:gl2-x8}

Let $E\simeq\mathbb P^3$ be the negative section of $Y_1$, and let
$Q^2=\{\det A=0\}\subset E$ be the smooth quadric of rank-one
matrices.  Set
\[
 \beta:X_8=\operatorname{Bl}_{Q^2}Y_1\longrightarrow Y_1,
 \qquad I_1'=\beta^{-1}(I_1).
\]
The center is disjoint from $I_1$, so $I_1'$ is smooth and
$\mathcal O(I_1')=\beta^*\mathcal O(I_1)$ is globally generated.
This is the smooth Fano compactification denoted by $X_8$.
Its complement has a proper morphism to matrix space:
\begin{equation}\label{sec9:gl2-x8-complement}
 X_8\setminus I_1'
 =\operatorname{Bl}_{Q^2}\bigl(\operatorname{Bl}_0\mathbb C^4\bigr)
 \longrightarrow\mathbb C^4.
\end{equation}
The moving anticanonical polytopes are
\begin{equation}\label{sec9:gl2-x8-polytopes}
 P_{8,k}=\{-k\le s\le1,\quad |t|+s\le4,
                                  \quad |t|+3s\le4\}.
\end{equation}
 $P_{8,k}$ is wonderful for every $k\ge1$, and the only
moving facet is the one corresponding to $I_1'$.

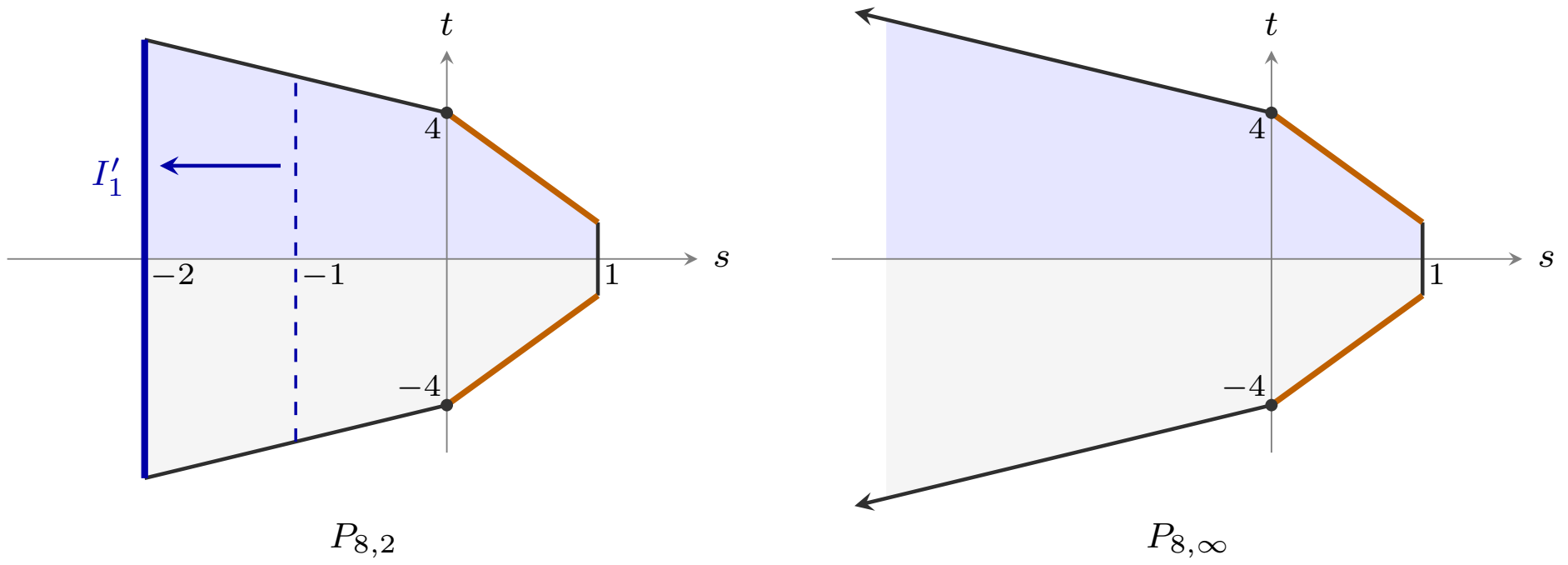
\begin{figure}[htbp]
\centering
\begin{tikzpicture}[
  x=1.2cm,y=.29cm,font=\footnotesize,
  axis/.style={->,>=stealth,draw=black!50,line width=.35pt},
  fixed/.style={draw=black!82,line width=.85pt},
  exceptional/.style={draw=orange!75!black,line width=1.35pt},
  moving/.style={draw=blue!65!black,line width=1.55pt},
  former/.style={draw=blue!65!black,dashed,line width=.65pt},
  ticktext/.style={font=\scriptsize,inner sep=1pt}]
  \foreach \panelshift/\islimit in {0/0,6.55/1}{
    \begin{scope}[xshift=\panelshift cm]
      \ifnum\islimit=0 \def\leftend{-2}\else \def\leftend{-2.55}\fi
      \pgfmathsetmacro{\topend}{4-\leftend}
      \fill[black!4] (\leftend,0) -- (\leftend,-\topend) -- (0,-4)
        -- (1,-1) -- (1,0) -- cycle;
      \fill[blue!10] (\leftend,0) -- (\leftend,\topend) -- (0,4)
        -- (1,1) -- (1,0) -- cycle;
      \draw[axis] (-2.91,0) -- (1.66,0) node[right,black] {$s$};
      \draw[axis] (0,-5.30) -- (0,5.70) node[above,black] {$t$};
      \draw[fixed] (1,-1) -- (1,1);
      \draw[exceptional] (0,4) -- (1,1) (0,-4) -- (1,-1);
      \ifnum\islimit=0
        \draw[fixed] (0,4) -- (-2,6) (0,-4) -- (-2,-6);
        \draw[former] (-1,-5) -- (-1,5);
        \draw[moving] (-2,-6) -- (-2,6);
        \draw[->,>=stealth,blue!65!black,line width=.85pt] (-1.10,2.55) -- (-1.90,2.55);
        \node[anchor=east,text=blue!65!black] at (-2.02,2.25) {$I_1'$};
        \foreach \sv in {-2,-1}{\node[ticktext,below right] at (\sv,0) {$\sv$};}
        \node at (-.55,-7.70) {$P_{8,2}$};
      \else
        \draw[fixed,->,>=stealth] (0,4) -- (-2.76,6.76);
        \draw[fixed,->,>=stealth] (0,-4) -- (-2.76,-6.76);
        \node at (-.55,-7.70) {$P_{8,\infty}$};
      \fi
      \fill[black!80] (0,4) circle[radius=1.4pt] (0,-4) circle[radius=1.4pt];
      \node[ticktext,below left] at (0,4) {$4$};
      \node[ticktext,above left] at (0,-4) {$-4$};
      \node[ticktext,below right] at (1,0) {$1$};
    \end{scope}
  }
\end{tikzpicture}
\caption{The deformation on $X_8$. The left panel shows the facet of
$I_1'$ at $s=-1$ (dashed) and $s=-2$ (solid blue). In the limit
the two sides $|t|+s=4$ extend to infinity. The two orange edges
$|t|+3s=4$ represent the exceptional divisor on the torus closure and
meet those sides at $(0,\pm4)$. Blue shading marks the chamber $t\ge0$.}
\label{sec9:gl2-x8-figure}
\label{sec9:gl2-divisor-figure}
\end{figure}

Figure~\ref{sec9:gl2-x8-figure} shows the cut induced by the
blowup. The fixed right edge represents the strict transform of the
negative section. The edges $|t|+s=4$ come from the strict transform
of $\pi_1^{-1}(\{\det A=0\})$, and $|t|+3s=4$ comes from the exceptional
divisor of $\beta$.

Put
\[
 h(s)=\begin{cases}4-s,&s\le0,\\4-3s,&0\le s\le1.\end{cases}
 \qquad
 I_k(a)=\int_{-k}^1s e^{as}h(s)^3\,ds.
\]
The central parameter $a_{8,k}$ is the unique zero of $I_k$.
Indeed, $I_k'(a)>0$, and direct integration gives
\[
 I_k(0)\le-\frac{228}{5}<0,\qquad
 I_k(3)>I_\infty(3)
       =\frac{1360}{81}-\frac{e^3}{9}>0.
\]
Hence $0<a_{8,k}<3$.  Write $\mathbb E_k$ for expectation with
density proportional to $e^{a_{8,k}s}h(s)^3$ on $[-k,1]$.
Since $\mathbb E_k s=0$ and $h(s)=4-s-2s_+$,
\begin{equation}\label{sec9:gl2-x8-root-moment}
 b_{t,8,k}=\frac34\mathbb E_k h(s)
            =3-\frac32\mathbb E_k s_+.
\end{equation}
To bound the last expectation, set
$f_a(s)=e^{as}(4-3s)^3$ on $[0,1]$.  For $0\le s<1/2$
and $a<3$,
\[
 \log\frac{f_a(s)}{f_a(1-s)}
 =3\log\frac{4-3s}{1+3s}-a(1-2s)
 \ge\left(\frac{18}{5}-a\right)(1-2s)>0.
\]
Pairing $s$ with $1-s$ shows that the conditional mean on
$[0,1]$ is less than $1/2$.  Consequently
$\mathbb E_k s_+<1/2$, and
\begin{equation}\label{sec9:gl2-x8-barycenter}
 b_{t,8,k}>\frac94>2,\qquad
 \operatorname{bar}_{v_{8,k}}(P_{8,k}^+)
 =\left(\frac{b_{t,8,k}}2,-\frac{b_{t,8,k}}2\right)
 \in2\rho+\Xi.
\end{equation}
This proves a strict margin independent of $k$.

As above, $a_{8,k}\uparrow a_{8,\infty}\in(0,3)$, where
$I_\infty(a_{8,\infty})=0$.  Exponential domination gives
convergence of every weighted polynomial moment, and the same
pairing argument gives $b_{t,8,\infty}>9/4$.
The limiting polyhedron is obtained from
\eqref{sec9:gl2-x8-polytopes} by omitting $s\ge-k$, and
\[
 C_{8,\infty}
 =\{(\delta s,\delta t):\delta s+|\delta t|\le0\}
 =\{\delta y_1\le0,\ \delta y_2\le0\}.
\]
It is again pointed and two-dimensional. Following the conjecture (\ref{conj-g}), we make the following claim: 
\begin{claim}\label{sec9:gl2-x8-smooth-limit}
For every integer $k\ge1$ there is a finite-parameter KR soliton
$\omega_{8,k}\in\Omega_k$ on $(X_8,I_1')$. There are integers
$k_j\to\infty$ and central automorphisms $\sigma_j$ such that
\[
 \sigma_j^*\omega_{8,k_j}\longrightarrow\omega_{8,\infty}
 \quad\text{in }C^\infty_{\mathrm{loc}}(X_8\setminus I_1').
\]
The limit is a complete $U(2)\times U(2)$-invariant gradient shrinking
KR soliton with torus-closure moment image $P_{8,\infty}$ and central
soliton vector $v_{8,\infty}=a_{8,\infty}(1,1)$.
\end{claim}


\begin{thebibliography}{99}
\let\l\textPolishL
\urlstyle{same}

\bibitem{Ab-JDG}
M. Abreu, \textit{K\"ahler metrics on toric orbifolds},
J. Differential Geom. \textbf{58} (2001), no.~1, 151--187.
\href{https://doi.org/10.4310/jdg/1090348285}{\nolinkurl{doi:10.4310/jdg/1090348285}}.

\bibitem{AB1}
V. Alexeev and M. Brion, \textit{Stable reductive varieties I: Affine varieties},
Invent. Math. \textbf{157} (2004), no.~2, 227--274.
\href{https://doi.org/10.1007/s00222-003-0347-y}{\nolinkurl{doi:10.1007/s00222-003-0347-y}}.

\bibitem{AB2}
V. Alexeev and M. Brion, \textit{Stable reductive varieties II: Projective case},
Adv. Math. \textbf{184} (2004), no.~2, 380--408.
\href{https://doi.org/10.1016/S0001-8708(03)00164-6}{\nolinkurl{doi:10.1016/S0001-8708(03)00164-6}}.

\bibitem{AK}
V. Alexeev and L. Katzarkov, \textit{On K-stability of reductive varieties},
Geom. Funct. Anal. \textbf{15} (2005), no.~2, 297--310.
\href{https://doi.org/10.1007/s00039-005-0507-x}{\nolinkurl{doi:10.1007/s00039-005-0507-x}}.

\bibitem{AL}
H. Azad and J.-J. Loeb, \textit{Plurisubharmonic functions and K\"ahlerian metrics on complexification of symmetric spaces},
Indag. Math. (N.S.) \textbf{3} (1992), no.~4, 385--390.
\href{https://doi.org/10.1016/0019-3577(92)90017-F}{\nolinkurl{doi:10.1016/0019-3577(92)90017-F}}.

\bibitem{BCCD}
R. H. Bamler, C. Cifarelli, R. J. Conlon and A. Deruelle, \textit{A new complete two-dimensional shrinking gradient K\"ahler--Ricci soliton},
Geom. Funct. Anal. \textbf{34} (2024), no.~2, 377--392.
\href{https://doi.org/10.1007/s00039-024-00668-9}{\nolinkurl{doi:10.1007/s00039-024-00668-9}}.

\bibitem{BB}
R. J. Berman and B. Berndtsson, \textit{Real Monge--Amp\`ere equations and K\"ahler--Ricci solitons on toric log Fano varieties},
Ann. Fac. Sci. Toulouse Math. (6) \textbf{22} (2013), no.~4, 649--711.
\href{https://doi.org/10.5802/afst.1386}{\nolinkurl{doi:10.5802/afst.1386}}.

\bibitem{BBEGZ}
R. J. Berman, S. Boucksom, P. Eyssidieux, V. Guedj and A. Zeriahi, \textit{K\"ahler--Einstein metrics and the K\"ahler--Ricci flow on log Fano varieties},
J. Reine Angew. Math. \textbf{751} (2019), 27--89.
\href{https://doi.org/10.1515/crelle-2016-0033}{\nolinkurl{doi:10.1515/crelle-2016-0033}}.

\bibitem{BLXZ}
H. Blum, Y. Liu, C. Xu and Z. Zhuang, \textit{The existence of the K\"ahler--Ricci soliton degeneration},
Forum Math. Pi \textbf{11} (2023), article e9, 28~pp.
\href{https://doi.org/10.1017/fmp.2023.5}{\nolinkurl{doi:10.1017/fmp.2023.5}}.

\bibitem{CLN}
B. Chow, S.-C. Chu, D. Glickenstein, C. Guenther, J. Isenberg, T. Ivey, D. Knopf, P. Lu, F. Luo and L. Ni, \textit{The Ricci Flow: Techniques and Applications. Part III: Geometric-Analytic Aspects},
Math. Surveys Monogr., vol.~163, American Mathematical Society, Providence, RI, 2010.
\href{https://doi.org/10.1090/surv/163}{\nolinkurl{doi:10.1090/surv/163}}.

\bibitem{Ci}
C. Cifarelli, \textit{Uniqueness of shrinking gradient K\"ahler--Ricci solitons on non-compact toric manifolds},
J. Lond. Math. Soc. (2) \textbf{106} (2022), no.~4, 3746--3791.
\href{https://doi.org/10.1112/jlms.12673}{\nolinkurl{doi:10.1112/jlms.12673}}.

\bibitem{CCD-type1}
C. Cifarelli, R. J. Conlon and A. Deruelle, \textit{On finite time Type I singularities of the K\"ahler--Ricci flow on compact K\"ahler surfaces},
J. Eur. Math. Soc. \textbf{28} (2026), no.~2, 463--504.
\href{https://doi.org/10.4171/JEMS/1485}{\nolinkurl{doi:10.4171/JEMS/1485}}.

\bibitem{CHTX}
T. C. Collins, G. Huang, F. Tong and Y. Xu, \textit{Shrinking K\"ahler--Ricci solitons on toric Fano fibrations},
preprint (2026), \href{https://arxiv.org/abs/2609.24884}{arXiv:2609.24884}.

\bibitem{CGSZ}
D. Coman, V. Guedj, S. Sahin and A. Zeriahi, \textit{Toric pluripotential theory},
Ann. Polon. Math. \textbf{123} (2019), no.~1, 215--242.
\href{https://doi.org/10.4064/ap180409-3-7}{\nolinkurl{doi:10.4064/ap180409-3-7}}.

\bibitem{CCD-classification}
R. J. Conlon, A. Deruelle and S. Sun, \textit{Classification results for expanding and shrinking gradient K\"ahler--Ricci solitons},
Geom. Topol. \textbf{28} (2024), no.~1, 267--351.
\href{https://doi.org/10.2140/gt.2024.28.267}{\nolinkurl{doi:10.2140/gt.2024.28.267}}.

\bibitem{CK}
D. Cordero-Erausquin and B. Klartag, \textit{Moment measures},
J. Funct. Anal. \textbf{268} (2015), no.~12, 3834--3866.
\href{https://doi.org/10.1016/j.jfa.2015.04.001}{\nolinkurl{doi:10.1016/j.jfa.2015.04.001}}.

\bibitem{DGSW}
V. Datar, B. Guo, J. Song and X. Wang, \textit{Connecting toric manifolds by conical K\"ahler--Einstein metrics},
Adv. Math. \textbf{323} (2018), 38--83.
\href{https://doi.org/10.1016/j.aim.2017.10.035}{\nolinkurl{doi:10.1016/j.aim.2017.10.035}}.

\bibitem{Del1}
T. Delcroix, \textit{K\"ahler--Einstein metrics on group compactifications},
Geom. Funct. Anal. \textbf{27} (2017), no.~1, 78--129.
\href{https://doi.org/10.1007/s00039-017-0394-y}{\nolinkurl{doi:10.1007/s00039-017-0394-y}}.

\bibitem{Del}
T. Delcroix, \textit{$K$-stability of Fano spherical varieties},
Ann. Sci. \'Ec. Norm. Sup\'er. (4) \textbf{53} (2020), no.~3, 615--662.
\href{https://doi.org/10.24033/asens.2430}{\nolinkurl{doi:10.24033/asens.2430}}.


\bibitem{DM}
T. Delcroix and P.-L. Montagard, \textit{Spherical actions on locally factorial Fano varieties of dimension $\leq 4$ and rank $\leq 2$},
extended version (2023), \href{https://arxiv.org/abs/2308.15858}{arXiv:2308.15858}. Shortened version accepted for publication in \textit{Algebraic Geometry}.

\bibitem{DeS20}
R. Dervan and G. Sz\'ekelyhidi, \textit{The K\"ahler--Ricci flow and optimal degenerations},
J. Differential Geom. \textbf{116} (2020), no.~1, 187--203.
\href{https://doi.org/10.4310/jdg/1599271255}{\nolinkurl{doi:10.4310/jdg/1599271255}}.

\bibitem{Do02}
S.~K. Donaldson, \textit{Scalar curvature and stability of toric varieties},
J. Differential Geom. \textbf{62} (2002), no.~2, 289--349.
\href{https://doi.org/10.4310/jdg/1090950195}{\nolinkurl{doi:10.4310/jdg/1090950195}}.

\bibitem{EMT}
J. Enders, R. M\"uller, and P.~M. Topping, \textit{On type-I singularities in Ricci flow},
Comm. Anal. Geom. \textbf{19} (2011), no.~5, 905--922.
\href{https://doi.org/10.4310/CAG.2011.v19.n5.a4}{\nolinkurl{doi:10.4310/CAG.2011.v19.n5.a4}}.

\bibitem{FIK}
M. Feldman, T. Ilmanen, and D. Knopf, \textit{Rotationally symmetric shrinking and expanding gradient K\"ahler--Ricci solitons},
J. Differential Geom. \textbf{65} (2003), no.~2, 169--209.
\href{https://doi.org/10.4310/jdg/1090511686}{\nolinkurl{doi:10.4310/jdg/1090511686}}.

\bibitem{Fig}
A. Figalli, \textit{The Monge--Amp\`ere Equation and Its Applications},
Zurich Lectures in Advanced Mathematics, European Mathematical Society, Z\"urich, 2017.
\href{https://doi.org/10.4171/170}{\nolinkurl{doi:10.4171/170}}.

\bibitem{FW}
A. Futaki, H. Ono, and G. Wang, \textit{Transverse K\"ahler geometry of Sasaki manifolds and toric Sasaki--Einstein manifolds},
J. Differential Geom. \textbf{83} (2009), no.~3, 585--635.
\href{https://doi.org/10.4310/jdg/1264601036}{\nolinkurl{doi:10.4310/jdg/1264601036}}.

\bibitem{GT}
D. Gilbarg and N.~S. Trudinger, \textit{Elliptic Partial Differential Equations of Second Order},
2nd ed., Grundlehren der mathematischen Wissenschaften, vol.~224, Springer-Verlag, Berlin, 1983, ISBN 978-3-540-13025-3.

\bibitem{Gu}
V. Guillemin, \textit{Moment Maps and Combinatorial Invariants of Hamiltonian $T^n$-spaces},
Progress in Mathematics, vol.~122, Birkh\"auser, Boston, MA, 1994.
\href{https://doi.org/10.1007/978-1-4612-0269-1}{\nolinkurl{doi:10.1007/978-1-4612-0269-1}}.

\bibitem{Ha82}
R.~S. Hamilton, \textit{Three-manifolds with positive Ricci curvature},
J. Differential Geom. \textbf{17} (1982), no.~2, 255--306.
\href{https://doi.org/10.4310/jdg/1214436922}{\nolinkurl{doi:10.4310/jdg/1214436922}}.

\bibitem{Ha95}
R.~S. Hamilton, \textit{The formation of singularities in the Ricci flow},
in \textit{Surveys in Differential Geometry, Vol.~II} (Cambridge, MA, 1993), International Press, Cambridge, MA, 1995, 7--136.
\href{https://doi.org/10.4310/SDG.1993.v2.n1.a2}{\nolinkurl{doi:10.4310/SDG.1993.v2.n1.a2}}.

\bibitem{HL}
J. Han and C. Li, \textit{On the Yau--Tian--Donaldson conjecture for generalized K\"ahler--Ricci soliton equations},
Comm. Pure Appl. Math. \textbf{76} (2023), no.~9, 1793--1867.
\href{https://doi.org/10.1002/cpa.22053}{\nolinkurl{doi:10.1002/cpa.22053}}.

\bibitem{Kol98}
S. Ko{\l}odziej, \textit{The complex Monge--Amp\`ere equation},
Acta Math. \textbf{180} (1998), no.~1, 69--117.
\href{https://doi.org/10.1007/BF02392879}{\nolinkurl{doi:10.1007/BF02392879}}.

\bibitem{KW-rigidity}
B. Kotschwar and L. Wang, \textit{Rigidity of asymptotically conical shrinking gradient Ricci solitons},
J. Differential Geom. \textbf{100} (2015), no.~1, 55--108.
\href{https://doi.org/10.4310/jdg/1427202764}{\nolinkurl{doi:10.4310/jdg/1427202764}}.

\bibitem{Legendre}
E. Legendre, \textit{Toric K\"ahler--Einstein metrics and convex compact polytopes},
J. Geom. Anal. \textbf{26} (2016), no.~1, 399--427.
\href{https://doi.org/10.1007/s12220-015-9556-z}{\nolinkurl{doi:10.1007/s12220-015-9556-z}}.

\bibitem{LTZ}
Y. Li, G. Tian, and X. Zhu, \textit{Singular K\"ahler--Einstein metrics on $\mathbb{Q}$-Fano compactifications of Lie groups},
Math. Eng. \textbf{5} (2023), no.~2, paper 028, 1--43.
\href{https://doi.org/10.3934/mine.2023028}{\nolinkurl{doi:10.3934/mine.2023028}}.

\bibitem{LZ}
Y. Li and B. Zhou, \textit{Mabuchi metrics and properness of the modified Ding functional},
Pacific J. Math. \textbf{302} (2019), no.~2, 659--692.
\href{https://doi.org/10.2140/pjm.2019.302.659}{\nolinkurl{doi:10.2140/pjm.2019.302.659}}.

\bibitem{LZZ}
Y. Li, B. Zhou, and X. Zhu, \textit{$K$-energy on polarized compactifications of Lie groups},
J. Funct. Anal. \textbf{275} (2018), no.~5, 1023--1072.
\href{https://doi.org/10.1016/j.jfa.2018.04.009}{\nolinkurl{doi:10.1016/j.jfa.2018.04.009}}.

\bibitem{LZhu}
Y. Li and X. Zhu, \textit{$G$-Sasaki manifolds and $K$-energy},
J. Geom. Anal. \textbf{31} (2021), no.~2, 1415--1470.
\href{https://doi.org/10.1007/s12220-019-00285-1}{\nolinkurl{doi:10.1007/s12220-019-00285-1}}.

\bibitem{LW}
Y. Li and B. Wang, \textit{On K\"ahler Ricci shrinker surfaces},
Acta Math. \textbf{236} (2026), no.~1, 1--50.
\href{https://doi.org/10.4310/ACTA.2026.v236.n1.a1}{\nolinkurl{doi:10.4310/ACTA.2026.v236.n1.a1}}.

\bibitem{MW-positive}
O. Munteanu and J. Wang, \textit{Positively curved shrinking Ricci solitons are compact},
J. Differential Geom. \textbf{106} (2017), no.~3, 499--505.
\href{https://doi.org/10.4310/jdg/1500084024}{\nolinkurl{doi:10.4310/jdg/1500084024}}.

\bibitem{Ni}
L. Ni, \textit{Closed type I ancient solutions to Ricci flow}, in \textit{Recent Advances in Geometric Analysis} (Y.-I. Lee, C.-S. Lin and M.-P. Tsui, eds.), Adv. Lect. Math. (ALM), vol.~11, International Press, Somerville, MA, 2010, 147--150.

\bibitem{Pe02}
G. Perelman, \textit{The entropy formula for the Ricci flow and its geometric applications}, preprint (2002), \href{https://arxiv.org/abs/math/0211159}{arXiv:math/0211159}.

\bibitem{SZ}
Y. Shi and X. Zhu, \textit{K\"ahler--Ricci solitons on toric Fano orbifolds}, Math. Z. \textbf{271} (2012), nos.~3--4, 1241--1251. \href{https://doi.org/10.1007/s00209-011-0913-8}{\nolinkurl{doi:10.1007/s00209-011-0913-8}}.

\bibitem{Siu}
Y.-T. Siu, \textit{Lectures on Hermitian--Einstein Metrics for Stable Bundles and K\"ahler--Einstein Metrics}, DMV Seminar, vol.~8, Birkh\"auser, Basel, 1987. \href{https://doi.org/10.1007/978-3-0348-7486-1}{\nolinkurl{doi:10.1007/978-3-0348-7486-1}}.

\bibitem{Ti97}
G. Tian, \textit{K\"ahler--Einstein metrics with positive scalar curvature}, Invent. Math. \textbf{130} (1997), no.~1, 1--37. \href{https://doi.org/10.1007/s002220050176}{\nolinkurl{doi:10.1007/s002220050176}}.

\bibitem{Ti15}
G. Tian, \textit{$K$-stability and K\"ahler--Einstein metrics}, Comm. Pure Appl. Math. \textbf{68} (2015), no.~7, 1085--1156. \href{https://doi.org/10.1002/cpa.21578}{\nolinkurl{doi:10.1002/cpa.21578}}.

\bibitem{TY1}
G. Tian and S.-T. Yau, \textit{Complete K\"ahler manifolds with zero Ricci curvature.~I}, J. Amer. Math. Soc. \textbf{3} (1990), no.~3, 579--609. \href{https://doi.org/10.1090/S0894-0347-1990-1040196-6}{\nolinkurl{doi:10.1090/S0894-0347-1990-1040196-6}}.

\bibitem{TY2}
G. Tian and S.-T. Yau, \textit{Complete K\"ahler manifolds with zero Ricci curvature.~II}, Invent. Math. \textbf{106} (1991), no.~1, 27--60. \href{https://doi.org/10.1007/BF01243902}{\nolinkurl{doi:10.1007/BF01243902}}.


 \bibitem{TZZZZ} 
 G. Tian, Q. S. Zhang, Z. Zhang, M. Zhu and X. Zhu,   \textit{Laplace comparison on Kaehler Ricci flow and convergence},  arXiv: 2509.14820; Peking Math. J., 2026, in press.

\bibitem{TZZZ}
G. Tian, S. Zhang, Z. Zhang and X. Zhu, \textit{Perelman's entropy and K\"ahler--Ricci flow on a Fano manifold}, Trans. Amer. Math. Soc. \textbf{365} (2013), no.~12, 6669--6695. \href{https://doi.org/10.1090/S0002-9947-2013-06027-8}{\nolinkurl{doi:10.1090/S0002-9947-2013-06027-8}}.

 \bibitem{TZh} G. Tian and  Z.L. Zhang, 
    \textit{Regularity of K\"ahler-Ricci flows on Fano manifolds},  Acta Math., \textbf{216} (2016), 127-176.

   

\bibitem{TZ1}
G. Tian and X. Zhu, \textit{Uniqueness of K\"ahler--Ricci solitons}, Acta Math. \textbf{184} (2000), no.~2, 271--305. \href{https://doi.org/10.1007/BF02392630}{\nolinkurl{doi:10.1007/BF02392630}}.

\bibitem{TZ02}
G. Tian and X. Zhu, \textit{A new holomorphic invariant and uniqueness of K\"ahler--Ricci solitons}, Comment. Math. Helv. \textbf{77} (2002), no.~2, 297--325. \href{https://doi.org/10.1007/s00014-002-8341-3}{\nolinkurl{doi:10.1007/s00014-002-8341-3}}.

\bibitem{TZ-applications}
G. Tian and X. Zhu, \textit{$H$-invariant and its applications}, Acta Math. Sin. (Engl. Ser.) \textbf{42} (2026), no.~3, 881--912. \href{https://doi.org/10.1007/s10114-026-5435-8}{\nolinkurl{doi:10.1007/s10114-026-5435-8}}.

\bibitem{TZ-survey}
G. Tian and X. Zhu, \textit{$H$-invariant on K\"ahler manifolds}, Peking Math. J. \textbf{9} (2026), no.~3, 565--575. \href{https://doi.org/10.1007/s42543-025-00108-z}{\nolinkurl{doi:10.1007/s42543-025-00108-z}}.

\bibitem{TZ24}
G. Tian and X. Zhu, \textit{Horosymmetric limits of K\"ahler--Ricci flow on Fano $\mathbf{G}$-manifolds}, J. Eur. Math. Soc. \textbf{28} (2026), no.~10, 4519--4558. \href{https://doi.org/10.4171/JEMS/1553}{\nolinkurl{doi:10.4171/JEMS/1553}}.

\bibitem{WT}
N. S. Trudinger and X.-J. Wang, \textit{The Monge--Amp\`ere equation and its geometric applications}, in \textit{Handbook of Geometric Analysis}, vol.~I (L. Ji, P. Li, R. Schoen and L. Simon, eds.), Adv. Lect. Math. (ALM), vol.~7, International Press, Somerville, MA, 2008, 467--524.

\bibitem{WZZ}
F. Wang, B. Zhou and X. Zhu, \textit{Modified Futaki invariant and equivariant Riemann--Roch formula}, Adv. Math. \textbf{289} (2016), 1205--1235. \href{https://doi.org/10.1016/j.aim.2015.11.036}{\nolinkurl{doi:10.1016/j.aim.2015.11.036}}.

\bibitem{WZ20}
F. Wang and X. Zhu, \textit{Tian's partial $C^0$-estimate implies Hamilton--Tian's conjecture}, Adv. Math. \textbf{381} (2021), article~107619. \href{https://doi.org/10.1016/j.aim.2021.107619}{\nolinkurl{doi:10.1016/j.aim.2021.107619}}.

\bibitem{WZ04}
X.-J. Wang and X. Zhu, \textit{K\"ahler--Ricci solitons on toric manifolds with positive first Chern class}, Adv. Math. \textbf{188} (2004), no.~1, 87--103. \href{https://doi.org/10.1016/j.aim.2003.09.009}{\nolinkurl{doi:10.1016/j.aim.2003.09.009}}.

\bibitem{Yau}
S.-T. Yau, \textit{On the Ricci curvature of a compact K\"ahler manifold and the complex Monge--Amp\`ere equation.~I}, Comm. Pure Appl. Math. \textbf{31} (1978), no.~3, 339--411. \href{https://doi.org/10.1002/cpa.3160310304}{\nolinkurl{doi:10.1002/cpa.3160310304}}.

\bibitem{Zhu}
X. Zhu, \textit{K\"ahler--Einsteinmetrics on toric manifolds and $\mathbf G$-manifolds}, Progress in Mathematics, Birkh\"auser, vol 333 (2020), 545--585. 

\end{thebibliography}
\end{document}